\documentclass[11pt,letterpaper]{amsart}  

\usepackage{amssymb,latexsym,bbm,graphicx,epsfig,epic,eepic,oldgerm,psfrag,comment}
\usepackage{amsmath,amsfonts,amscd,mathrsfs,amsthm,nicefrac}
\usepackage[outdir=./]{epstopdf}
\usepackage{mathtools}
\usepackage{mathtools}
\usepackage{tikz}
\usepackage{xcolor}
\usepackage{soul,xcolor}

\usepackage{times}
\usepackage{tikz,wasysym}
\usepackage{marvosym}                  
\usepackage{hyperref} 

\usepackage{xypic} 
\input xy
\xyoption{all}

\theoremstyle{plain}
\newtheorem{thm}{Theorem}
\newtheorem*{thm*}{Theorem}
\newtheorem{theorem}{Theorem}[section]
\newtheorem{lemma}[theorem]{Lemma}                           
\newtheorem{proposition}[theorem]{Proposition}

\newtheorem*{remark*}{Remark}
\newtheorem*{remarks*}{Remarks}
\newtheorem{remark}[theorem]{Remark}
\newtheorem{remarks}[theorem]{Remarks}

\newtheorem{example}[theorem]{Example}
\newtheorem{examples}[theorem]{Examples}
\newtheorem*{example*}{Example}
\newtheorem*{examples*}{Examples}
\newtheorem{definition}[theorem]{Definition}
\newtheorem*{definition*}{Definition}

\newtheorem{question}[theorem]{Question}
\newtheorem*{question*}{Question}

\newtheorem*{claim*}{Claim}

\numberwithin{figure}{section}
\numberwithin{equation}{section}

\newcommand{\proofend}{\hspace*{\fill} $\square$\\}
\newcommand{\diam}{\hspace*{\fill} $\diamond$}

\def\SDB{\operatorname{SDB}}

\def\Op{\operatorname{Op}}

\def\G{ {\operatorname{G}} }

\def\FS{{\operatorname{FS}}}

\newcommand{\cdcirc}{\overset%
{\raisebox{-.3ex}[0ex][-.3ex]{\mbox{$\scriptscriptstyle\circ$}} \mskip-2mu}\cd}

\newcommand{\DDcirc}{\overset%
{\raisebox{-.3ex}[0ex][-.3ex]{\mbox{$\scriptscriptstyle\circ$}} \mskip+2mu}\DD}

\newcommand{\Sigmapcirc}{\overset%
{\raisebox{-.3ex}[0ex][-.3ex]{\mbox{$\scriptscriptstyle\circ$}} \mskip6mu}{\Sigma'}}

\newcommand{\Xcirc}{\overset%
{\raisebox{-.2ex}[0ex][-.3ex]{\mbox{$\scriptscriptstyle\circ$}} \mskip-3mu}X}

\def\1{\:\!}
\def\2{\;\!}
\def\s{\smallskip}
\def\m{\medskip}
\def\eps{\varepsilon}

\def\im{\operatorname {im}}

\def\Diffc0{\operatorname{Diff^c_0}}

\def\Sympc0{\operatorname{Symp^c_0}}

\def\Int{\operatorname{Int}}

\def\supp{\operatorname{supp}}
\def\idd{\operatorname{id}}

\def\area{\operatorname{area}}
\def\sing{\operatorname{Sing}}

\def\PD{\operatorname{PD}}
\def\H{\operatorname{H}}

\def\const{\operatorname{const}}

\def\U{\operatorname{U}}

\def\ga{\alpha}

\def\gd{\delta}  

\def\gve{\varepsilon}
\def\gf{\varphi}

\def\gl{\lambda}
\def\go{\omega}

\def\gs{\sigma}

\def\T{\operatorname{T}}
\def\Z{\operatorname{Z}}

\def\ca{{\mathcal A}}
\def\cb{{\mathcal B}}

\def\cd{{\mathcal D}}

\def\cf{{\mathcal F}}

\def\cl{{\mathcal L}}

\def\co{{\mathcal O}}
\def\cp{{\mathcal P}}
\def\cq{{\mathcal Q}}

\def\cs{{\mathcal S}}

\def\CC{\mathbb{C}}
\def\DD{\mathbb{D}}

\def\NN{\mathbb{N}}

\def\RR{\mathbb{R}}

\def\ZZ{\mathbb{Z}}

\def\C{\mathbb{C}}

\def\H{\mathbb{H}}

\def\R{\mathbb{R}}
\def\T{\mathbb{T}}
\def\Z{\mathbb{Z}}

\def\RP{\operatorname{\mathbb{R}P}}
\def\CP{\operatorname{\mathbb{C}P}}

\def\st{{\operatorname{st}}}

\def\pp{\partial}

\def\fm{{\mathfrak m}}

\def\ni{\noindent}
\def\b{\bigskip}
\def\m{\medskip}

\def\id{\mbox{id}}

\def\proof{\noindent {\it Proof. \;}}

\newcommand{\nbd}{neighbourhood }
\newcommand{\fonction}[5]
{$$ 
\begin{array}{rcccl}
 #1 & : & #2 & \longrightarrow &#3 \\
    &   & #4 & \longmapsto &#5 
\end{array}
$$}

\newcommand{\dom}{\textnormal{Dom}\,}
\newcommand{\res}{\textnormal{Res}\,}

\newcommand{\priv}{\backslash}
\newcommand{\lra}{\longrightarrow}
\newcommand{\hra}{\hookrightarrow}
\newcommand{\hraa}{\lhook \joinrel \xrightarrow{\2\ga_{\st}\,}}

\newcommand{\om}{\omega}
\newcommand{\lag}{\textnormal{Lag}}

\newcommand{\sdb}{\textnormal{\small SDB}}

\newcommand{\what}[1]{\widehat{#1}}

\newcommand{\sbull}{{\tiny $\bullet$ }}

\newcommand{\Om}{\Omega}
\newcommand{\its}{\item[\sbull]}

\newcommand{\op}{\textnormal{Op}\,}
\newcommand{\skel}{\textnormal{Skel}}

\definecolor{marron}{rgb}{0.64,0.16,0.16}

\newcommand{\bfsigma}{{\bf \Sigma}}

\makeatletter
\renewcommand*\env@matrix[1][\arraystretch]{%
  \edef\arraystretch{#1}%
  \hskip -\arraycolsep
  \let\@ifnextchar\new@ifnextchar
  \array{*\c@MaxMatrixCols c}}
\makeatother

\def\bul{\raisebox{.2ex}{\sbull}}

\definecolor{amber}{rgb}{1.0, 0.75, 0.0}
\definecolor{banan}{rgb}{1.0, 0.88, 0.21}
\definecolor{gold}{rgb}{1.0, 0.84, 0.0}                
\definecolor{goldenpoppy}{rgb}{0.99, 0.76, 0.0}
\definecolor{forestgreen}{rgb}{0.13, 0.55, 0.13}
\definecolor{green}{rgb}{0.4, 0.69, 0.2}

\def\red#1{{\textcolor{red}{#1}}}
\def\blue#1{{\textcolor{blue}{#1}}}

\title[Liouville polarisations]{Liouville polarisations and the rigidity of their Lagrangian skeleta in dimension~$4$}

\author{Emmanuel Opshtein}
\address{(E.\ Opshtein) 
IRMA, 7 rue Ren\'e Decartes, 
67084 Strasbourg,
France}
\email{opshtein@unistra.fr}

\author{Felix Schlenk}
\address{(F.~Schlenk)
Institut de Math\'ematiques,
Universit\'e de Neuch\^atel,
Rue \'Emile Argand~11,
2000 Neuch\^atel,
Switzerland}
\email{schlenk@unine.ch}

\date{\today}

\begin{document}

\maketitle

\begin{abstract}
In 2000, Biran introduced polarisations of closed symplectic manifolds 
and showed that their Lagrangian skeleta exhibit remarkable rigidity properties. 
He found, in particular, that their complements contain only small balls. 
In this paper, we introduce so-called Liouville polarisations of certain open 4-dimensional symplectic manifolds.
This leads to several symplectic embedding results,
that in turn lead to new Lagrangian non-removable intersections
and a novel phenomenon of Legendrian barriers.

We show, for instance, that given any connected symplectic 4-manifold $(M,\go)$
and a 4-ball of smaller volume, there exists an explicit finite union of Lagrangian 
discs in the 4-ball such that their complement symplectically embeds into~$(M,\go)$,
extending a result by Sackel--Song--Varolgunes--Zhu and Brendel. 
Other applications are new Lagrangian intersection results and relative versions of the Arnold chord conjecture.
\end{abstract}

\bigskip
\noindent\textbf{Titre.}
Polarisations de Liouville et rigidité de leurs squelettes lagrangiens en dimension~4

\medskip
\noindent\textbf{Résumé.}
En 2000, Biran a introduit les polarisations de variétés symplectiques compactes et a montré que leurs squelettes lagrangiens présentent des propriétés remarquables de rigidité. 
Il a notamment constaté que leurs complémentaires ne contiennent que des petites boules. 
Dans cet article, nous introduisons des polarisations dites de Liouville de certaines variétés symplectiques non compactes de dimension~4.
Cela conduit à plusieurs résultats de plongements symplectiques, 
qui à leur tour conduisent à de nouvelles intersections lagrangiennes persistantes et 
un nouveau phénomène de barrières legendriennes.

Nous montrons par exemple qu'étant donné une variété symplectique connexe~$(M,\omega)$ de dimension~4
et une boule de dimension~4 de plus petit volume, il existe une union finie explicite de disques lagrangiens dans la 4-boule 
dont le complémentaire se plonge symplectiquement dans~$(M,\omega)$, 
ce qui étend un résultat de Sackel--Song--Varolgunes--Zhu et Brendel.
D'autres applications sont de nouveaux résultats d'inter\-sections lagrangiennes et des versions relatives de la conjecture des cordes d'Arnold.

%\newpage

\setcounter{tocdepth}{1}
\tableofcontents

\section{Introduction and main results}

A polarisation of a {\it closed}\/ symplectic manifold $(M,\omega)$ is a decomposition 
of the symplectic class~$[\omega]$ along positive real multiples of the Poincar\'e-dual classes 
of codimension~2 symplectic submanifolds. 
Closed symplectic manifolds always admit polarisations, 
and these have proven useful to study the flexible side of symplectic embedding problems~\cite{buhi,Op13b}. 
In this work, we introduce polarisations of {\it open}\/ symplectic manifolds 
and show that this concept offers new opportunities, at least in dimension~$4$.
Referring 
to~$\S$\2\ref{ss:open} for the definition,  
we focus in this introduction on the applications.

\subsection{Symplectic embeddings} \label{ss:emb}
As usual, $D(a)$, $B^4(a)$, $C^4(a)$, and $Z^4(a)$ denote 
the open round disc of area~$a$, the open 4-ball of capacity $\pi r^2=a$, the open ``cube" $D(a) \times D(a)$, 
and the open cylinder $D(a) \times \CC \subset \CC \times \CC$, all endowed with the standard symplectic form 
$\omega_{\st} = \sum_{j=1}^2 dx_j \wedge dy_j$ on~$\RR^4$
with primitive $\ga_{\st} \,=\, \frac 12 \sum_{j=1}^2 x_j dy_j - y_j dx_j$.

A cornerstone of modern symplectic geometry is Gromov's non-squeezing theorem from~\cite{Gr85}, 
stating in dimension four that $B^4(a)$ symplectically embeds into~$Z^4(1)$ only for~$a \leq 1$.
In~\cite{SSVZ24}, Sackel, Song, Varolgunes, and Zhu 
investigated a refinement of this result, asking how large a set 
(in the sense of the Hausdorff dimension, for instance) 
one has to remove from~$B^4(a)$ such that the complement symplectically embeds into~$Z^4(1)$. 
They first proved that it does not suffice to remove a set
of lower Minkowski dimension~$<2$, and then showed optimality of this result
by proving that for $a \leq 2$ it suffices to remove a Lagrangian coordinate disc from~$B^4(a)$.
Later, Brendel in the appendix of~\cite{SSVZ24} showed that for~$a<3$ it suffices to remove
three Lagrangian pinwheels (certain simple Lagrangian CW-complexes) and a symplectic torus.
In their Question~3, these authors then asked:

\begin{question} \label{q:Umut}
What is the largest $a$ such that the complement of a 2-dimen\-sional subset of~$B^4(a)$
symplectically embeds into $Z^4(1)$\2?
\end{question}

We denote by $\Gamma_{\!\frac 1k}$ the union of the
$k$~half-lines emanating from the origin in~$\C$ and cutting
the plane into isometric sectors, as partly shown in (i), (ii), and~(iii) of
Figure~\ref{fig-grids} 
on page~\pageref{fig-grids} 
for $k = 2,3,4$,
and write 
$$
\Delta_k := \Gamma_{\!\frac 1k} \times \Gamma_{\!\frac 1k} \subset \C^2 .
$$ 
Then $\Delta_k$ is the orbit of the Lagrangian quadrant $\RR_{\geq 0}^2 \subset \CC^2$ 
under the subgroup~$G$ of~$\U (2)$ generated by $(e^{2\pi i /k},1)$ and~$(1,e^{2\pi i /k})$.
When $k$ is even, then $\Delta_k$ is also the $G$-orbit of the 
standard Lagrangian plane $\R^2 \subset \C^2$.

A symplectic embedding $\gf \colon (M, d\alpha) \to (M', d\alpha')$
between exact symplectic manifolds is {\bf $(\alpha, \alpha')$-exact}
if $\gf^* \alpha' = \alpha + df$ for a smooth function~$f$ on~$M$.
Since $B^4(1) \subset C^4(1)$, our first symplectic embedding result
give a positive answer to Question~\ref{q:Umut} for every~$a$.

\begin{thm} \label{t:Umut}
There exists an $(\alpha_{\st}, \alpha_{\st})$-exact symplectic embedding 
\begin{equation*}  
C^4(1) \priv \Delta_k  \, \to \,  Z^4 \bigl(\tfrac 2k \bigr). 
\end{equation*}
\end{thm}

\begin{remarks}
{\rm
(i)
There exists a symplectic embedding of $D(\frac 1k) \times D(\frac 1k)$ into
$C^4(1) \setminus \Delta_k$, so $C^4(1) \setminus \Delta_k$ does not symplectically embed
into $Z^4(A)$ for $A < \frac 1k$ by the nonsqueezing theorem.
Thus Theorem~\ref{t:Umut} is sharp up to a possible factor of~$2$.
We do not know the exact size of the thinnest cylinder into which $C^4(1) \priv \Delta_k$
or $B^4(1) \priv \Delta_k$ symplectically embed,
except for the optimal embedding $B^4(1) \setminus \Delta_2 \to Z^4(\frac 12)$ from~\cite{SSVZ24}.
For a different such embedding see Remark~\ref{rem:volumefill}.

\s
(ii)
A different positive answer to Question~\ref{q:Umut}  
was independently given by Haim-Kislev, Hind, and Ostrover in~\cite{HHO22},
and it is interesting to compare the two constructions, see $\S$\2\ref{ss:comp}.
}
\end{remarks}

Our construction also applies to unbounded domains. 
For instance, consider the union of Lagrangian planes $\Gamma \times \Gamma$ in~$\RR^4$,
where $\Gamma \subset \R^2$ is the square grid 
$$
\Gamma = \bigcup_{(n,m) \in \Z^2} \{n\} \times \R \cup \R \times \{m\}.  
$$

\begin{thm}  \label{t:Remb}
There exists an $(\ga_{\st}, \ga_{\st})$-exact symplectic embedding      
$$
\R^4 \setminus (\Gamma \times \Gamma) \, \to \,  Z^4 (1) .
$$ 
\end{thm}

This answers in particular the question asked by Claude Viterbo to Paul Biran in~1998, 
whether the Gromov width of $\R^4 \setminus (\Gamma \times \Gamma)$ is infinite or not, 
\cite{Vit98}.

We can construct similar embeddings when the source or the target have topology. 
If we take the 4-ball as domain, we can fill ``anything" after removing sufficiently many 
Lagrangian planes: 

\begin{thm} \label{t:anything}
Let $(M,\omega)$ be a connected symplectic 4-manifold of finite volume. 
Let $B^4(a)$ be the ball of the same volume.  
Then for every $\gve >0$ there exist an even number $k \in \NN$ 
such that $B^4(a-\gve) \setminus \Delta_k$ symplectically
embeds into~$(M,\omega)$. 
\end{thm}

%PP
%In the sequel~\cite{OS24} to this work, 
%we combine our embedding method (Theorem~\ref{t:embliouville})
%with the work of Donaldson~\cite{Do96} and Giroux~\cite{Gi17}
%to generalise the previous theorems.
%Referring to Definition~\ref{def:pol} for the notion of polarizing divisor,
%we recall that every symplectic 4-manifold admits polarizing divisors.
%
%\begin{thm} \label{t:galemb} 
%Let $(\Om,\om=d\alpha)$ and $(\Om',\om')$ be 
%exact symplectic 4-manifolds that each can be represented 
%as the complement of a polarizing divisor of a closed
%symplectic 4-manifold, and assume that $\vol \Om < \vol \Om'$.
%
%Then there exists a properly embedded Lagrangian CW complex $\Delta \subset \Om$ such that 
%$\Om \priv \Delta$ embeds into~$\Om'$ by an exact symplectic embedding. 
%\end{thm}
In fact a much more general statement holds. 
Any affine part of a closed symplectic 4-manifold
(defined as the complement of a polarizing divisor, see Definition~\ref{def:pol})
verifies the same kind of flexibility:
After removing a suitable Lagrangian CW-complex, 
it can be embedded into any other affine part of a closed symplectic 4-manifold of larger volume. 
We shall give the proof of this result in a separate paper, \cite{OS24}, 
because in addition to the embedding method developed here
the proof relies on a relative version of Giroux's result in~\cite{Gi17} 
and requires specific techniques. 
%
%PP (PP veut dire: dans OS24 faire) 
%The proof of Theorem~\ref{t:galemb} also shows that given any closed rational
%symplectic 4-manifold $(M,\go)$, the complement of a 2-dimensional
%CW-complex (namely the union of the symplectic polarizing divisor $\Sigma$ and a
%Lagrangian CW-complex~$\Delta$) embeds into~$Z^4(1)$.
%
In contrast to the results of this paper, the removed CW-complex 
cannot be explicitly described, however.

\subsection{Rigidity properties of Lagrangian skeleta} 
Following Biran~\cite{Bi99}, we notice that
if the complement of a Lagrangian CW-complex~$\Delta$ in a domain~$U$ 
symplectically embeds into a cylinder, 
then any domain $V \subset U$ that does not symplectically embed
into that cylinder must intersect~$\Delta$.
The results of the previous paragraph therefore have implications on non-removable intersections 
of~$\Delta$ with various subsets,
%je suis d'accord avec "moin lourd", mais il n'est pas claire ce qui est un plongement symplectique
%d'une sous-variété Lagrangienne. J'ai donc repris le "even". 
%oui mais on parle de non-removable intersections, qui ont lieu avec les sous-variétés Lagrangiennes. Mais ça me va comme ça aussi. 
and even with Lagrangian submanifolds. 
We shall remove Lagrangian CW-complexes more general than the~$\Delta_k$, 
namely products of arbitrary connected graphs.

\begin{definition} \label{def:grid}
{\rm
A {\bf grid}~$\Gamma$ in $D(A)$ is the part in $D(A)$ of a connected graph $\Gamma \cup \pp D(A)$ 
in~$\overline{D(A)}$ that 
has no 1-valent vertex, 
contains the boundary of~$D(A)$, 
and has smooth edges. 
}
\end{definition}

\begin{figure}[h]   
 \begin{center}
  \psfrag{y}{your turn}
	\psfrag{1}{(i)}
  \psfrag{2}{(ii)}
	\psfrag{3}{(iii)}
  \psfrag{4}{(iv)}
	\psfrag{5}{(v)}
  \psfrag{6}{(vi)}
	\psfrag{7}{(vii)}
  \psfrag{8}{(viii)}
	\psfrag{9}{(ix)}
  \leavevmode\includegraphics{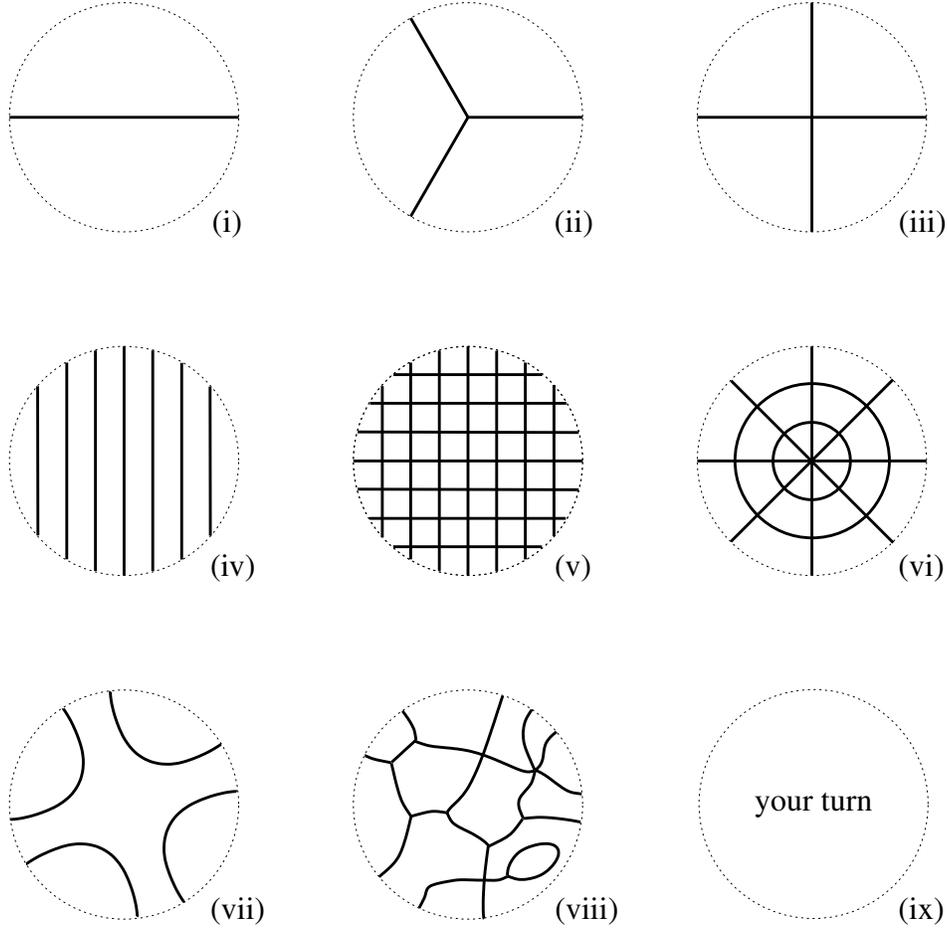}
 \end{center}
 \caption{Examples of grids.}  \label{fig-grids}
\end{figure}

Figure~\ref{fig-grids} shows examples of grids. 
Note that the complement of a grid consists of topological discs,
and that for two grids $\Gamma$ and~$\Gamma'$
the Lagrangian CW-complex $\Gamma \times \Gamma'$ is smooth if and only if 
both $\Gamma$ and~$\Gamma'$ have no vertex, as in~(i) and~(iv).

Cieliebak and Mohnke in~\cite{CM18} defined for every closed Lagrangian submanifold 
$L \subset \RR^4$ its {\bf minimal symplectic area} by 
\begin{equation} \label{e:areamin}
\ca_{\min}(L) := \inf \left\{ \int_\gs \omega_{\st} \mid [\gs] \in \pi_2 (\RR^4, L), \, \int_\gs \omega_{\st} >0 \right\} 
\;\in\; [0,\infty]. 
\end{equation}

\begin{thm} \label{t:lagrig} 
Let $\Gamma_{\!\leq a} \subset D(A)$ and $\Gamma_{\!\leq b} \subset D(B)$ be any two grids 
whose complements are a union of topological discs of area~$\leq \!a$ and~$\leq \!b$, respectively.
Then a closed Lagrangian submanifold~$L$ of $D (A) \times D(B)$ with 
$$\ca_{\min}(L)\geq a+b$$
cannot be mapped to 
$\bigl( D(A) \times D(B) \bigr) \setminus \bigl( \Gamma_{\!\leq a} \times \Gamma_{\!\leq b} \bigr)$ 
by a Hamiltonian diffeo\-morphism of~$\RR^4$.
\end{thm}

\begin{example*}
{\rm
For every Markov triple $\fm$ there exists a monotone Lagrangian torus $\widehat L(\fm)$
in the complex projective plane endowed with the Fubini--Study form integrating to~$3$
over complex lines.
Each such torus gives rise to a Lagrangian torus $L(\fm) \subset B^4(3)$,
and $L(\fm)$ is Hamiltonian isotopic to~$L(\fm')$ in~$B^4(3)$ only if $\fm = \fm'$,
see~\cite{Via17}.
For these Lagrangians, $A_{\min}(L(\fm)) = 1$.
If $\phi_H$ is a Hamiltonian diffeomorphism of~$\R^4$ 
such that $\phi_H(L(\fm)) \subset B^4(A) \setminus \Delta_k$, then
$A \geq \frac k2$.
\diam
}
\end{example*}

%PP also corollary from Theorem 4 ??? (L must intersect $\Sigma \cup \Delta$)

\subsection{Legendrian barriers} 
Lagrangian CW-complexes ~$\Delta$ with symplectically small complement were called 
Lagrangian barriers in~\cite{Bi01} (cf.\ $\S$\2\ref{ss:comp}).
A further consequence of the phenomenon of Lagrangian barriers, 
that apparently went unnoticed, 
is a similar Legendrian barrier phenomenon in contact dynamics. 
A grid is {\bf radial} if it consists of straight rays emanating from the origin.

\begin{thm} \label{t:legbarriers1}
Let $\Gamma_{\! \leq \delta_1}$ and $\Gamma_{\! \leq \delta_2}$ be two radial grids
that divide $D(1)$ into sectors of area~$\leq \delta_1$ and area~$\leq \delta_2$,
respectively. 
Let $U \subset C^4(1)$ be a starshaped domain with smooth boundary~$S$, 
endowed with the usual contact form $\gl_S = \ga_{\st}|_S$,
and consider the connected Legendrian graph 
$$
\Lambda_\gd \,=\, \bigl( \Gamma_{\! \leq \delta_1} \! \times \Gamma_{\! \leq \delta_2} \bigr) 
\cap S .
$$
Then any Legendrian knot $\Lambda$ in~$S$ has a Reeb chord from $\Lambda$ to 
$\Lambda \cup \Lambda_\gd$ of length~$\leq \gd_1+\gd_2$ for the contact form~$\gl_S$. 
\end{thm}

\begin{remarks}              
{\rm 
(i)
While the proof of Theorem~\ref{t:legbarriers1} relies on tools from hard symplectic
geometry, the case of the round sphere can be shown by elementary geometric arguments, 
see~Section~\ref{s:legbar}.

\s
(ii)
The smooth arcs of $\Lambda_\gd$ are indeed Legendrian curves on~$S$,
because $\iota_{r \pp_r} \omega_{\st} = \ga_{\st}$ and $\pp_r$ is tangent to 
$\Gamma_{\! \leq \delta_1} \! \times \Gamma_{\! \leq \delta_2}$.
The theorem says that for any Legendrian knot $\Lambda$, 
every embedded Reeb cylinder of length $\geq \gd_1+\gd_2$ 
over~$\Lambda$ must intersect~$\Lambda_\gd$,
hence the name Legendrian barrier for~$\Lambda_\gd$. 
It follows, in particular, that for $a>\frac 2k$ there exists no exact contact embedding 
$$
\bigl( D(a) \times \RR/\ZZ, \ga_{\st} + dt \bigr) \,\hra\, 
\bigl( S \priv \Delta_k, \lambda_S \bigr) .
$$ 
It is then natural to ask whether there is a purely contact topological statement 
of this kind. 
For instance, using the Reeb flows one constructs exact contact embeddings
$$
\bigl( D(a) \times \RR/\ZZ, \ga_{\st} + dt \bigr) \,\hra\, 
\bigl( S^3(1), \lambda_S \bigr) 
$$ 
for all $a \leq 1$.
Is it true that for $a \in (\frac 2k, 1]$ the image of these embeddings cannot be displaced 
from $\Delta_k \cap S^3(1)$ by a (not necessarily exact) contact isotopy? 

\s 
(iii) 
In fact, the proof will show that there exist both forwards and backwards Reeb chord from $\Lambda$ to $\Lambda\cup \Lambda_\delta$ of length $\leq \delta_1+\delta_2$. 

\s
(iv)
(Improved return time)
Theorem~\ref{t:legbarriers1} states an alternative: 
There must be short Reeb chords from $\Lambda$ to $\Lambda$ {\it or}\/ $\Lambda_\gd$.
This is a relative version of the 
solution of the Arnold chord conjecture by Cieliebak--Mohnke in~\cite{Mo01,CM18}.
They proved that for $\Lambda \subset S = \pp U$ as in the theorem, 
there must be a Reeb chord {\it from $\Lambda$ to $\Lambda$} of length at most
$$
\min \left\{ e(U), \tfrac 12 \right\}
$$
where $e(U)$ denotes the displacement energy of $U$.
Theorem~\ref{t:legbarriers1} improves this bound on the return time 
for all those $\Lambda$ that lie in the complement
of the 
%PP backward  etwas unglücklich, da man sich an (iii) erinnern muss 
Reeb cylinder of~$\Lambda_\gd$ of length~$\gd_1+\gd_2$,
whenever $\gd_1+\gd_2 < \min \left\{ e(U), \tfrac 12 \right\}$.
}
\end{remarks}

\subsection{A general embedding result.}
The upper bounds on the cylindrical capacity of the complement of our
Lagrangian CW-complexes all rely on a single embedding result. 
To state it, we need to introduce the subclass of regular grids.

\begin{definition*}
{\rm
A grid $\Gamma$ on $D(A)$ (as defined in \ref{def:grid})
is {\bf regular} if 
every vertex of $\Gamma \cup \pp D(A)$
has a Darboux chart in which $\Gamma$ consists of radial rays that cut the disc into equal sectors. When the vertex belongs to the boundary, the Darboux chart must be understood with image $D(\eps) \cap \{\im z\geq 0\} \subset \H$. 
}
\end{definition*}

In Figure~\ref{fig-grids},
all grids are regular except for (viii) -- two of the three non-smooth points of $\Gamma \cup \partial D(A)$
have no Darboux chart as required -- 
and possibly~(ix).

\begin{thm} \label{t:main}
Let $\Gamma_{\!\leq a} \subset D(A)$ and $\Gamma_{\!\leq b} \subset D(B)$ be any two regular grids 
whose complements are a union of topological discs of area~$\leq a$ and $\leq b$, respectively.
Then there exists an $(\ga_{\st}, \ga_{\st})$-exact symplectic embedding 
$$
\bigr( D(A) \times D(B) \bigr) \setminus \bigl( \Gamma_{\!\leq a} \! \times \Gamma_{\!\leq b} \bigr)
\, \hra \, Z^4(a+b) .
$$ 
\end{thm}

Theorem~\ref{t:Umut} is a special case of this theorem,
and Theorem~\ref{t:Remb} is a generalisation to a non-compact setting. 
The exactness in this statement is essential to derive 
Theorems~\ref{t:lagrig} and~\ref{t:legbarriers1} on Lagrangian rigidity and Legendrian barriers.

\subsection{Related results and references} \label{ss:comp} $ $
%In this section
%we first compare our approach and our results with those of Biran in~\cite{Bi99, Bi01}, 
%and then put Theorems~\ref{t:Umut} and~\ref{t:lagrig} 
%into the context of other results.

\m \ni
{\bf Biran's Lagrangian barriers.}
As we have already alluded to, the present paper is in many ways a sequel of Biran's work~\cite{Bi01}, 
that is useful to have in mind while reading our paper. 
We therefore briefly recall it, integrating Giroux's results from~\cite{Gi02,Gi17} 
that allow to generalise Biran's results from the K\"ahler to the symplectic setting. 
Let $(M,\om)$ be a closed symplectic mani\-fold with integral symplectic class 
(meaning that $[\om] \in H^2(M,\Z)$). 
For every sufficiently large~$k$, Donaldson~\cite{Do96} produced a symplectic hypersurface~$\Sigma$ 
Poincaré-dual to~$k[\om]$,
called a polarisation of degree~$k$ of~$(M,\omega)$.
The complement $M \setminus \Sigma$ was proven to be Weinstein by Giroux~\cite{Gi17}. 
Hence there exists an isotropic CW-complex (called the skeleton) $\Delta \subset M \priv \Sigma$ 
to which $M \priv \Sigma$ retracts.
Furthermore, any compact subset of $M \priv \Delta$ embeds into  
a ruled symplectic manifold,
namely a symplectic sphere bundle over~$\Sigma$ whose fibers have area~$\frac 1k$, 
see~\cite{Bi99,Bi01}. 
A fundamental example to have in mind along this discussion is the following:

\begin{example} \label{ex:Biran}
{\rm (\cite[$\S$\23.1.3]{Bi01})
In the complex projective space~$\CP^n$ with Fubini--Study form integrating to~1 over complex lines, 
the hypersurface 
\begin{equation} \label{e:S}
\Sigma_{2k} \,=\, \left\{ \sum_{j=0}^{n} z_j^{2k}=0 \right\}
\end{equation}
is a symplectic polarisation of degree~$2k$, 
and one can arrange the skeleton of~$\CP^n \setminus \Sigma_{2k}$ to be the union
$\Delta_{2k}(\CP^n)$ of $k^n$ equally distributed images of~$\RP^n$ by projectivised unitary maps. 
The complement of~$2k$ fibers of the symplectic disc bundle over~$\Sigma_{2k}$ is symplectomorphic to
$B^{2n}(1) \setminus \Delta_{2k}^n$, where 
$$
\Delta_{2k}^n := (\Gamma_{2k})^n  
$$
is the union of $k^n$~unitary images of~$\R^n$. 
For $n=2$ this is our set $\Delta_{2k}^2 = \Delta_{2k}$. 
}
\end{example}

Ruled symplectic manifolds are known to have many symplectic invariants bounded by 
the area of their sphere fibres, 
which allowed Biran to get interesting 
non-removable intersections with the skeleton.  
For instance, for a polarisation of degree~$k$ with skeleton~$\Delta$,
every symplectic ball of capacity~$\geq \frac 1k$ in~$M$ must intersect~$\Delta$. 
He therefore called these skeleta Lagrangian barriers.
In particular, such balls cannot be displaced from~$\Delta$ by a Hamiltonian isotopy. 
Combined with the work~\cite{CM18} by Cieliebak and Mohnke, 
Biran's decomposition result also implies non-displaceability of small Lagrangian submanifolds from the skeleton. 
Namely, applying the neck-stretching in~\cite{CM18} to the set-up in~\cite{Bi99,Bi01}, 
one finds that the general Cieliebak--Mohnke Lagrangian width
$$
\ca_{\min}(L,M) := \inf \left\{ \int_\gs \omega \mid [\gs] \in \pi_2 (M, L), \, \int_\gs \omega >0 \right\} 
$$
can be estimated as follows.

\begin{thm}[Biran--Cieliebak--Mohnke] \label{t:bicimo} 
Let $\Delta$ be the skeleton of a polarisation of degree~$k$ of $(M,\omega)$. 
Then for any closed Lagrangian submanifold $L \subset M \priv \Delta$
it holds that
$$
\ca_{\min}(L,M) < \frac 1k .
$$
In other words, a closed Lagrangian submanifold with $\ca_{\min}(L,M) \geq \frac 1k$ cannot be displaced 
from the skeleton~$\Delta$ by a Hamiltonian isotopy. 
\end{thm}

See \cite[Theorem~3]{Op25} for a proof.
For instance, this theorem combined with a geometric argument explained in~\cite{OS24}
implies that  
for $x \geq \frac 1{2k}$ and $2k>n$ 
the product torus $\T^n(x) \subset B^{2n}(1)$ 
whose factors are circles enclosing area~$x$
cannot be displaced from~$\Delta_{2k}^n$ inside~$B^{2n}(1)$.
Note that for $n=2$ our Theorem~4 gives the same result under the stronger hypothesis $x \geq \frac 1k$.

\medskip
Again, Theorem \ref{t:bicimo} has a Legendrian corollary:
\begin{thm} \label{t:legbarn} 
Let $U \subset B^{2n}(1)$ be a starshaped domain with smooth boundary~$S$, 
endowed with the usual contact form $\gl_S = \ga_{\st}|_S$,
and consider the Legendrian CW-complex $\Lambda_{2k} = \Delta_{2k}^n \cap S$.
Then for every Legendrian submanifold $\Lambda \subset S$ 
there exists a Reeb chord from $\Lambda$ to $\Lambda \cup \Lambda_{2k}$ 
of length $T \leq \frac{1}{2k}$. 
\end{thm}

Biran's decomposition theorem thus implies the same kind of applications that 
we have drawn from our embedding results, in any dimension and with better constants. 
We should therefore explain what this paper adds beyond the embedding results 
from Section~\ref{ss:emb}. 
Here is a first answer. 
Biran's intersection result and Theorem~\ref{t:bicimo} state that {\it some}\/ symplectic capacities
of $M \priv \Delta_k$ are smaller than $\frac 1k$, namely Gromov's ball embedding capacity 
$$
c_B(M) \,:=\, \sup \left\{ a \mid \mbox{$B^{2n}(a)$ symplectically embeds into $M$} \right\} ,
$$
which is the smallest normalised symplectic capacity, 
and the Cieliebak--Mohnke capacity 
$$
c_\lag(M) \,:=\, \sup \left\{ \ca_{\min} (L,M) \mid L \subset M \right\} 
$$ 
where the supremum runs over all closed Lagrangian submanifolds. 
Since the cylindrical capacity 
$$
c^Z(U) \,:=\, \inf \left\{ A \mid \mbox{$U$ symplectically embeds into $Z^{2n}(A)$} \right\}
$$ 
is the largest normalised capacity,  Theorem~\ref{t:Umut}
states instead  that {\it all}\/ symplectic capacities of~$B^4(1) \priv \Delta_k$ are small, 
and our work~\cite{OS24} generalises this property to the affine part of any closed symplectic manifold. 
There are many different capacities, whose estimates have different implications. 
For instance, applied to the Hofer--Zehnder capacity, Theorem~\ref{t:Umut} says that 
the flow of any compactly supported Hamiltonian function on~$C^4(1)$ 
that vanishes on a neighbourhood of 
$\Delta_k$ and whose maximum is at least~$\frac 2k$ has a non-trivial closed orbit of period~$\leq 1$.

Another, more important, gain is that our approach allows to study the rigidity properties of 
certain {\it singular}\/ polarisations. 
The role of singular polarisations for symplectic embeddings has been observed earlier, 
and they were used in several works~\cite{buhi,buhiop,Op13b,Op15}. 
But the study of their skeleta was left aside in these works.
The question of rigidity of Lagrangian skeleta associated to 
singular polarisations is more subtle than
for smooth polarisations.
Indeed, the complement of the skeleton of a singular polarisation
is not ruled anymore,
and the Gromov width of this complement is not related anymore to the degree
but may be arbitrarily close to the Gromov width of the ambient manifold, 
~\cite[Theorem~1]{Op15}.
The present work opens a way to identify singular polarisations
for which the complement of the skeleton has small width.
It may also be worth noting that some very natural Lagrangian CW-complexes are
the skeleta of singular polarisations and not of smooth ones.
For instance, the Clifford torus in~$\CP^2$ is the skeleton of a polarisation by
three lines, but cannot be the skeleton of a smooth polarisation, 
see~$\S$\2\ref{ss:polclosed}.
And even if the skeleton of a smooth polarisation is computable
(like $\Delta_{2k}$ for $\Sigma_{2k}$), it may be much easier to obtain it
from a singular polarisation, see e.g.\ $\S$\2\ref{ss:2x2}.

\b
\ni
{\bf Around Theorems~\ref{t:lagrig} and \ref{t:legbarriers1}.}
Traditionally, Lagrangian intersection results were proven for smooth closed monotone Lagrangians. 
For non-smooth Lagrangians, the Floer machinery does not work directly (see however \cite{EnPo09}), 
and for non-monotone Lagrangians it is more difficult~\cite{FOOO}.
Our Lagrangian skeleta $\Gamma_{\!\leq a} \times \Gamma_{\!\leq b}$ may be smooth or not, are not closed, 
and they may be very far from monotone.

Theorem~\ref{t:lagrig} says that even rather small Lagrangians must intersect 
$\Gamma_{\!\leq a} \times \Gamma_{\!\leq b}$. 
Such Lagrangian rigidity at a small scale 
has been observed recently in~\cite{MS21, PS23}. 
In these papers, the role of our~$\Delta$ is played by Lagrangian submanifolds of~$S^2 \times S^2$
that are products of one circle in the first factor with a collection of circles in the second factor which, 
in contrast to our grids, decompose~$S^2$ 
into components that may have arbitrary topology. 
Another difference is that these Lagrangian submanifolds are secretly monotone 
(when lifted to a symmetric product), 
whereas our skeleta do not need to be monotone in any sense: 
Already each disc $D(A)$ or~$D(B)$ may contain several discs bounded 
by $\Gamma_{\! \leq a}$ or~$\Gamma_{\! \leq b}$ 
whose areas are completely independent. 
Theorem~\ref{t:lagrig} therefore suggests that rigidity at a small scale is not intrinsically related to monotonicity. 

Singular Lagrangian were first studied in \cite{EnPo06, EnPo09}, 
where it is shown that the product of the 1-skeleton of a fine enough triangulation of the 2-sphere is super-heavy 
and hence must intersect any heavy Lagrangian. 
Recall that heavy sets are non-displaceable. 
Our intersection condition on the minimal action is thus very often much weaker.

%J'enleverais toute la suite (jusqu'à % suivant)
%D'accord, c'était un peu avare.
%Theorem \ref{t:bicimo} is stronger than Theorem~\ref{t:lagrig} 
%in the case of the special radial product
%$\Delta_k = \Gamma_{\! \frac 1k} \times \Gamma_{\frac 1k}$.
%On the other hand, 
%$\Delta_k$ is the only explicitly known 
%Lagrangian skeleton in the ball from a smooth polarisation of~$\CP^2$,
%while Theorem~\ref{t:lagrig} applies to all skeleta $\Gamma_{\! \leq a} \times \Gamma_{\! \leq b}$.
%Similarly, in Theorem~\ref{t:legbarriers1} we can take any Legendrian graph
%$(\Gamma_{\delta_1} \times \Gamma_{\delta_2}) \cap S$
%as Legendrian barrier, 
%while Theorem~\ref{t:legbarn} gives a better upper bound for the hitting time
%in case $\delta_1 = \delta_2 = \frac 1k$.

\m \ni
{\bf Theorem \ref{t:main} versus results from~\cite{SSVZ24,HHO22,HHO24}.}
In \cite{SSVZ24} it is shown that the complement $B^4(1) \setminus \Delta_2$ 
of the Lagrangian plane~$\Delta_2$ symplectically embeds into~$Z^4(\frac 12)$.
While Theorems~\ref{t:Umut} and~\ref{t:main} do not provide such an embedding, 
a modification of our construction yields a symplectic embedding
$B^4(1) \setminus \Delta_2 \to Z^4(\frac 12+\gve)$ for every $\gve >0$, 
see Remark~\ref{rem:volumefill}.

A different positive answer to Question~\ref{q:Umut} 
was given by Haim-Kislev, Hind, and Ostrover afterwards in~\cite{HHO22}.
They showed that for every~$a$
one can remove from the ball~$B^4(a)$ a finite number (depending on~$a$) 
of parallel {\it symplectic}\/ planes such that the complement symplectically embeds into~$Z^4(1)$.
Their embedding has very different properties from ours: 
The removed set is symplectic instead of Lagrangian, 
and the construction is rigid in the sense that one cannot alter the position of the planes, 
while our grids are rather arbitrary in view of Theorem~\ref{t:main}. 
The key difference is that in contrast to theirs, our embeddings are exact, 
a property that we crucially need for deriving Theorems~\ref{t:lagrig} and~\ref{t:legbarriers1}.

In \cite{HHO24}, the same authors compute the cylindrical capacity of the complement $B^4(1) \setminus P_t$ of a plane with K\"ahler angle $t \in [0,1]$ 
to be~$\frac{1+t}2$, 
interpolating between the previously known case of the Lagrangian plane~$\Delta_2 = P_0$ from~\cite{SSVZ24}
and the flexible case of a complex plane.
By contrast, our embedding method leads only to barriers which are Lagrangian skeleta. 
It would be interesting to know whether these Lagrangian skeleta also belong to natural 1-parameter families 
of sets which do not remain Lagrangian, but still retain some barrier property.

%Now take the union of equidistant Lagrangian planes $\Gamma \times \Gamma$
%with $\Gamma$ the grid in Figure~\ref{fig-grids}~(iv).
%If $a$ is the largest area of a component of $D(1) \setminus \Gamma$, then Theorem~\ref{t:main}
%shows that $B^4(1) \setminus (\Gamma \times \Gamma)$ symplectically embeds into $Z^4(2a+\gve)$ for every $\gve >0$.
%In particular, $c \bigl(B^4(1) \setminus (\Gamma \times \Gamma) \bigr) \to 0$ as $a \to 0$ for every symplectic capacity~$c$.
%Such an asymptotic result has been obtained in~\cite[Theorems 1.3 and~1.4]{HHO24}
%for every K\"ahler angle $t \in [0,1]$.

\medskip
\noindent 
{\bf Organisation of the paper.}
The paper is organised as follows. 
In Section~\ref{s:pol} we recall the results from~\cite{Op13a, Op13b} on symplectic polarisations 
in dimension four, 
introduce Liouville polarisations and explain their main relations with symplectic embeddings. 
In Section~\ref{s:explicit} we compute some explicit pairs polarisations/skeleta. 
Combining these tools, we prove our embedding results in Section~\ref{s:sympemb}, 
deduce the Lagrangian rigidities in Section~\ref{s:lagrig} and the Legendrian ones in Section~\ref{s:legbar}.

\medskip
\noindent
{\bf Notation and conventions.}
Since the paper is long, we list here some of our terminology and conventions, 
to which the reader may return when necessary.

\begin{itemize}
\item[\bul]
A symplectic embedding $\gf \colon (M, d\alpha) \to (M', d\alpha')$
between exact symplectic manifolds is {\bf $(\alpha, \alpha')$-exact}
if $\gf^* \alpha' = \alpha + df$ for a smooth function~$f$ on~$M$.
Equivalently, $\int_{\gamma} \alpha = \int_{\gf (\gamma)} \alpha'$ 
for all closed curves~$\gamma$ in~$M$.
It is enough to check this equality of actions on a set of closed curves that generate~$H_1(M;\RR)$.
If $H^1(M;\RR) = 0$, then $(\alpha, \alpha'$)-exactness is automatic;
and otherwise, the notion of $(\alpha, \alpha'$)-exactness may depend on the choice of primitives~$\alpha,\alpha'$.
We sometimes write $\hraa$ to abbreviate ``there exists an $(\ga_{\st}, \ga_{\st})$-exact symplectic embedding".

\s
\item[\bul]
A {\bf Liouville form} on a symplectic manifold $(M,\omega)$
is just a primitive~$\lambda$ of the symplectic form: $d\lambda = \omega$.
The associated {\bf Liouville vector field}~$X_\lambda$ is defined by 
$\iota_{X_\lambda} d\lambda = \lambda$.
Its flow, the {\bf Liouville flow}, is conformally symplectic:
$\left( \phi_{X_\lambda}^t \right)^* \omega = e^t\, \omega$.

\s
\item[\bul]
By a {\bf symplectic curve} in a symplectic 4-manifold we mean a 2-dimen\-sion\-al
embedded symplectic submanifold,
which is closed if the ambient symplectic manifold is, or properly embedded if it is open.  
A {\bf normal crossing} between symplectic curves $\Sigma, \Sigma'$
is an intersection point $p \in \Sigma \cap \Sigma'$ such that $T_p\Sigma$ and~$T_p\Sigma'$
are $\omega$-orthogonal.

\s
\item[\bul]
A {\bf symplectic multi-curve with normal crossings}~$\Sigma$ is a union of symplectic curves~$\Sigma_i$ 
(called the components of~$\Sigma$) whose pairwise intersections are all normal crossings. 
The singularities of a multi-curve is the set of intersection points between components and is denoted $\sing(\Sigma)$. 
The complement of the singular locus is called the regular part of~$\Sigma$. 

\s
\item[\bul]
Usual objects of differential geometry are generalised to  multi-curves
  by just taking a collection of such objects on each component. 
For instance, a differential $k$-form $\alpha$ on $\Sigma=\bigcup \Sigma_i$ is a collection of differential $k$-forms on the $\Sigma_i$. 
The morphisms between 
multi-curves
 $\Sigma, \Sigma'$ are smooth maps from $\Sigma$ to $\Sigma'$ 
(in the sense that they are smooth when restricted to each components), 
and send the regular part and singular locus of~$\Sigma$ to those of~$\Sigma'$. 
A symplectic embedding between multi-curves is an injective such morphism that pulls back $\om'|_{\Sigma'}$ to $\om|_{\Sigma}$. 

\s
\item[\bul]
A {\bf weighted} symplectic multi-curve $\bfsigma := \{(\Sigma_i,\mu_i)\}$ 
with normal crossings 
is a collection of symplectic curves $\Sigma_i$ 
weighted by real numbers $\mu_i$, whose total space $\Sigma:=\bigcup \Sigma_i$ is a symplectic multi-curve  with normal crossings. 

\s
\item[\bul]
A {\bf symplectic morphism} from $\bfsigma$ to $\bfsigma'$
is a symplectic embedding of $\Sigma$ into $\Sigma'$ that sends each component to a component of the same weight and sends the regular part and singular locus of~$\Sigma$ to those of $\Sigma'$.
We write $\bfsigma \to \bfsigma'$.

\s
\item[\bul] 
Symplectic multi-curves will also be called {\bf symplectic divisors}, or simply divisors.

\s
\item[\bul]
Given a subset $S$ of a manifold $M$, we write $\Op(S,M)$ or just $\Op(S)$
instead of ``some open neighbourhood of~$S$ in~$M$".
\end{itemize}

%PP \section{Main idea of the proof of Th 1.}
%Roughly, we look at collections of properly embedded symplectic hypersurfaces $\Sigma_j \subset (M,\omega)$
%such that on $M \setminus \bigcup \Sigma_j$ the symplectic form~$\omega$ admits a primitive 1-form $\lambda$
%that is well-behaved near $\bigcup \Sigma_j$ and whose Liouville flow is ``tangent to the boundary".

\m
\ni
{\bf Acknowledgment.}
EO thanks Kai Cieliebak for an interesting discussion during a visit to Strasbourg University. 
FS is grateful to
Viktor Ginzburg, Basak G\"urel, Marco Mazzucchelli, and Alfonso Sorrentino
for organizing the workshop Symplectic Dynamics at Rome in May~2023, 
where some of the results of this work were presented, 
and to Umut Varolgunes for many interesting discussions and a great week
in June~2023 at Istanbul. 
We also thank Umut for pointing out a gap in an earlier version,
and Paul Biran for pointing out reference~\cite{Vit98}. 
FS also cordially thanks Yasha Eliashberg and the Department of Mathematics at Stanford University
for their wonderful hospitality during the spring term~2024.

\section{Polarisations of symplectic manifolds} \label{s:pol}

In this section we review those results from \cite{Bi01,Op13a,Op13b,Op15} 
that we need in this paper, in a form useful for us. 
For simplicity, we stick to dimension~$4$, 
although most of the statements can be adapted to higher dimensions.

\subsection{Neighborhoods of symplectic curves in dimension~$4$}\label{s:sdb}

Let $\pi \colon \cl \to (\Sigma,\tau)$ be a complex line bundle with first Chern class~$c_1$ over a symplectic curve 
of area $\ca = \int_\Sigma \tau$. 
The multiplication by~$e^{2\pi i\theta}$ defines a vector field $\frac\partial{\partial \theta}$ which, 
together with the complex structure on the fiber, 
provides a closed $1$-form~$d\theta$ on the punctured fibers~$F \setminus \{0\}$. 
We now in addition endow~$\cl$ with a Hermitian metric. 
This provides a radial coordinate~$r$, 
and throughout the paper we write $R := r^2$. 
Chern--Weil theory guarantees the existence of a connection $1$-form with curvature 
$\frac {c_1}\ca \, \tau$. 
This is a $1$-form~$\Theta$ on~$\cl \setminus \Sigma$ 
that satisfies
$$
\Theta |_{F \setminus \{0\}} = d\theta \,\mbox{ for each fiber $F$ }
\quad \mbox{ and } \quad
d\Theta = -\frac {c_1} \ca \, \pi^* \tau \,\mbox{ on } \cl . 
$$
The $2$-form 
\begin{equation} \label{e:c1RA}
\om_0 \,:=\, \pi^*\tau + d(R\Theta) 
\,=\, \left( 1-\frac{c_1R}\ca \right) \pi^* \tau + dR \wedge \Theta
\end{equation}
is defined on all of $\cl$, is closed, and is symplectic on $\bigl\{ R<\frac\ca{c_1} \bigr\}$. 
Furthermore, $\omega_0$ is exact on $\cl \setminus \Sigma$ 
(where we identify ${\bf 0}_\cl$ with $\Sigma$),
with Liouville form 
\begin{equation} \label{e:la0}
\lambda_0 := \left( R-\frac\ca{c_1} \right) \Theta.
\end{equation}
The associated Liouville vector field is 
\begin{equation} \label{e:Li0}
\left( R-\frac{\ca}{c_1} \right) \frac\partial{\partial R} .
\end{equation}

The 1-form $d\theta$ is defined on the punctured fibers $F \setminus \{0\}$,
and in a local trivialisation of~$\cl$ over $U \subset \Sigma$ 
it is defined on~$\cl |_U \setminus U$,
but not on all of~$\cl \setminus \Sigma$.
A transition function between two local trivialisations
is of the form $(z,w)\mapsto (z,g(z)\2w)$, 
where $g$ is a smooth $\CC \setminus \{0\}$-valued function, 
and $d\theta$  is pulled back by this map to $d\theta + \im \frac{dg}g$. 
We can therefore define the {\it angular class} $\{ d \theta \}$
as the set of $1$-forms on $\cl \setminus \Sigma$ that locally are equal to $d \theta$
up to a smooth 1-form on~$\cl$.
The form $\Theta$ is a representative of~$\{d\theta\}$. 
Recall that we have fixed a Hermitian metric on~$\cl$.

\begin{definition} 
The \emph{symplectic disc bundle} over $(\Sigma,\tau)$ with Chern class $c_1$ is 
$$
\sdb (\Sigma,\tau,c_1) \,:=\, \left( \left\{ R < \frac \ca{c_1} \right\}, \om_0, \lambda_0 \right) 
\,\subset\, \cl . 
$$
It comes with a distinguished class of $1$-forms on $\sdb(\Sigma,\tau,c_1) \setminus \Sigma$ 
modulo smooth forms on~$\sdb(\Sigma,\tau,c_1)$, called the angular class. 
Any connection $1$-form provides a representative of the angular class. 
For $\eps \leq \frac \ca{c_1}$, we denote   
$$
\sdb_\eps (\Sigma,\tau,c_1) \,:=\, \big(\{R<\eps\},\om_0,\lambda_0\big) \,\subset\, \cl .
$$
\end{definition}

By the symplectic neighbourhood theorem, 
for any symplectic curve $\Sigma \subset (M^4,\om)$ and for $\eps >0$ sufficiently small, 
an open neighbourhood of~$\Sigma$ is symplectomorphic to
$$
\sdb_\eps(\Sigma,M) \,:=\, \sdb_\eps \bigl( \Sigma,\om |_{\Sigma}, c_1 (T\Sigma^{\om}) \bigr)
$$
by a symplectomorphism that lifts any given symplectomorphism between 
the zero section of~$\sdb_\eps(\Sigma,M)$ and $\Sigma \subset M$.

\subsection{Polarisations (closed case)}
\label{ss:polclosed}
In \cite{Bi01}, Biran introduced the notion of polarisation for closed symplectic manifolds with rational symplectic class. 
This definition was later extended in~\cite{Op13b} as follows.

\begin{definition} \label{def:pol}
A \emph{polarisation} ${\bf \Sigma} := \{ (\Sigma_i,\mu_i) \}$ 
of a closed $4$-dimensional symplectic manifold $(M,\omega)$ 
is a finite collection of closed 2-dimensional symplectic submanifolds~$\Sigma_i$ 
that intersect $\omega$-orthogonally, 
weighted by real coefficients~$\mu_i>0$, such that
\begin{equation}\label{eq:wpd}
[{\bf \Sigma}]  \,:=\, \sum_i \mu_i \2 [\Sigma_i] \,=\, \PD ([\omega])  \,\in\,  H_2(M;\RR) . 
\end{equation}
In other terms, a polarisation is a weighted
symplectic divisor with normal crossings Poincar\'e-dual to~$[\om]$ 
in the sense of~\eqref{eq:wpd}.
We write $\Sigma := \bigcup_i \Sigma_i$.
We say that the polarisation is \emph{smooth} if it consists of a single component. 
In this case, $\mu$ is determined by~$\Sigma$ and~$[\omega]$, and we say that $\Sigma$ has degree $d = 1/\mu$. 
\end{definition}

Let $(M,\om,{\bf \Sigma})$ be a polarised closed symplectic manifold. 
The symplectic form~$\omega$ is exact on the complement of~$\Sigma$,
hence $\omega$ has primitives on~$M \setminus \Sigma$.
Throughout the paper, we consider only Liouville forms that satisfy a regularity assumption near~$\Sigma$, 
that we discuss now.

When ${\bf \Sigma}$ is smooth, a neighbourhood can be modeled on 
$$
\sdb_\eps (\Sigma,M) \,:=\, \sdb_\eps \big(\Sigma,\om |_{\Sigma}, c_1(T\Sigma^{\om})\big)
$$ 
for $\eps$ small enough, and is therefore endowed with a radial coordinate~$R$ 
and a connection $1$-form~$\Theta$, as in Section~\ref{s:sdb}. 
When ${\bf \Sigma}$ consists of several pairwise $\omega$-orthogonal components~$\Sigma_i$, 
we follow~\cite[Section 3.1]{Op13b} and make a neighbourhood of~$\Sigma$ a symplectic plumbing 
of the $\eps$-disc bundles over the $\Sigma_i$, which we still denote $\sdb_\eps(\Sigma,M)$: 
a deleted neighbourhood of each~$\Sigma_i$ can be endowed with local coordinates $(R_i,\theta_i)$,
$d\theta_i$ belonging to the angular class around $\Sigma_i$,
such that on a neighbourhood of an intersection point of $\Sigma_i$ and~$\Sigma_j$ 
we have $\Sigma_i=\{R_j=0\}$, $\Sigma_j=\{R_i=0\}$ and $\om=d(R_id\theta_i+R_jd\theta_j)$.  

\begin{definition}  \label{def:tame}
A $1$-form $\lambda$ on $M \setminus \Sigma$ is {\bf tame} 
\begin{itemize}
\item[$\bullet$] 
{\bf at a regular point} $p \in \Sigma_i \setminus \bigcup_{j \neq i} \Sigma_j$,
if there exists a real number~$a_i$  
and a bounded $1$-form~$\lambda'$ on $\op(p) \setminus \Sigma_i$
such that $\lambda = a_i \, d \theta_i + \lambda'$ on $\op(p) \setminus \Sigma_i$ 
($\lambda'$ is then smooth on $\op(p)\priv \Sigma_i$).
In other terms, $\lambda$ locally represents a multiple of a connection form modulo bounded forms.

\s
\item[$\bullet$] 
{\bf at an intersection point} $p \in \Sigma_i \cap \Sigma_j$, 
if there exist real numbers~$a_i, a_j$ 
such that 
$\lambda = (R_i+a_i) \,d\theta_i + (R_j+a_j) \,d\theta_j$ 
on $\op(p) \setminus (\Sigma_i \cup \Sigma_j)$ in the plumbing coordinates.

\s
\item[$\bullet$] 
{\bf along $\Sigma$}, if it is tame at each point of $\Sigma$.  
\end{itemize}
\end{definition}

\begin{remarks} \label{rk:tamedepends} 
{\rm
(i)
The class of tame forms around $\Sigma$ depends on 
the choice of the Hermitian metrics on the~$\cl_i$,
and on the identification of our neighbourhood of~$\Sigma$ in~$M$ 
with the plumbing of the $\eps$-disc bundles over~$\Sigma_i$. 
Given a polarisation, we fix these choices,  
%a fixed such identification,
whether implicitely or after a construction. 
The notion of tameness will always refer to these choices. 

\s
(ii)   
Liouville forms in the complement of a normal crossing divisor $\Sigma \subset M$ have been used in many works, for instance in \cite{McL12, BSV22}. 
In these works, the main interest was the symplectic homology of $M \setminus \Sigma$, 
whence it sufficed to understand that the Liouville flow is pointing towards~$\Sigma$.
Variants of the notion of tameness, that is crucial here, 
were introduced in~\cite{Op15}.
\diam
}
\end{remarks}

It is not hard to check that the number $a_i$ depends only on the component~$\Sigma_i$ 
and not on the point~$p$ on this component
(see~\cite[Lemma~4.1]{Op15}).
We call $a_i$ the residue of~$\lambda$ at $\Sigma_i$ and denote it $\res(\lambda,\Sigma_i)$. 
By \cite[Lemma~3.2]{Op13b} plumbings of symplectic disc bundles have tame 
Liouville forms with arbitrary residues. 
When, moreover,
the homological polarizing condition~\eqref{eq:wpd} holds for $\mu_i=-\res(\lambda,\Sigma_i)$, 
then these tame Liouville forms extend to Liouville forms
on~$M \setminus \Sigma$, see the proof of \cite[Lemma~4.1~(iii)]{Op15}.
For further reference we state:

\begin{proposition} \label{p:singliouville} 
Let $(M,\om,{\bf \Sigma})$ be a polarised closed symplectic manifold 
with ${\bfsigma =\{(\Sigma_i,\mu_i)\}}$. 
Then
\begin{itemize}
\item[{\rm (i)}] 
There exists a tame Liouville form $\lambda$ on $M \setminus \Sigma$ 
with residues $\res(\lambda,\Sigma_i) = - {\mu_i}$.

\s
\item[{\rm (ii)}] 
Let $U$ be an open subset of $M$. Two Liouville forms on~$U$ tame along~$\Sigma$ with the same residues 
differ on $U\priv \Sigma$ by a bounded smooth closed 1-form~$\vartheta$. 

\s
\item[{\rm (iii)}] 
Let $\lambda$ be a tame Liouville form on $M \priv \Sigma$.
For every small enough $\gve >0$ the vector field $X_{\lambda}$ is pointing inwards along
the boundary of~$\SDB_{\gve} (\Sigma,M)$,
and for every point $p \in \SDB_{\gve} (\Sigma,M) \setminus \Sigma$ 
the flow line of $X_{\lambda}$ starting at~$p$
hits $\Sigma$ in finite positive time. 
The Liouville field~$X_\lambda$ is therefore backward complete, meaning that any trajectory is defined on~$\R_{\leq 0}$.
\end{itemize}  
\end{proposition}

\proof
(i) has been discussed above
and (ii) holds by definition of tameness. 
We prove~(iii). Let $\lambda$ be a Liouville form tame along~$\Sigma$. 
Near a regular point of~$\Sigma_i$, $\lambda$ is then equal to $-\mu_i \1 \Theta_i$ 
modulo a bounded one-form on $M \priv \Sigma_i$, and we have seen that $\Theta_i$ is 
$\om$-dual to $\frac \partial{\partial R_i}$, so 
$$
X_\lambda \,=\, -\mu_i \, \frac\partial{\partial R_i} + Z
$$
for a bounded vector field $Z$. 
Since $\frac \partial{\partial R_i}$ has norm of order $\frac 1{r_i}$ and since $\mu_i>0$, 
$X_\lambda \cdot R_i<0$ near $\Sigma_i=\{R_i=0\}$. 
A similar argument applies near a singular point.
\proofend

\begin{definition}[Biran decomposition]
Given a polarised closed symplectic manifold $(M,\om,\bfsigma)$ 
and a tame Liouville form~$\lambda$ on $M \setminus \Sigma$, we denote 
\begin{eqnarray*}
\cb (\bfsigma,\lambda,M) &:=& \Sigma \cup \{p \in M \mid \mbox{$\exists$ $t(p)>0$,\; } 
\lim_{t \nearrow t(p)} \phi^t_{X_\lambda}(p) \in \Sigma\}, 
\\
\skel(\bfsigma,\lambda,M) &:=& M \setminus \cb(M,\bfsigma,\lambda) . 
\end{eqnarray*}
More generally, for any subset $X\subset \Sigma$, we define 
$$
\cb(X,\lambda,M):=X\cup \{p \in M \mid \mbox{$\exists$ $t(p)>0$,\; } 
\lim_{t \nearrow t(p)} \phi^t_{X_\lambda}(p) \in X \} .
$$
There is a natural continuous projection map $\pi \colon \cb(X,\lambda,M) \to X$ that associates to~$x$ the point on~$X$ 
to which its flow line arrives. 
\end{definition}

By (iii) of Proposition~\ref{p:singliouville}, the set $\cb(\bfsigma,\lambda,M)$
is open. It is the {\bf basin of attraction} of~$\Sigma$ under the Liouville flow of $\lambda$.
The closed set $\skel(\bfsigma,\lambda,M)$ is the maximal compact subset of~$M \setminus \Sigma$
that is invariant under the Liouville flow. 
It is called the {\bf skeleton} of~$(M,\om)$ with respect to 
the polarisation~$\bfsigma$ and the Liouville form~$\lambda$. 
Finally, point~(iii) also guarantees that $\overline{\cb(X,\lambda,M)} \cap \Sigma = \overline X$. 

\m
The following examples should be useful for readers unfamiliar with these notions.
  
\begin{example} \label{ex:galpolB} 
{\rm 
If $(\Sigma,\mu)$ is a smooth polarisation of $(M,\omega)$, then the basin
$\cb (\bfsigma,\lambda,M)$ is symplectomorphic to the symplectic disc bundle
$\sdb(\Sigma,M) = \sdb_\mu \bigl(\Sigma,\om |_{\Sigma}, c_1(T\Sigma^{\om}) \bigr)$.
It therefore depends only very mildly on~$M$ or~$\lambda$ from a symplectic perspective.
}
\end{example}

\begin{examples}\label{ex:polp2}
{\rm 
As before, endow
the complex projective plane $\CP^2$ with the Fubini--Study symplectic form~$\om_\FS$ 
that integrates to~$1$ over every projective line. 
Every smooth algebraic curve in $\CP^2$ is symplectic and (since $H_2(\CP^2)$ has rank~$1$)
is a polarisation when weighted by  the inverse of its algebraic degree. For instance:

\begin{itemize}
\item[$\bullet$] 
When $\Sigma$ is a line, then $(\CP^2 \setminus \Sigma,\om_\FS)$ is symplectomorphic to 
the ball $(B^4(1),\om_{\st})$. 
Using the usual toric description of $(\CP^2,\omega_{\FS})$, 
one checks that the Liouville form $\lambda_0 = \ga_{\st} = R_1 \,d\theta_1 + R_2 \,d\theta_2$ 
defined on $B^4 (1) \cong \CP^2 \priv \Sigma$ is tame along~$\Sigma$. 
The skeleton of~$\lambda_0$ is a single point, namely the origin of~$B^4(1)$.

\s
\item[$\bullet$] 
When $\Sigma = \Sigma_k := \{ z_0^{k} + z_1^{k} + z_2^{k} = 0\} \subset \CP^2$, then
$\Sigma$ is the vanishing locus of a holomorphic section~$s_k$ of the Hermitian line bundle~$\co(k)\to \CP^2$ of curvature $-k \2 \om_{\FS}$. 
The form $-\frac 1kd^c\log\|s_k\|$, defined on $\CP^2 \priv \Sigma$, 
provides a tame Liouville form with residue~$-\frac 1k$. 
%Its skeleton can be computed to be~$\Delta_k(\CP^2)$.
%FF ajouté:
For $k$ even, its skeleton is the set $\Delta_k(\CP^2)$ described in Example~\ref{ex:Biran},
and for all $k$ the intersection of the skeleton with $B^4(1) = \CP^2 \setminus \CP^1$
is~$\Delta_k$.
\end{itemize}

For these examples in~$\CP^2$, notice the discrepancy between the ``algebraic degree" of a curve 
and its degree as a polarisation (the ``symplectic degree"). 
These two degrees coincide when the Fubini--Study form is normalised so that 
the area of a line is~$1$, 
but they differ for different normalisations. 
For instance, if $\CP^2$ is obtained as the compactification of~$B^4(2)$, 
so that the symplectic form becomes $2 \,\om_\FS$, a curve of algebraic degree~$2$ 
is a polarisation of degree~$1$. 
}
\end{examples}

\begin{examples}  $ $
{\rm
\begin{itemize}
\item[$\bullet$]
$\left( \CP^1 \times \CP^1, a \2 \om_\FS \oplus b \2 \om_\FS \right)$ 
is polarised by 
$$
{\bf \Sigma_1} \,:=\,
\Bigl\{ 
        \left( \CP^1 \times \{0\},b \right) , 
        \left( \{0\} \times \CP^1,a \right) 
\Bigr\} .
$$ 
For $\lambda = (R_1-a)\,d\theta_1 + (R_2-b)\,d\theta_2$
the skeleton is the single point~$\{(\infty,\infty)\}$.
%F check sign!
 
\m
\item[$\bullet$] 
$\left( \CP^1 \times \CP^1, a \2 \om_\FS \oplus b \2 \om_\FS \right)$ is polarised 
by ${\bf \Sigma_2} :=$ 
$$
\Bigl\{ \left(\CP^1 \times \{0\}, \tfrac b2 \right) , 
        \left( \{0\} \times \CP^1, \tfrac a2 \right),
        \left( \CP^1 \times \{\infty\}, \tfrac b2 \right),
				\left( \{\infty \} \times \CP^1, \tfrac a2 \right) 
\Bigr\} .
$$
For a suitable choice of $\lambda$
the skeleton is the Clifford torus, i.e.\ the product of the two equators. 

\m
\item[$\bullet$] 
$\left( \CP^2, \om_\FS \right)$ is polarised by 
$$
\bfsigma \,:=\, 
\Bigl\{
  \left( \{ z_0=0\},\tfrac 13 \right) , 
  \left( \{z_1=0\},\tfrac 13 \right) , 
	\left( \{z_2=0\},\tfrac 13 \right)
\Bigr\} . 
$$
For a suitable choice of $\lambda$ the skeleton is the Clifford torus. 
\end{itemize}
}
\end{examples}

We end this section by keeping a promise made in the introduction:

\begin{proposition}
There exists no smooth polarisation of $\CP^2$ whose skeleton is the Clifford torus.
\end{proposition}

\proof
Assume that the Clifford torus $L$ is the skeleton of a smooth polarisation
of~$\CP^2$ of degree~$k$.
By~\cite{Bi99}, the Gromov width of the complement is $c_{\G}(\CP^2 \setminus L) \leq \frac 1k$.
On the other hand, 
$\CP^2 \setminus L$ is symplectomorphic to the standard symplectic disc bundle
$\SDB_{\frac 1k} (\Sigma, \tau, k)$ over a closed surface of area~$k$ with disc fibres 
of area~$\frac 1k$. 
By \cite[Proposition 1.3]{Op13a} 
this disc bundle contains a symplectic ellipsoid $E(k, \frac 1k)$
and hence a ball $B^4(\frac 1k)$. It follows that 
$c_{\G}(\CP^2 \setminus L) = \frac 1k$.
This is in contradiction to the fact that $c_{\G}(\CP^2 \setminus L) = \frac 23$, 
proven in \cite[p.\ 2887]{BC09}.
\proofend

\subsection{Liouville polarisations (open case)} \label{ss:open}

Let $(\Om,\om)$ be a connected symplectic manifold 
(without boundary) and 
with an exact symplectic form~$\omega$.
Let $\Sigma \subset \Om$ be a symplectic divisor with normal crossings. 
This now means that each of the finitely many components~$\Sigma_i$ is a 
properly embedded symplectic surface 
and that the singular locus $\sing(\Sigma)$ is compact.

\s \ni
{\bf Tameness} of a 1-form~$\lambda$ on~$\Omega \setminus \Sigma$ along~$\Sigma$ is 
defined as in the closed case (Definition~\ref{def:tame}),
in terms of Hermitian line bundles $\cl_i \to \Sigma_i$ modeled on the symplectic normal bundles, 
and an identification of a neighbourhood of~$\Sigma$ in~$\Omega$ with the plumbing of disc bundles 
over the~$\Sigma_i$.

\begin{definition} \label{def:polliouville} 
{\rm
A \emph{\bf Liouville polarisation} $(\bfsigma,\lambda)$ of~$(\Om, \om)$ 
consists of $\bfsigma = \{(\Sigma_i,\mu_i)\}$, where $\Sigma = \bigcup_i \Sigma_i \subset \Omega$ 
is a symplectic divisor as above and $\mu_i >0$,
together with a Liouville form~$\lambda$ on~$\Om \setminus \Sigma$ 
that is tame along~$\Sigma$, with residues $-\mu_i$ at~$\Sigma_i$,
and such that the Liouville flow~$\phi_{X_\lambda}^t$ on $\Omega \setminus \Sigma$ 
is backward complete and ``forward complete up to hitting~$\Sigma$":
For every point $p \in \Omega \setminus \Sigma$, the flow line starting at~$p$ is defined
for $t \in (-\infty, t^+(p))$, where either $t^+(p) = + \infty$
or $\lim_{t \nearrow t^+(p)} \phi_{X_\lambda}^t(p) \in \Sigma$. 

\s
The points $p$ with $t^+(p) < \infty$ form the basin of attraction $\cb (\bfsigma, \lambda,\Om)$.
The complement of the basin of attraction of~$\Sigma$ under the Liouville flow is  
the \emph{\bf skeleton} of~$(\Om,\omega, \bfsigma, \lambda)$. 
This is the set of points whose Liouville flow line is defined for all times.
If $\bfsigma$ has only one component $\Sigma$, we call the polarisation \emph{\bf smooth}, 
and in analogy with the closed case say that it has \emph{\bf degree~$d=1/\mu$}. 
}
\end{definition}

\begin{remarks} 
{\rm 
(i) In the context of closed manifolds, we have called polarisation the sole weighted divisor~$\bfsigma$, which could be easily separated from the Liouville form. In the context of open manifolds, 
the situation is more intricate because the homological condition $[\bfsigma]=\PD([\om])$ 
does not make sense anymore. 
This condition is replaced by the residue and integrability conditions on~$\lambda$ and its Liouville flow, 
making impossible to separate $\bfsigma$ from~$\lambda$ in the above definition. 

\s
(ii)
Whether or not an exact symplectic manifold $(\Omega, \omega)$ has a Liouville polarisation 
does not seem to be easy to decide, in general. 
We think that if $\Omega$ is the interior of a compact manifold with smooth boundary, 
this only holds when $(\Omega,\omega)$ is the interior of a very special kind of Liouville domain 
$(\overline \Omega, d \beta)$. 
Note that for Liouville domains, the Liouville vector field~$X_\beta$ is {\it transverse}\/ to the boundary, pointing outwards. 
This is in sharp contrast to our Liouville vector fields~$X_\lambda$, 
that in view of the integrability condition on~$X_\lambda$ are ``{\it tangent}\/ at infinity".

\s 
(iii) 
The reader may wonder why we insist on assuming $\Om$ open, 
and propose a definition that may look a bit artificial.
Indeed, if $\Om$ were assumed compact with smooth boundary, we could as well define our Liouville polarisation in terms of a Liouville form on $\Om \priv \Sigma$ whose Liouville vector field is tangent to~$\partial \Om$. 
This makes a perfectly sound definition in this case, that would however discard important examples that we want to consider: 
open manifolds whose boundaries have corners, like the bidisc and, more generally, 
symplectic disc bundles over open surfaces as in Example~\ref{ex:ext} below, 
and~$\R^4$, which is not compact.
}
\end{remarks}

\begin{figure}[h]   
 \begin{center}
  \psfrag{O}{$\Omega$}
	\psfrag{Oh}{$\widehat \Omega$}
	\psfrag{S}{$\Sigma$}
	\psfrag{Sh}{$\widehat \Sigma$}
	\psfrag{X}{$X_{\widehat \lambda}$}
  \leavevmode\includegraphics{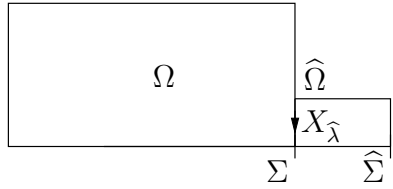}
 \end{center}
 \caption{An extension $(\widehat \Omega, \widehat \Sigma)$ of $(\Omega, \Sigma)$.}  \label{fig-OO}
\end{figure}

\begin{definition}\label{d:extendable}
{\rm
A Liouville polarisation $(\bfsigma,\lambda)$ of~$(\Om, \om)$ is {\bf extendable} if
there exists an exact symplectic manifold $(\widehat \Omega, \widehat \omega)$
with $\Omega \subset \widehat \Omega$ and $\widehat \omega |_\Omega = \omega$,
a symplectic divisor $\widehat \Sigma \subset \widehat \Omega$ with
$\Sigma \subset \widehat \Sigma$ and such that the closure $\overline \Sigma$ of~$\Sigma$ 
in~$\widehat \Sigma$ is a compact surface with boundary and $\widehat \Sigma$ is a collar neighbourhood of~$\overline \Sigma$ and a tame Liouville form $\widehat \lambda$ on~$\widehat \Omega\priv \what \Sigma$ with $\widehat \lambda |_{\widehat \Omega \setminus \widehat \Sigma} = \lambda |_{\Omega \setminus \Sigma}$
such that
the boundary $\partial \Omega$ of~$\Omega$ in~$\widehat \Omega$ is a smooth manifold
to which $\widehat \Sigma$ is transverse
(cf.~Figure~\ref{fig-OO}). 
For an extendable polarisation $(\bfsigma,\lambda)$, we denote $\overline \Sigma$ 
the closure of $\Sigma$ in $\what \Sigma$. 
\diam
}
\end{definition}

The completeness assumption on~$X_\lambda$ implies that
$X_{\widehat \lambda}$ is tangent to 
the boundary of $\Omega$ in~$\widehat \Omega$.

In \cite{OS24} we give a large class of examples of extendable Liouville polarisations.
For the moment, we look at the following model example.

\begin{example} \label{ex:ext}
{\rm
Let $(X,\tau)$ be a compact symplectic surface with smooth non-empty boundary, and let $\mu >0$. 
Consider an extension $(\Sigma,\tau)$ of~$(X,\tau)$ to a closed symplectic surface, 
whose area~$\ca_\tau(\Sigma)$ is an integral multiple~$c_1$ of~$\mu$. 
Let $\pi \colon \bigl( \sdb( \Sigma,\tau,c_1), \lambda_0 \bigr) \to \Sigma$ be 
the symplectic disc bundle defined in Section~\ref{s:sdb}, 
together with its Liouville form of residue $-\mu=-\frac{\ca_\tau(\Sigma)}{c_1}$. 
Then $X_{\lambda_0}$ vanishes on the boundary of $\sdb(\Sigma,\tau,c_1)$ by~\eqref{e:Li0}
and is tangent to the fibers.
Hence $\bigl( (\Xcirc,\mu), \lambda_0 \bigr)$ is a Liouville polarisation of 
$$
\sdb( \Xcirc, \tau, \mu) := \pi^{-1}(\Xcirc) \;\subset\; \sdb(\Sigma,\tau,c_1) .
$$ 
The basin of attraction is the whole manifold, and the skeleton is empty.
Taking an open collar neighbourhood of $X \subset \Sigma$ we see that this Liouville 
polarisation is extendable.

From the differentiable view point, $\sdb(\Xcirc,\tau,\mu)$ is just $\Xcirc \times \DDcirc$ 
with a twisted symplectic form,
with the particularity of having a Liouville vector field tangent to the fibers and vanishing at their boundaries. 
\diam
}
\end{example}

The above construction does not rule out the possibility that different extensions 
$(\Sigma,\tau,c_1)$ 
of~$(X,\tau)$ lead to different symplectic manifolds~$\sdb (\Xcirc, \tau, \mu)$. 
This is not so. 
While this fact is not crucial for us, its proof is a good warm-up for the main result 
of this section (Theorem~\ref{t:embliouville}), 
and it leads to two lemmas needed later on.

\begin{lemma} \label{le:A1}
Let $(X,\tau)$ be a compact symplectic surface with boundary, and let $(\widehat X, \widehat \tau)$
be an open collar neighbourhood of~$X$.
Let $W$ be a neighbourhood of~$\widehat X$ in $\SDB(\widehat X, \hat \tau, \mu)$,
and let $\vartheta$ be a bounded closed 1-form on~$W \setminus \widehat X$.
Let $\lambda$ be a Liouville 1-form for $\SDB(\widehat X, \hat \tau, \mu)$ 
as constructed in Example~\ref{ex:ext}. 

\s
{\rm (i)}
There exist neighbourhoods $V,V'$ of~$X$ in~$W$
and a $C^1$-smooth symplectomorphism $f \colon (V,\omega) \to (V',\omega)$
such that $f |_X = \idd$ and $f^* (\lambda + \vartheta) = \lambda$ on~$V \setminus \widehat X$.

\s
{\rm (ii)}
If, moreover, $\vartheta$ vanishes near a union $\partial' X$ of components of the boundary~$\partial X$,
then we can choose $V,V'$ and~$f$ such that $f=\idd$ on a neighbourhood of~$\pi^{-1}(\partial' X) \cap V$.
\end{lemma}

\begin{figure}[h]   
 \begin{center}
  \psfrag{X}{$X$}
	\psfrag{Xh}{$\widehat X$}
	\psfrag{V}{$V$}
	\psfrag{V'}{$V'$}
	\psfrag{f}{$f$}
	\psfrag{W}{$W$}
  \leavevmode\includegraphics{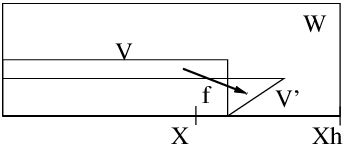}
 \end{center}
 \caption{The symplectomorphism $f \colon V \to V'$ with $f^* (\lambda + \vartheta) = \lambda$.}  \label{fig-VV}
\end{figure}

\proof
This is proven in \cite[Lemma~A.1]{Op13a} in the case that $\widehat X$ is a closed surface.
The same proof yields Lemma~\ref{le:A1}.
\proofend

\begin{proposition} \label{p:unsdbX} 
The symplectic manifold $\sdb(\Xcirc,\tau,\mu)$ does not depend on the choice of 
the extension $(\Sigma,\tau,c_1)$. 
\end{proposition}

\proof 
For $i=0,1$,
let $\bigl( \sdb_i(X,\tau,\mu), \lambda_i \bigr)$ be defined as above by the inclusion 
$\iota_i \colon (X,\tau) \to (\Sigma_i,\tau_i)$.
We need to show that $\sdb_0 ( \Xcirc,\tau,\mu)$ is symplectomorphic to $\sdb_1( \Xcirc,\tau,\mu)$. 
Let $\widehat X_i \subset \Sigma_i$ be open collar neighbourhoods of $X \subset \Sigma_i$.
We can assume that there exists a symplectomorphism $\varphi \colon (\widehat X_0, \tau_0) \to (\widehat X_1, \tau_1)$
that is the identity on~$X$. 
By the symplectic neighbourhood theorem, 
there exist open neighbourhoods~$U_i$ of~$\widehat X_i$ in $\SDB(\widehat X_i,\tau_i,\mu)$
and a symplectic diffeomorphism $\psi \colon U_0 \to U_1$
that extends~$\varphi$. 

\begin{figure}[h]   
 \begin{center}
  \psfrag{X}{$X$}
	\psfrag{X1h}{$\widehat X_0$}
	\psfrag{X2h}{$\widehat X_1$}
	\psfrag{S1}{$\Sigma_0$}
	\psfrag{S2}{$\Sigma_1$}
	\psfrag{f}{$\varphi$}
  \psfrag{p}{$\psi$}
  \psfrag{U0}{$U_0$}
	\psfrag{U1}{$U_1$}
  \leavevmode\includegraphics{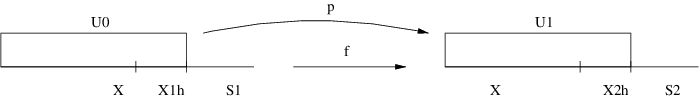}
 \end{center}
 \caption{The extension $\psi \colon U_0 \to U_1$ of $\varphi \colon \widehat X_0 \to \widehat X_1$.}  \label{fig-fp}
\end{figure}

By 
Proposition~\ref{p:singliouville}~(ii) 
there exists a bounded smooth closed one-form~$\vartheta$ on $\im \psi \priv \what X_1$ 
such that $\psi_* \lambda_0-\lambda_1=\vartheta$. 
Let $\widetilde X_1$ be a closed collar neighbourhood of~$X$ that lies in~$\widehat X_1$.
Applying Lemma~\ref{le:A1}~(i) with $(X,\widehat X, W) := (\widetilde X_1, \widehat X_1, \im \psi )$,
we obtain open neighbourhoods $V,V'$ of~$X$ in~$\im \psi$  
and a $C^1$-smooth symplectomorphism~$f \colon V \to V'$ 
that equals the identity on~$\widehat X_1$ and takes 
$\lambda_1 + \vartheta$ to~$\lambda_1$.
Choose $\eps >0$ so small that 
$\psi \bigl(\sdb_{0,\eps}(X,\tau,\mu) \bigr) \subset \dom f$. Then the map 
$$
\varphi := f \circ \psi \colon \sdb_{0,\eps} (X,\tau,\mu) \,\hra\, \sdb(\Sigma_1,\tau_1,c_1^1)
$$  
is well-defined, covers $\iota_1 \circ \iota_0^{-1}$, and verifies $\varphi_*\lambda_0 = \lambda_1$. 
Define 
\fonction{\Phi}{\sdb_0(X,\tau,\mu) }{\sdb_1(X,\tau,\mu)}{x}{ \phi_{X_{\lambda_1}}^{-t} \circ \varphi \circ \phi_{X_{\lambda_0}}^t (x) ,}
where for $x \in X$ we take $t=0$ and for $x \notin X$ we take any $t \geq 0$ such that
$\phi_{X_{\lambda_0}}^t (x) \in \text{Dom}(\varphi) \setminus X$.
Then $\Phi (x)$ is well-defined because $\varphi_* \lambda_0 = \lambda_1$, and 
 a $C^1$-smooth symplectomorphism (see e.g.\ \cite[Section~2.1]{Op13a} if needed). 
By restriction, we obtain a $C^1$-smooth symplectomorphism 
$\sdb_0 (\Xcirc,\tau,\mu) \to \sdb_1(\Xcirc,\tau,\mu)$. 
The claim in the $C^\infty$-category
now follows from the Smoothing Lemma~\ref{le:Zehnder} relegated to the end of this section.
\proofend

Almost the same proof yields the following statement.

\begin{lemma} \label{le:SiOm}
Let $((\Sigma, \mu), \lambda)$ be a smooth and extendable Liouville polarisation 
of~$(\Omega, \omega)$.
Then there exists $\gve >0$ and a $C^1$-smooth symplectic embedding 
$$
\Phi \colon \SDB_{\gve}(\Sigma, \omega|_\Sigma, \mu) \,\to\, \cb (\Sigma, \lambda,\Omega)
$$ 
that is onto a neighbourhood of $\Sigma$ in $\cb (\Sigma, \lambda, \Omega)$ and is
such that $\Phi |_\Sigma = \idd_\Sigma$ and $\Phi^* \lambda = \lambda_0$.
\end{lemma}

\proof 
The previous proof shows the existence of a symplectic embedding
$$
\Phi \colon \sdb_\eps(\Sigma,\om|_{\Sigma},\mu) \,\to\, \cb(\what \Sigma,\hat\lambda,\what \Om),
$$
where $(\what \Sigma,\hat\lambda,\what \Om)$ is an extension of the Liouville 
polarisation $(\Sigma,\lambda)$, with $\Phi|_\Sigma=\id$ and $\Phi^*\widehat \lambda=\lambda_0$. 
The map~$\Phi$ is constructed by a dynamical conjugacy procedure that guarantees that 
$\im \Phi\subset \cb(\Sigma,\widehat \lambda,\what \Om)$. 
Since the Liouville flow associated to~$\widehat \lambda$ is tangent to~$\partial \Om$ 
in a \nbd of $\what \Sigma\cap \partial \Om$, this basin of attraction lies in~$\Om$ and 
coincides with $\cb(\Sigma,\lambda,\Om)$.  
\proofend

Putting Lemmas \ref{le:A1} and \ref{le:SiOm} together, we obtain the following lemma.

\begin{lemma} \label{l:transferliouville2} 
Let $(\bfsigma,\lambda)$ be an extendable polarisation of an exact symplectic manifold~$\Omega$.
Let $(\widehat \Omega, \widehat \Sigma, \widehat \lambda)$ be an extension and let 
$\vartheta$ be a bounded closed $1$-form on a deleted neighbourhood of~$\widehat \Sigma$ in 
$\widehat \Omega \setminus \widehat \Sigma$ that vanishes near~$\sing (\Sigma)$. 
Then there exists a $C^1$-smooth symplectomorphism~$\Phi$ between neighbourhoods of~$\overline \Sigma$ 
in~$\widehat \Omega$  
such that: 
\begin{itemize}
\its $\Phi$ is the identity on $\overline \Sigma$,

\s
\its $\Phi^*(\widehat \lambda + \vartheta) = \widehat \lambda$. 
\end{itemize}
\end{lemma}

\proof
For each $i$ choose small open discs $\coprod_j D_i^j$ around the points in $\Sigma_i \cap \sing (\Sigma)$
whose closures are disjoint from the support of~$\vartheta$,
and take the compact surface $X_i := \widetilde \Sigma_i \setminus \coprod_j D_i^j$,
where $\widetilde \Sigma_i$ is a closed collar neighbourhood of~$\overline \Sigma_i$ 
in~$\widehat \Sigma_i$.
Then $X := \coprod_i X_i$ is a disjoint union, and so is the basin 
$\cb(X,\widehat \lambda) = \sqcup_i \cb(X_i,\widehat \lambda)$.

\begin{figure}[h]   
 \begin{center}
  \psfrag{X1}{$\red{X_1}$}
	\psfrag{X2}{$\red{X_2}$}
	\psfrag{S1}{$\widehat \Sigma_1$}
	\psfrag{S2}{$\widehat \Sigma_2$}
  \psfrag{O}{$\Omega$}
	\psfrag{su}{$\blue{\supp \vartheta}$}
	\leavevmode\includegraphics{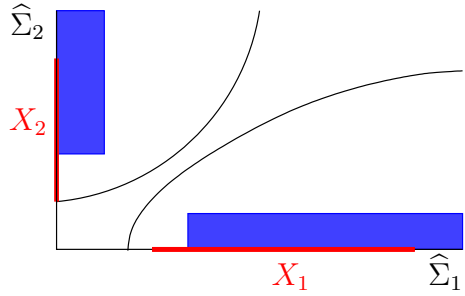}
 \end{center}
 \caption{$\red{X} \subset \widehat \Sigma$ and \blue{$\supp \vartheta$}.}  \label{fig-supp}
\end{figure}

By Lemma~\ref{le:SiOm} we can identify a neighbourhood of~$X_i$ in~$\cb(X_i,\widehat \lambda)$
with~$\SDB_{\gve_i} (X_i,\widehat \lambda)$.
By Lemma~\ref{le:A1}~(ii) there exists a   
$C^1$-smooth symplectomorphism~$\Phi_i$ between neighbourhoods of~$X_i$ in~$\widehat \Omega$ 
such that $\Phi_i^* (\widehat \lambda + \vartheta) = \widehat \lambda$ and $\Phi_i = \id$ 
on~$X_i$ and near the boundaries of the discs~$D_i^j$.
Extend $\Phi_i$ by the identity to an open neighbourhood of all of~$\widetilde \Sigma_i$.
Then the map~$\Phi$ obtained by gluing together the~$\Phi_i$ is as required.
\proofend

Our main interest in the definition of Liouville polarisations is the next result, 
that allows to use the techniques set up to study the symplectic embedding problems 
in~\cite{Op13a,Op13b} in the context of exact symplectic manifolds. 

\begin{definition}
A {\bf symplectic morphism} $\phi \colon \bfsigma \to \bfsigma'$ 
between weighted symplectic divisors
$\bfsigma$, $\bfsigma'$ is an injective continuous map $\phi \colon \Sigma \to \Sigma'$
with the following two properties:
\begin{itemize}
\item[$\bullet$] 
For every component $(\Sigma_i,\mu_i)$ of $\bfsigma$ there exists a component $(\Sigma_{i'}',\mu_{i'})$ of~$\bfsigma'$
with $\mu_i = \mu_{i'}$ such that $\phi |_{\Sigma_i} \colon \Sigma_i \to \Sigma_{i'}'$ is a symplectic embedding.

\s
\item[$\bullet$]
$\phi |_{\Sigma_i} \colon \Sigma_i \to \Sigma_{i'}'$ takes smooth
points to smooth points, i.e., 
$$
\phi \bigl( \Sigma_i \setminus \bigcup_{j \neq i} \Sigma_j \bigr) \,\subset\,
         \Sigma_{i'}' \setminus \bigcup_{j' \neq i'} \Sigma_{j'}' .
$$
\end{itemize}
\end{definition}

%Recall that an embedding~$\gf$ is $(\alpha,\alpha')$-exact if $\gf^* \alpha' = \alpha + df$ 
%for a smooth function~$f$. 

\begin{theorem} \label{t:embliouville} 
Let $(\Om,\om=d\alpha)$ and $(\Om',\om'=d\alpha')$
be exact symplectic mani\-folds with Liouville polarisations 
$(\bfsigma,\lambda)$ and $(\bfsigma',\lambda')$, respectively. 
Assume that $(\bfsigma,\lambda)$ is extendable,
with extension $(\widehat \bfsigma, \widehat \lambda)$.
Assume further that there exists a symplectic morphism 
$\phi \colon \widehat \bfsigma \to \bfsigma'$
that is $(\alpha|_\Sigma, \alpha'|_{\Sigma'})$-exact
(meaning that $\int_\gamma \alpha|_\Sigma = \int_{\phi(\gamma)}\alpha'|_{\Sigma'}$ for all 
piecewise smooth closed curves~$\gamma$ in~$\Sigma$). 

\s
Then there exists an $(\alpha,\alpha')$-exact symplectic embedding
$$
\Phi \colon \cb(\bfsigma,\lambda) = \Om \setminus \skel(\Om,\bfsigma,\lambda) \,\to\, 
          \cb(\bfsigma',
          \lambda') = \Om'\setminus \skel(\Om',\bfsigma',
          \lambda') . 
$$ 
\end{theorem}

\proof
Let $(\widehat \Omega, \widehat \Sigma, \widehat \lambda)$ be the extension of~$(\Omega, \Sigma, \lambda)$.
Let $\widetilde \Sigma$ be a closed collar neighbourhood of~$\overline \Sigma$ in~$\widehat \Sigma$.
By the symplectic neighbourhood theorem, there exists an open neighbourhood~$V$ of~$\widetilde \Sigma$
in~$\widehat \Omega$ and a symplectic extension 
$\psi \colon V \hra \cb(\bfsigma',\lambda') \subset\Om'$ of~$\phi$ 
that preserves the symplectic plumbing structure 
(i.e., the Darboux coordinates $(R_1,\theta_1,R_2,\theta_2)$ near the intersections). 
Since $\psi^*\lambda'$ and~$\widehat \lambda$ are tame and have the same residues, 
$$
\psi^*\lambda' - \widehat \lambda = \vartheta
$$ 
where 
$\vartheta$ is a smooth bounded closed $1$-form on $V \priv \what \Sigma$. 
Since $\psi$ preserves the plumbing structure, we also have 
$\vartheta=0$ on~$\op(\sing \Sigma)$. 

By Lemma~\ref{l:transferliouville2} there exists a 
$C^1$-smooth symplectomorphism~$\Phi_1$ between open neighbourhoods of~$\overline \Sigma$ 
in~$\cb(\widetilde \bfsigma, \hat\lambda)$ 
such that 
$\Phi_1^* \widehat \lambda = \widehat \lambda + \vartheta$ and 
$\Phi_1 |_{\widetilde \Sigma} = \id$. 
Then the map $\Phi := \psi \circ \Phi_1^{-1}$ provides a symplectic extension of~$\phi$ such that 
$\Phi^* \lambda' = \widehat \lambda$.
As in the proof of Proposition~\ref{p:unsdbX} define 
\fonction{\Psi_1}{\cb(\bfsigma,\lambda)}{\cb(\bfsigma',\lambda')}{x}{\phi_{X_{
\lambda'}}^{-t} \circ \Phi \circ \phi_{X_\lambda}^t(x),}
where for $x \in \Sigma$ we take $t=0$ and for $x \notin \Sigma$ we take any $t \geq 0$ such that
$\phi_{X_{\lambda}}^t (x) \in \text{Dom}(\Phi) \setminus \Sigma$.
Then $\Psi_1$ is a $C^1$-smooth symplectic embedding. 
%FFFF
Since $\cb(\bfsigma,\lambda)$ retracts to~$\Sigma$, 
its first homology is generated by cycles in~$\Sigma$, 
on which $\Psi_1 = \phi$. 
Hence, if $\phi$ is $(\alpha |_\Sigma,\alpha'|_{\Sigma'})$-exact, 
$\Psi_1$ is $(\alpha,\alpha')$-exact. 
The claim in the smooth category 
now follows from Lemma~\ref{le:Zehnder}.
\proofend

\begin{lemma} \label{le:Zehnder}
Let $\Psi_1 \colon (M,\omega) \to (M',\omega')$ be a $C^1$-smooth symplectomorphism
of symplectic manifolds without boundary.
Then there exists a $C^\infty$-smooth symplectomorphism $\Psi \colon (M,\omega) \to (M',\omega')$.
Furthermore, if $\Psi_1$ is $(\alpha, \alpha')$-exact, then $\Psi$ can be taken $(\alpha, \alpha')$-exact.
\end{lemma}

\proof
It is shown in \cite[pp.\ 831--836]{Zeh77} that there exists a $C^\infty$-smooth symplectomorphism
$\Psi \colon (M,\omega) \to (M',\omega')$ which is $C^1$-close to~$\Psi_1$.
The proof shows that $\Psi$ is $(\alpha,\alpha')$-exact whenever $\Psi_1$ is $(\alpha,\alpha')$-exact:
Take Darboux charts $\gf_i \colon B^4(a_i) \to U_i \subset M$ 
such that $\gf_i \bigl( B^4(a_i/2) \bigr) =: U_i'$ form a locally finite covering of~$M$.
Smoothing a suitable generating function, one replaces $\Psi_1$ by $\Psi_1^1$
that is $C^\infty$-smooth on~$U_1'$ and agrees with $\Psi_1$ outside~$U_1$.
In fact, convex interpolation of the two generating functions yields a smooth symplectic isotopy $\Psi_1^t$
from $\Psi_1$ to~$\Psi_1^1$ such that $\Psi^{-1}_1 \circ \Psi_1^t$ is supported in~$U_1$.
Since $U_1 = \gf_1(B^4(a_1))$ is simply-connected, this isotopy is Hamiltonian,
and hence $\Psi^{-1}_1 \circ \Psi_1^1$ is $(\alpha,\alpha)$-exact,
see e.g.~\cite[Proposition 9.3.1]{McSa15}.
Since $\Psi_1$ is $(\alpha,\alpha')$-exact, we find that $\Psi^1_1$ is also $(\alpha,\alpha')$-exact.
In the same way, smoothen $\Psi_1^1$ to~$\Psi_1^2$ on~$U_2'$, $\dots$,  
$\Psi_1^k$ to~$\Psi_1^{k+1}$ on~$U_{k+1}'$, $\dots$.
The limit map $\Psi$ is then $(\alpha,\alpha')$-exact. 
\proofend

%%%%%%%%%%%%%%%%%%%%%%%%%%%%%%%%%%%%%%%%%%%%%%%%%%%%%%%%%%%%%%%%%%%%%%%%%%%%%%%%%%%%%%%

\section{Some explicit polarisations and skeleta} \label{s:explicit}

We start our computation of explicit polarisations with the simplest case, in dimension $2$. 
This case is then used in $\S$\2\ref{ss:2x2} to construct Liouville polarisations on bidiscs.
In $\S$\2\ref{ss:surger} we show how certain singular Liouville polarisations can be 
surgered to smooth Liouville polarisations.

\subsection{Polarisations in dimension $2$} \label{ss:dim2}

Let $(S,\omega)$ be a 2-dimensional compact symplectic surface with non-empty boundary~$\partial S$. 
We define a Liouville polarisation of~$S$ to be a finite set of points~$p_i$ in~$\Int S$
with weights $\mu_i >0$
together with a tame Liouville form~$\lambda$ defined on $S \setminus \bigcup_i p_i$  
such that $X_\lambda$ is tangent to $\partial S$.
Tameness now means that for each~$i$ there are symplectic polar coordinates~$(R,\theta)$
near~$p_i$ such that $\lambda = -\mu_i d\theta + \lambda'$ on~$\Op (p_i) \setminus \{p_i\}$,
where $\lambda'$ is a bounded 1-form on~$\Op (p_i) \setminus \{p_i\}$.
One easily checks that $\sum_i \mu_i = \ca_\om (S)$.
Conversely, 
any finite set of positively weighted points $\{(p_i,\mu_i)\}$ on~$(S,\om)$ 
that satisfies $\sum_i \mu_i=\ca_\om(S)$ 
admits a Liouville polarisation, as is easy to see. 
Because of this very direct link between the residues and the area in dimension~$2$, 
we switch notation from $\mu_i$ to~$a_i$ throughout this section. 
Its main purpose is to explain that in dimension~$2$, the skeleton can be prescribed, 
which seems much harder in higher dimensions. 
We recall the notion of regular grid 
in a slightly more general setting.

\begin{definition}\label{d:reggrid} 
{\rm
A {\bf regular grid} $\Gamma$ on a compact surface $S$ with boundary is a finite 
connected 
graph $\Gamma \subset S$ with smooth edges, whose edges cover~$\partial S$ and such that each vertex has a Darboux chart on which $\Gamma$ is a union of radial rays that cut the disc into equal sectors. When the vertex belongs to the interior of $S$, the Darboux chart takes values in $D(\eps)\subset \C$, while it takes values in $D(\eps) \cap \{\im z \geq 0\}\subset \H$ if the vertex belongs to~$\partial S$. 
}
\end{definition}

\begin{proposition}\label{p:poldim1}
Let $(S,\omega)$ be a compact surface with boundary, 
and let $\Gamma$ be a regular grid that decomposes~$S$ 
into $m$~discs of area~$a_i$, in each of which we chose a point~$p_i$. 
Then there exists a Liouville polarisation $\bigl( \{(p_i,a_i)\}, \lambda \bigr)$ of~$S$
whose skeleton is~$\Gamma$. 
This Liouville form can be required to coincide with $(R-a_i)d\theta$ in Darboux coordinates 
near each~$p_i$.  
\end{proposition}

In the case of a closed disc $\overline D(A)$ (which will be the example we will use later on)
one may attempt to prove the proposition as follows:
Just define a first tame Liouville form with residues~$-a_i$ at~$p_i$ by 
$\lambda := \ga_{\st} -\sum a_i \,d\theta_i$, 
where $\theta_i$ is the pull-back by the translation of the vector~$-p_i$ 
of the standard angular coordinate in~$\C \setminus \{0\}$. 
Then the periods of~$\lambda$ vanish on each loop of~$\Gamma$, 
so we can correct~$\lambda$ by adding the differential of a smooth function in order to make 
$\lambda$ vanish on~$\Gamma$. 
Then the Liouville flow of~$\lambda$ fixes~$\Gamma$ pointwise, 
and so $\bigl( \{(p_i,a_i)\}, \lambda \bigr)$ provides a Liouville polarisation of~$\overline D(A)$
whose skeleton contains~$\Gamma$. 
Unfortunately, this proof does not constrain the skeleton enough: 
it could be strictly bigger than~$\Gamma$. 
We therefore proceed in a different manner. 
Recall that $X_\lambda$ is tangent to $\ker \lambda$, 
since $\lambda (X_\lambda) = d\lambda(X_\lambda, X_\lambda) =0$.
Our strategy is therefore to first construct a (singular) foliation~$\cf$ 
to which we wish $X_\lambda$ to be tangent to. 
Once this ``vanishing foliation" is constructed, 
the following lemma will readily provide the Liouville form itself.

\begin{lemma}\label{l:foltoliouville}
Let $\gamma_\theta \colon [0,1] \to \C$, $\theta \in [0,\alpha]$, be a family of smooth rays 
emanating from the origin of~$\C$, 
and set $U := \left\{\gamma_\theta(t) \mid \theta \in [0,\alpha],\, t\in[0,1] \right\}$. 
Fix $a\in \R$. 
We assume that 
\begin{itemize}
\item[$\bul$] 
$(t,\theta) \mapsto \gamma_\theta(t)$ is smooth and a diffeomorphism except at $t=0$.
Thus the set of curves $\{\im \gamma_\theta\}$ provides a foliation $\cf$  of $U\setminus \{0\}$.
 
\item[$\bul$] 
There exists $\eps>0$ such that $\gamma_\theta(t) = t \1 e^{i\theta}$ for $t \leq \eps$. 
\end{itemize}

Then there exists a unique smooth Liouville form $\lambda_a$ on $(U \setminus \{0\},\om_{\st})$ such that
\begin{itemize}
\item[$\bul$] 
$\ker \lambda_a=T\cf$,

\item[$\bul$] 
$\lambda_a(\gamma_\theta(t)) = (R-a) d\theta$ \, for $t$ sufficiently small. 
\end{itemize}
\end{lemma}

\proof
This follows at once from Stokes' theorem.
\proofend

\noindent
{\it Proof of Proposition \ref{p:poldim1}:}  
We write the proof for the disc~$\overline D(A)$. The adaptation to a general surface~$(S,\omega)$
will be clear.
We first set some notation. 
We write $D$ instead of~$\overline D(A)$.
For each point $p_i$ let $\cd_i$ be the closure in~$\RR^2$ of the connected component of~$p_i$ 
in~$D \setminus \Gamma$ (recall that $\Gamma$ contains~$\partial D$ by assumption). 
Write $\sing (\Gamma)$ for the set of non-smooth points of~$\Gamma$, and 
$$
\cs_i := \sing(\Gamma) \cap \cd_i
$$
for the set of singularities of $\Gamma$ in the boundary of~$\cd_i$.
%For each~$p_i$ and a sufficiently small~$\eps$ the translation
For each~$p_i$ we choose a disc~$D(p_i)$ centred at~$p_i$ that is contained in~$\cd_i$.
The translation $(D(p_i),p_i) \to (D(\eps_{p_i}),0)$ to the disc of the same area centred at the origin 
is a Darboux chart that induces coordinates $(R,\theta)$ on~$D(p_i)$. 
Similarly, by our assumption on~$\Gamma$ we can consider for each 
$q \in \sing(\Gamma) \setminus \partial D$ 
symplectic coordinates $\phi_q = (R,\theta) \colon D(q) \to D(\eps_q)$ on a disk around~$q$
such that 
$$
\Gamma \cap D(q) \,=\, \Bigl\{ \theta = \frac k{m_q}, \, k \in [0,m_q-1] \Bigr\},
$$
where $m_q$ is the number of branches of~$\Gamma$ at~$q$. 
The sectors in~$D(q)$ delimited by the rays $\{\theta=\frac k{m_q}\}$ are denoted by 
$$
S_i(q) := \cd_i \cap D(q). 
$$
Also, when $q \in \sing(\Gamma) \cap \partial D$, we have symplectic coordinates~$(R,\theta)$ 
on a \nbd $D^+(q)$
of $q$ in~$D$, with values in $\{ R<\eps_q,\theta \in [0,\pi] \}$
such that 
$$
\Gamma \cap D^+(q) \,=\, \Bigl\{ \theta = \frac k{2m_q}, \, k \in [0,m_q-1] \Bigr\} ,
$$
where $m_q$ is the number of branches of $\Gamma$ at~$q$. 
The sectors  delimited by the rays $\theta=\frac k{2m_q}$ are again denoted by 
$$
S_i(q) := \cd_i \cap D^+(q). 
$$
The radial ray that cuts the sector $S_i(q)$ into two equal sectors is called a local separatrix at~$q$ 
and is denoted~$s_i(q)$. 
We choose the various (half) disks $D(p_i)$, $D(q)$, $D^+(q)$ so small that they are mutually disjoint. 

\begin{figure}[h]   
 \begin{center}
  %\psfrag{De}{$D_\eps(q_1)$}  %FF changé à  D(\eps)
	\psfrag{De}{$D(\eps_{q_1})$}
	\psfrag{Dep1}{$D(p_1)$}
	\psfrag{Dq1}{$D(q_1)$}
  \psfrag{pq1}{$\phi_{q_1}$}
	\psfrag{rq1}{$\rho_{q_1}$}
  \psfrag{D1}{$\cd_1$}
  \psfrag{p1}{$p_1$}
	\psfrag{p2}{$p_2$}
  \psfrag{p3}{$p_3$}
  \psfrag{p4}{$p_4$}
  \psfrag{p5}{$p_5$}
	\psfrag{q1}{$q_1$}
  \psfrag{q2}{$q_2$}
	\psfrag{q6}{$q_6$}
	\psfrag{Sq1}{$S_3(q_1)$}
	\psfrag{Sq6}{$S_5(q_6)$}
  \psfrag{s1}{$s_1(q_1)$}
  \leavevmode\includegraphics{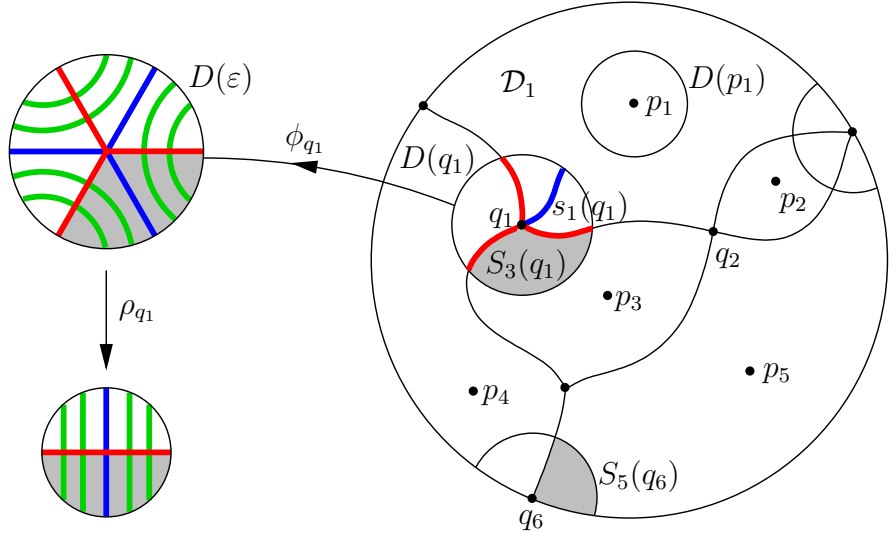}
 \end{center}
 \caption{Notation, and the foliation on $D(\eps_{q_1})$ for $m_{q_1}=3$.}  \label{fig-notation}
\end{figure}

Finally, we fix a symplectic diffeomorphism $f_i \colon \cd_i \setminus \cs_i \to \overline D(a_i) \setminus \cp_i$ 
where $\cp_i$ is a finite set of points, 
and set $\lambda_i^f :=f_i^* \ga_{\st}$. 
Then $\lambda_i^f$ is positive on $\partial \cd_i \setminus \cs_i$. 
The notations are illustrated on the right part of Figure~\ref{fig-notation}.

\m \ni
{\bf Step 1: Constructing the vanishing foliation near the $p_i$ and $\sing(\Gamma)$ 
(see Figure~\ref{fig-notation})
}  

\s \ni
In $D(p_i)$, we fix $\cf$ to be the radial foliation $\{\theta = \const \}$.

\s
For a point $q \in \sing(\Gamma) \setminus \partial D$, let 
$\rho_q \colon D(\eps_q) \to \C$ be the ``ramified covering"  
$$
\rho_q (R,\theta) \,:=\, \left( \frac {2R}{m_q},\frac{m_q}2 \theta \right).
$$  
Note that $\rho_q$ is smooth and symplectic on 
    $D(\eps_q) \setminus \{0\}$, 
but only continuous at~$0$. 
The sectors $S_i(q)$ around $q$ are the preimages under $\rho_q \circ \phi_q$ of $\{ y \geq 0\}$ or $\{y \leq 0\}$, 
and the separatrices~$s_i(q)$ are sent by~$\rho_q \circ \phi_q$ to~$\{x=0\}$. 
We define the foliation in~$D(q)$ to be the pull-back by~$\rho_q \circ \phi_q$ of the vertical foliation:
$$
\cf |_{D(q)} \,:=\, (\rho_q \circ \phi_q)^* \{x=c\}.
$$
This local foliation has remarkable properties:

\begin{enumerate}
\item[(i)]
The local separatrices emanating from $q$ are leaves of $\cf$; 
they are the pull-back by~$\rho_q \circ \phi_q$ of the curves~$\{x=0\}$.
 
\item[(ii)]
Since $\rho_q$ sends concentric circles around the origin to concentric circles,
and since $\rho_q \circ \phi_q$ is a diffeomorphism between $S_i(q) \setminus \{q\}$ and $\{y \geq 0\} \setminus\{0\}$ or $\{ y \leq 0\} \setminus \{0\}$, 
the leaves $\{x=c\}$ with $c \neq 0$ are pulled back by~$\rho_q \circ \phi_q$ to smooth curves that intersect $\Gamma$  
orthogonally (in the chart $\phi_q$). 

\item[(iii)] 
As a result, the foliation $\cf$ can be extended to a smooth foliation on $\op(\Gamma) \priv \sing (\Gamma)$ 
that is transverse to the regular part of $\Gamma$. 
\end{enumerate}

\begin{figure}[h]   
 \begin{center}
  \psfrag{pq1}{$\phi_{q_1}$}
	\psfrag{G}{$\Gamma$}
  \psfrag{D1}{$\cd_1$}
  \psfrag{p1}{$p_1$}
  \psfrag{p3}{$p_3$}
	\psfrag{q1}{$q_1$}
  \psfrag{q2}{$q_2$}
	\psfrag{q3}{$q_3$}
	\psfrag{q0}{$q_0$}
	\psfrag{G}{\textcolor{red}{$\Gamma$}}
  \leavevmode\includegraphics{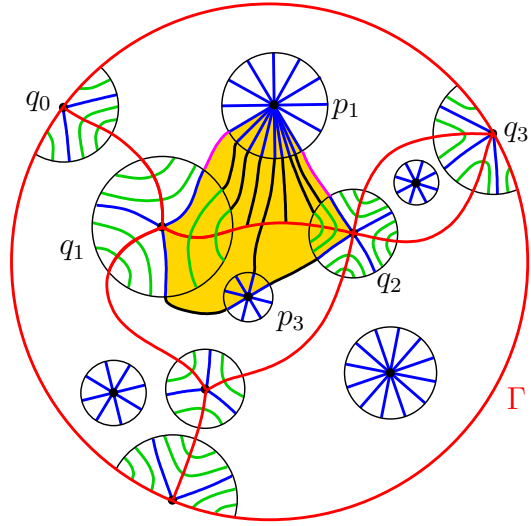}
 \end{center}
 \caption{The foliation $\cf$ on the $D(p_i)$, $D(q_j)$, and $D^+(q_k)$ (blue and green), 
          part of the foliation~$\cf_1$, and $\cd(p_1,p_3)$ (yellow).}  \label{fig-fol}
\end{figure}

Finally, for a point $q \in \sing(\Gamma) \cap \partial D$, the model is the same. 
We leave it to the reader to complete the picture, cf.\ Figure~\ref{fig-fol}.

\m \ni
{\bf Step 2: Interpolating the separatrices (see Figure \ref{fig-fol})}

\s \ni 
Fix $i \in \{1, \dots, m\}$.
(In Figure~\ref{fig-fol}, $i=1$.)
Fix a reference point $q_0 \in \cs_i$, denote $q_1, \dots, q_\ell$ the points of~$\cs_i$ enumerated by 
going around~$\partial \cd_i$ in the anti-clockwise sense from~$q_0$, and let $\partial \cd_i(q_j \to q_{j+1})$ 
be the smooth arc in~$\Gamma$ joining $q_j$ to~$q_{j+1}$ on~$\partial \cd_i$. 
Now interpolate between the local separatrix $s_i(q_0)$ and the ray $\{\theta=0\}$ in~$D(p_i)$ 
by a smoothly embedded curve. 
Set 
$$
\beta_j \,:=\, \frac 1{a_i} \int_{\partial \cd_i(q_j \to q_{j+1})} \lambda_i^f, \qquad j=0, \dots , |\cs_i|-1, 
$$
and 
$$
\theta_j \,:=\, \sum_{k=0}^{j-1} \beta_k.
$$
Interpolate now inductively between the local separatrix $s_i(q_j)$ and the ray 
$\{\theta=\theta_j\}$ 
by disjoint smoothly embedded curves (for different $j$) in~$\cd_i$ such that the area enclosed by 
the three arcs
$\partial \cd_i(q_j \to q_{j+1})$, $s(p_i\to q_j)$, and $s(p_i\to q_{j+1})$ 
is~$\beta_j \1 a_i$.  
Here, $s(p_i\to q_j)$, $j=1,\dots, \ell$, denote the global separatrices thus constructed. 
These interpolations can be found because the total area enclosed by~$\partial \cd_i$ is~$a_i$.  
Interpolations $s(p_1 \to q_1)$ and $s(p_1 \to q_2)$ are drawn in pink in Figure~\ref{fig-fol}.
We declare the $s(p_i\to q_j)$ to be leaves of~$\cf$ in~$\cd_i$.

\m \ni
{\bf Step 3: Interpolating between the local foliations in $\cd_i$} 

\s \ni
The foliation~$\cf_i$ on~$\cd_i$ is then taken to be any (singular) foliation in~$\cd_i$ that verifies the following properties:

\begin{enumerate}
\item The points $p_i$ and $q_j \in \cs_i$ are the only singularities of $\cf |_{\cd_i}$, 
so $\cf |_{\cd_i \setminus \{p_i\}}$ is smooth.

\item The curves $s(p_i\to q_j)$ are leaves.

\item 
On the $D(p_i)$, $D(q_j)$, $D^+(q_k)$, and $\op(\Gamma) \priv \sing (\Gamma)$
the leaves of~$\cf$ are the ones described in Step~1. 
 
\item All leaves join $p_i$ to $\partial \cd_i$ and are transverse to 
$\partial \cd_i \setminus \cs_i$.

\item For $c \in (\beta_j, \beta_{j+1})$ the leaf extending the ray 
$\{\theta=c\} \subset D(p_i)$ joins 
$p_i$ to the unique point $q \in \partial \cd_i (q_j\to q_{j+1})$ such that 
$$
\int_{\partial \cd_i(q_j \to q_{j+1})} \lambda_i^f \,=\, c-\beta_j . 
$$
\end{enumerate}
Figure~\ref{fig-fol} shows in black seven curves that belong to leaves of~$\cf_1$
from $p_1$ to~$\pp \cd_1(q_1 \to q_2)$.
By (3) and~(4) we can successively construct foliations $\cf_1, \cf_2, \dots, \cf_m$ such that all
non-singular leaves are $C^\infty$-smooth.

\medskip \noindent
{\bf Step 4: Joining the foliations $\cf_i$} 

\s \ni
Since $\cf_i$ is a smooth foliation on~$\cd_i \setminus \bigl( \{p_i\} \cup \cs_i \bigr)$ 
and by the smooth fitting of the leaves across the~$\partial \cd_i(q_j \to q_{j+1})$,
the foliation~$\cf$ on $D \setminus \bigl( \bigcup \{p_i\} \cup \sing(\Gamma) \bigr)$ 
defined by 
$$
\cf |_{\cd_i \setminus \left( \{p_i\} \cup \cs_i \right)} \,:=\, \cf_i
$$
is $C^\infty$-smooth, 
and has a unique extension to a singular foliation on the whole disc~$D$ whose singularities 
are $\bigcup \{p_i\} \cup \sing(\Gamma)$. 
The different separatrices $s(p_i \to q)$, $q \in \cs_i$, are smooth leaves. 
The leaves of~$\cf$ that pass through a point in $\Gamma \setminus \sing(\Gamma)$ 
that separates between $\cd_i$ and~$\cd_j$ join $p_i$ to~$p_j$. 
Finally, the leaves through a sequence of regular points on~$\cd_i \cap \cd_j$ 
that converge to a singular point $q \in \sing(\Gamma)$ break into the two separatrices 
$s(p_i \to q)$ and~$s(p_j\to q)$. 

We denote by $\cd(p_i,p_j)$ the closure of the union of leaves joining $p_i$ to~$p_j$. 
Figure~\ref{fig-fol} shows $\cd(p_1,p_3)$ in yellow.

\m \ni
{\bf Step 5: Fixing the Liouville form itself}

\s \ni 
By applying Lemma~\ref{l:foltoliouville} twice,
first to the curves of the foliation~$\cf$ of the interior~$\cdcirc (p_i,p_j)$ 
that emanate from $p_i$ with $a:=a_i$,
and then to the same set of curves but seen as emanating from~$p_j$ with $a:=a_j$,   
we obtain on~$\cdcirc (p_i,p_j)$ 
two Liouville forms $\lambda_i, \lambda_j$.
We claim that they coincide. To see this, we first check that these forms vanish on~$\Gamma$. 
Indeed, let $q \in \Gamma \cap \cdcirc (p_i,p_j)$. Since $\lambda_i$ vanishes along~$\cf$, 
we already have $\lambda_i(q)v=0$ for $v \in T_q \gamma$, where $\gamma$ is the leaf through~$q$. 
On the other hand, parametrise the leaves through~$q \in \Gamma \cap \cdcirc (p_i,p_j)$ by curves $\gamma(\theta,t)$, 
where $\theta$ represents the angle of the leaf at~$p_i$ 
and the parameter~$t$ verifies $\gamma_\theta(1) \in \Gamma$, and put 
$\frac \partial {\partial \theta}(q) := \frac {\partial \gamma} {\partial \theta}(\theta(q),1)\in T_q \Gamma$. 
Since $\lambda_i$ vanishes along~$\cf$, and by Stokes' theorem and the second property of~$\lambda_i$ 
from~Lemma~\ref{l:foltoliouville}, 
\begin{equation} \label{e:van}
\lambda_i(q) \, \frac \partial {\partial \theta}(q) \,=\, \frac {d\ca_i}{d\theta}(\theta(q))-a_i,
\end{equation}
where $\ca_i(\theta)$ is the area of the disc bounded by 
$\Gamma$, $\gamma(\theta,[0,1])$, and $\gamma(\theta(q),[0,1])$. 
By the choice of our foliation (point~$5$ in Step~3), the right hand side of~\eqref{e:van} vanishes, 
so $\lambda_i(q)$ vanishes. 
Now, since $\lambda_i$ and~$\lambda_j$ coincide on~$\cf$ and on~$\Gamma$, 
Stokes' theorem shows that they coincide on~$\cdcirc(p_i,p_j)$. 
We therefore define 
$$
\lambda |_{\cdcirc(p_i,p_j)} \,:=\, \lambda_i = \lambda_j. 
$$
Since $\lambda_i$ is defined and smooth on $\cdcirc_i \setminus \{p_i\}$, 
and not only on~$\cdcirc(p_i,p_j)$, 
we can glue the $\lambda_i$ to get a smooth Liouville form~$\lambda$ on $D \setminus \sing(\Gamma)$. 
Notice that for $\gamma(\theta,t)$ as defined above, we have by the same argument
$$
\lambda(\gamma(\theta,t)) \, \frac \partial {\partial \theta} < 0 \qquad \forall \, t<1 . 
$$
Hence on $\cdcirc_i \setminus \{p_i\}$ we have 
$d \lambda (X_\lambda, \frac{\partial}{\partial \theta}) = \lambda (\frac{\partial}{\partial \theta}) < 0$,
which shows that
the coordinate $t$ is (negative) gradient-like for the Liouville flow $X_\lambda$ of~$\lambda$, 
and so $\cb(p_i,\lambda) \supset \cd_i$. 
On the other hand, $\lambda$ vanishes on $\Gamma$, so $\Gamma$ is invariant under the Liouville flow of~$\lambda$. Summarizing, we have therefore proved that 
\begin{itemize}
\item[$\bul$] 
$\lambda$ is a smooth Liouville form on $D \setminus \bigl( \sing \Gamma \cup \{p_1,\dots,p_m \} \bigr)$,

\item[$\bul$] 
$\lambda$ is tame at each $p_i$ with residue $-a_i$, in fact $\lambda = (R-a_i) d\theta$ near~$p_i$,

\item[$\bul$] 
$X_\lambda$ vanishes on $\Gamma\priv \sing \Gamma$,

\item[$\bul$] 
every trajectory starting in $D \setminus \Gamma$ converges to some $p_i$. 
\end{itemize}

\m \ni
{\bf Step 6: Smoothening $\lambda$ at $\sing(\Gamma)$}

\s \ni 
We only treat the case of $q \in \sing(\Gamma) \setminus \partial D$, the other case being similar.  
We recall that a neighbourhood $D(q)$ of~$q$, that is identified with~$D(\eps_q)$, 
covers a neighbourhood of $0 \in \C$ 
via the map 
$\rho_q (R,\theta) = (\frac{2R}{m_q}, \frac{m_q}2 \theta)$, 
and that in these coordinates, 
$\cf$ is the foliation $\rho_q^* \{ x = \text{const}\}$. 
This foliation is thus also tangent to the kernel of the local Liouville form 
$\lambda_q := \rho_q^* (-\tfrac 1\pi \2 y \1 dx)$.
Away from $q$, this form~$\lambda_q$ also vanishes on $\rho_q^* \{y=0\} = \Gamma \cap D(q)$. 
Away from $q$, both forms~$\lambda$ and~$\lambda_q$ are thus smooth and vanish along~$\cf$ and on~$\Gamma$. 
Stokes' theorem therefore guarantees that they coincide on $D(q) \priv \{q\}$. 
None of the two forms smoothly extends to~$q$, however, since at $(x,y) \neq (0,0)$,
\begin{eqnarray*}
\lambda_q &=& \rho_q^* (-\tfrac 1\pi \2 y \1 dx) \\
 &=& \rho_q^* \big(\tfrac1{2\pi} (x \1 dy - y \1 dx) - \tfrac 1{2\pi} (x \1 dy + y \1 dx) \big) \\
 &=& \rho_q^* \big( R \2 d \theta - \tfrac 1{2\pi} d(xy) \big) \\
 &=& \rho_q^* \big( R \2 d \theta - \tfrac 1{2\pi} d(R \sin (2\pi \theta) \cos (2\pi\theta)) \big) \\
 &=& R \2 d \theta - \tfrac 1{4\pi} \rho_q^*d ( R \sin( 4\pi \theta) ) \\
 &=& R \2 d \theta - \tfrac 1{2\pi m_q} d (R \sin(2\pi m_q \theta) )
\end{eqnarray*}
and the function $R \sin (2\pi m_q \theta)$ is not smooth at the origin for $m_q \geq 2$.
We thus alter $\lambda_q$ in~$D(q)$ 
so as to make it smooth at~$q$. 

The function $r^{m_q} \sin (2 \pi m_q \theta)$ is smooth on~$\RR^2$, 
since it is the imaginary part of $z \mapsto z^{m_q}$.
Let $\chi \colon [0,2\eps_q] \to \R_{\geq 0}$ be a function that coincides with~$r^{m_q}$ near~$0$
and with~$R$ on~$[\eps_q,2\eps_q]$, is positive and smooth except at~$0$, 
and meets $\chi(R) < R$ on~$(0,\eps_q)$. 
Define the smooth 1-form~$\lambda'$ on~$D$ by $\lambda' = \lambda$ on~$D \setminus \bigcup_{q \in \sing \Gamma} D(q)$
and
$$
\lambda' \,=\,  
R \2 d \theta -\frac 1{2\pi m_q}  d \bigl( \chi(R) \sin(2 \pi m_q \theta) \bigr)  
      \; \text{ on }\, D(q). 
$$
Then $\lambda'$ is a smooth Liouville form on $D \setminus \{p_1, \dots, p_m \}$, 
tame at~$p_i$ with residue~$-a_i$. 
We claim that the skeleton of $\lambda'$ is $\Gamma$, as required. 
To see this, notice that $X_{\lambda'} = X_{\lambda}$ on 
$D \setminus \bigcup_{q \in \sing \Gamma} D(q)$, and 
$$
X_{\lambda'} =(R-\chi(R) \cos (2 \pi m_q \theta)) \frac \partial{\partial R} + \frac 1{2\pi m_q} \chi'(R) \sin(2 \pi m_q \theta) \frac\partial{\partial \theta} 
$$
on $D(q)$. 
We first look at the $\frac{\pp}{\pp \theta}$-component of~$X_{\lambda'}$. 
Since 
$$
\Gamma \,=\, \{\sin (\pi m_q \theta) = 0 \} \,\subset\, \{\sin (2 \pi m_q \theta) = 0 \} ,
$$ 
$X_{\lambda'}$ is radial on $\Gamma \cap D(q)$ (hence its flow preserves $\Gamma$), 
and the trajectories of the points in $D(q) \setminus \Gamma$ flow away from~$\Gamma$. 
From the $\frac{\pp}{\pp R}$-component, 
and from 
$\Gamma \supset \{ \cos (2\pi m_q \theta) =1 \}$ and $\chi(R) < R$
on~$(0,\eps_q)$, we infer that the trajectories of the points in $D(q) \setminus \Gamma$ 
all leave~$D(q)$, 
and hence reach the set $\{\lambda=\lambda'\} \setminus \Gamma$ (since $\Gamma$ is invariant). 
It follows that $\skel(\lambda')=\Gamma$.
\proofend

The following easy result will be needed in the next subsection.

\begin{lemma} \label{l:doublingdim2} 
Let $(S,\om)$ be a compact symplectic surface with or without boundary. 
Let $\bigl( \cp := \{ (p_i,a_i), \; i=1,\dots, n \},\lambda \bigr)$ 
with $p_i \in S \setminus \partial S$
be a polarisation of~$S$ 
(where $\lambda$ is a tame Liouville form with residues~$a_i$ at~$p_i$ 
and with flow tangent to~$\partial S$),
and set $\Gamma := \skel(S,\cp,\lambda)$. 

For $k \leq n$ and any collection $\{(a_i^1,a_i^2), i=1,\dots,k\}$ of pairs of 
positive real numbers with $a_i=a_i^1+a_i^2$ 
there exist:
\begin{itemize}
\item[$\bul$] 
a collection of open embedded discs~$D_i$ around~$p_i$,
in any prescribed neighbourhood of~$p_i$

\s
\item[$\bul$] 
for $i \leq k$:
points $p_i^1,p_i^2 \in D_i$

\s
\item[$\bul$] 
and a tame Liouville form~$\lambda'$ with residues~$a_i^1$ at~$p_i^1$ and~$a_i^2$ at~$p_i^2$ for $i \leq k$ and $a_i$ at~$p_i$ for $i \in \{k+1, \dots, n\}$
\end{itemize}
such that 
$$
\cb(D_i,\lambda') \,:=\, 
\left\{ x \in S \mid \exists \, t>0 \mbox{ such that } \phi_{\lambda'}^t(x) \in D_i \right\}
\,=\, \cb(p_i,\lambda). 
$$
\end{lemma}

\proof
Since $\lambda$ is tame at $p_i$, its dual vector field points towards~$p_i$ on a neighbourhood of~$p_i$, 
so there exist arbitrarily small disjoint open discs~$D_i \subset S \setminus (\partial S \cup \Gamma)$ around~$p_i$ with 
$$
\cb(D_i,\lambda)=\cb(p_i,\lambda). 
$$
For $i \leq k$ take any pair of points
$p_i^j \in D_i$, $j=1,2$, and a  
Hamiltonian diffeomorphisms $\phi_i^j$ with compact support in~$D_i$ 
such that $\phi_i^j(p_i^j) = p_i$.
Define 
$$
\lambda'(x) \,:=\, 
\left\{
\begin{array}{ll}
\lambda(x) & \text{ if } x \in S \setminus \bigcup_{i=1}^k D_i,\\ [0.2em]
\frac {a_i^1}{a_i} \, \phi_i^{1*} \lambda(x) + \frac {a_i^2}{a_i} \, \phi_i^{2*} \lambda (x) & \text{ if } x \in D_i \mbox{ for some $i \in \{1, \dots, k\}$}.
\end{array}
\right.
$$
Then $\lambda'$ is a Liouville form tame at $p_i^j$ with residue $a_i^j$. 
Since it coincides with~$\lambda$ on~$\Gamma \cup \partial S$, its flow preserves both~$\Gamma$ and~$\partial S$, 
so $\left( \cp' := \left\{ (p_i^j,a_i^j) \right\}, \lambda' \right)$ is a Liouville polarisation of~$S$. 
Since $\lambda'$ coincides with~$\lambda$ outside the discs~$D_i$, 
$$
\cb(D_i,\lambda') = \cb(D_i,\lambda) = \cb(p_i,\lambda) .
$$ 
The lemma follows. 
\proofend

\subsection{Bidiscs} \label{ss:2x2} 

As before we denote by $\overline D(A)$ the closed disc of area $A$, and by $D(A)$
its interior. 

\begin{lemma} \label{l:poldim2} 
Let $\Gamma_1 \subset D(A)$ and $\Gamma_2 \subset D(B)$ be two regular grids containing $\partial D(A), \partial D(B)$ whose complements are a union of discs of area $a_i, b_j$ on which we choose points $p_i, q_j$, respectively. 
Then 
$$
\bfsigma := 
\left\{ \bigl( {p_i} \times D (B), a_i \bigr), \bigl( D (A) \times {q_j}, b_j \bigr) \right\}
$$
is an extendable Liouville polarisation of $D (A) \times D (B)$, 
and there is a tame Liouville form~$\lambda$ 
whose skeleton is~$\Gamma_1 \times \Gamma_2$.  
\end{lemma}

\proof
By Proposition~\ref{p:poldim1} there exist Liouville forms $\lambda_1, \lambda_2$ on $D(A), D(B)$ that are tame along 
the $p_i, q_j$, with residues $a_i, b_j$ and skeleton $\Gamma_1, \Gamma_2$, respectively. 
Then the restriction of the form $\lambda := \pi_1^*\lambda_1 \times \pi_2^* \lambda_2$ 
to~$\bfsigma$ is a Liouville form tame along~$\bfsigma$ 
with the correct residues. 
Since the associated Liouville flow is simply the product of the flows in each factor, it clearly preserves 
$$
\partial \bigl( D(A) \times D(B) \bigr) \,=\, 
\bigl( \partial D(A) \bigr) \times D(B) \cup D(A) \times \bigl( \partial D(B) \bigr) .
$$ 
Moreover, for $x \notin \Gamma_1 \times \Gamma_2$, we have $\pi_1(x) \in D(A) \setminus \Gamma_1$ or 
$\pi_2(x) \in D(B) \setminus \Gamma_2$, so the trajectory of $\pi_1(x)$ or~$\pi_2(x)$ under the Liouville flow 
of~$\lambda_1$ or~$\lambda_2$ is forward attracted by one of the~$p_i$ or~$q_j$, 
which shows that the complement of~$\Gamma_1 \times \Gamma_2$ in $D (A) \times D (B)$ is forward attracted by~$\bfsigma$, 
whence $\skel(\bfsigma,\lambda) \subset \Gamma_1 \times \Gamma_2$. 
Finally, for a point in $\Gamma_1 \times \Gamma_2$, 
both components remain in $\Gamma_1, \Gamma_2$, 
so $\Gamma_1 \times \Gamma_2$ is invariant under the Liouville flow and 
$\skel(\bfsigma,\lambda) \supset \Gamma_1 \times \Gamma_2$. 
The extendability is obvious: take $\what \Om:=D(A+\eps) \times D(B+\eps)$
and $\what \Sigma := 
\left\{  p_i \times D(B+\eps),  D(A+\eps) \times q_j \right\}$.
\proofend

Combining Lemma~\ref{l:doublingdim2} and Lemma~\ref{l:poldim2} we obtain: 

\begin{lemma} \label{l:double} 
Let $\bigl( \left\{(p_i,a_i) \right\}_{i=1,\dots,m}, \lambda_A \bigr)$ and 
    $\bigl( \left\{(q_j,b_j) \right\}_{j=1,\dots,n}, \lambda_B \bigr)$ 
be Liouville polarisations of $D(A)$ and~$D(B)$ 
(so $\sum a_i=A$ and $\sum b_j=B$) 
with skeleton $\Gamma_A$ and $\Gamma_B$, respectively. 
Let $m'\leq m$, $n' \leq n$, and assume that for $i \leq m'$ and $j \leq n'$
we are given decompositions
$$
a_i=a_i^1+a_i^2 \quad \text{and} \quad b_j=b_j^1+b_j^2,
\quad \mbox{ with }\; a_i^\ell, b_j^\ell > 0 .
$$

Then there exist:
\begin{itemize}
\item[$\bullet$] 
open discs~$D_i^A$, $D_j^B$ for $i \leq m$, $j \leq n$
around~$p_i$, $q_j$,
which lie in any prescribed neighbourhood of~$p_i$, $q_j$

\s
\item[$\bullet$] 
for $i \leq m'$, $j \leq n'$:
points $p_i^1,p_i^2 \in D_i^A$ and $q_j^1,q_j^2 \in D_j^B$

\s
\item[$\bullet$] 
and a tame Liouville form~$\lambda'$ on the complement of  
\begin{eqnarray*}
\qquad \quad
\bfsigma' &\!\!:=\!\!&
\bigl\{ 
\bigl( p_i^\ell \times D(B),a_i^\ell \bigr)_{i \leq m'}^{\ell =1,2}, \,
\bigl( D(A)\times q_j^\ell,b_j^\ell \bigr)_{j \leq n'}^{\ell =1,2}, \\
&&
\phantom{ \bigl\{ }
\bigl( p_i \times D(B),a_i \bigr)_{i > m'}, \;
\bigl( D(A) \times q_j,b_j \bigr)_{j > n'} 
\bigr\}
\end{eqnarray*}
that makes $(\bfsigma',\lambda')$ a Liouville polarisation of $D (A) \times D(B)$ 
\end{itemize}
such that 
$$
\cb \Bigl( \bigcup_i D_i^A \times D(B) \cup \bigcup_j D(A) \times D_j^B,\lambda' \Bigr)
\,=\, \bigl( D(A) \times D(B) \bigr) \setminus (\Gamma_A \times \Gamma_B) .
$$
\end{lemma}

\proof
Take small open discs $D_i^A$ around $p_i$, for $i \leq m'$ take points $p_i^1,p_i^2 \in D_i^A$, 
and let $\lambda_A'$  be the tame Liouville form on
$$
D(A) \setminus \left( 
\{ p_i^\ell \}_{i \leq m'}^{\ell =1,2} 
\cup
\{ p_i \}_{i > m'} 
\right)
$$
provided by Lemma~\ref{l:doublingdim2}, 
with residues $a_i^\ell$ at~$p_i^\ell$ for $i \leq m'$ 
and $a_i$ at~$p_i$ for $i>m'$, 
and such that 
$$
\cb \Bigl( \bigcup_i D_i^A,\lambda_A' \Bigr) \,=\, D(A) \priv \Gamma_A .
$$
Construct similarly  $\lambda_B$ on 
$$
D(B) \setminus \left( 
\{ q_j^\ell \}_{i \leq n'}^{\ell =1,2} 
\cup
\{ q_j \}_{j > n'} 
\right) .
$$
As we have seen in the proof of Lemma~\ref{l:poldim2},
the 1-form $\lambda' := \lambda_A' \oplus \lambda_B'$ is tame along $\bfsigma'$, 
and 
$$
\cb \Bigl( \bigcup_i D_i^A \times D(B) \cup \bigcup_j D(A) \times D_j^B,\lambda' \Bigr)
\,=\, \bigl( D(A) \times D(B) \bigr) \setminus (\Gamma_A \times \Gamma_B) .
$$
The lemma is proved.
\proofend

\subsection{Surgery on Liouville polarisations}
\label{ss:surger}
The aim of this section is Proposition~\ref{p:surgerform} that allows to glue certain components of a Liouville polarisation, to obtain a ``smoother" Liouville polarisation. 
We will need the following rather obvious statement.

\begin{lemma} \label{l:localtame}
Let $\Sigma \subset B^4(1)$ be a smooth symplectic curve that near the boundary of~$B^4(1)$
agrees with a collection of complex lines (1-dimensional complex vector spaces). 
Let $\mu \in \R$. 
Then there exists a Liouville form~$\lambda$ on $B^4(1) \priv \Sigma$, tame along~$\Sigma$, 
with residue~$\mu$. 
(We do not claim that $(B^4(1),\Sigma,\lambda)$ is a Liouville polarisation.) 
\end{lemma}

\proof
Since $\Sigma$ coincides with a union of complex lines near $\partial B^4(1)$,  
compactifying $B^4(1)$ into~$\CP^2(1)$ provides a smooth symplectic curve $\Sigma'$ in~$\CP^2$. 
This is a polarisation, so by Proposition~\ref{p:singliouville}~(i)
there exists a Liouville form~$\lambda'$ on~$\CP^2 \priv \Sigma'$ that is tame along~$\Sigma'$. Restricting $\lambda'$ to $B^4(1) \subset \CP^2(1)$ provides a Liouville form~$\lambda$ on 
$B^4(1) \priv \Sigma$ tame along~$\Sigma$.  
It has a well-defined residue which is a constant $\nu$ because $\Sigma$ is assumed to be smooth. Then the form 
$$
\lambda_{\st}+\frac \mu\nu (\lambda-\lambda_{\st})
$$
is a Liouville form on $B^4(1) \priv \Sigma$, 
tame along~$\Sigma$ because it coincides with a multiple of~$\lambda$ up to a smooth form, 
with residue~$\mu$.
\proofend

\begin{proposition}\label{p:surgerform}
Let $(M,\om,\bfsigma,\lambda)$ be a polarised symplectic manifold, either closed or Liouville polarised. 
Assume that the weights $\mu_i,\mu_j$ of two components~$\Sigma_i,\Sigma_j$ of~$\bfsigma$ coincide. 
Then for any neighbourhood~$U$ of $\Sigma_i \cap \Sigma_j$ 
there exist:
\begin{itemize}
\item[$\bullet$] 
A smooth symplectic curve $\Sigma_{ij}$ that coincides with $\Sigma_i\cup \Sigma_j$ 
in the complement of $U$. We then define 
$$
\bfsigma':=\{(\Sigma_k,\mu_k'=\mu_k)_{k\neq i,j}, (\Sigma_{ij},\mu_{ij}'=\mu_i=\mu_j)\}. 
$$
\item[$\bullet$] 
A Liouville form $\lambda'$ on $M \priv \Sigma'$ tame along $\Sigma'$ with residues $\mu_k'$ that coincides with~$\lambda$ on the complement of~$U$ and is such that 
$(M,\om,\bfsigma',\lambda')$ is again a (Liouville) polarised symplectic manifold. 
\end{itemize}
\end{proposition}

\begin{figure}[h]   
 \begin{center}
  \psfrag{e}{$\eps$}
	\psfrag{2e}{$2\eps$}
  \psfrag{S1}{$\Sigma_1$}
	\psfrag{S2}{$\Sigma_2$}
 \leavevmode\includegraphics{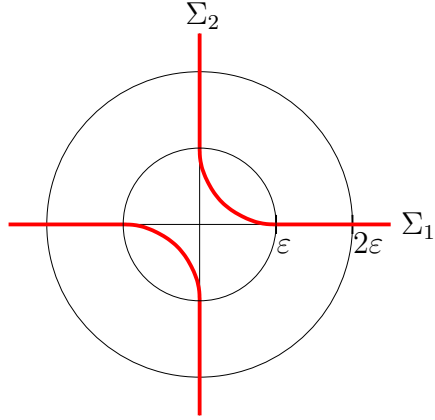}
 \end{center}
 \caption{From $\Sigma_1 \cup \Sigma_2$ to $Z_\gd$, schematically.}  \label{fig:interZ}
\end{figure}

\proof
We first consider a model in $\RR^4$.
Let $\Sigma_1 = \{z_2=0\}$ and $\Sigma_2 = \{z_1=0\}$ be the complex coordinate lines in~$\CC^2$.
For $\gve >0$ let $\chi \colon [0,2\gve] \to [0,1]$ be a smooth function that is~$1$ near~$0$
and vanishes on $[\gve, 2\gve]$.
For $\gd >0$ define the smooth cylinder
$$
Z_\delta := \{ z_1 z_2 = \delta \1 \chi(R) \} .
$$
Note that $Z_\gd = \Sigma_1 \cup \Sigma_2$ outside $B^4(\gve)$, see Figure~\ref{fig:interZ}.
As is well known, $Z_\gd$ is symplectic for $\delta = \gd (\gve)$ small enough.

Consider the Liouville form 
$$
\gl_0 \,=\, \gl_{\st} - \mu (d\theta_1+d\theta_2) 
$$
on $\RR^4 \setminus (\Sigma_1 \cup \Sigma_2)$, 
tame along $\Sigma_1 \cup \Sigma_2$ with residue~$-\mu$.
Since $Z_\gd = \Sigma_1 \cup \Sigma_2$ on $B^4(2\gve) \setminus B^4(\gve)$, 
Lemma~\ref{l:localtame} provides a tame Liouville form $\gl_\gd$ on~$B^4(2\gve) \setminus Z_\gd$
with residue~$-\mu$.
Consider the shell 
$$
V \,=\, \{ R \in (\gve,2\gve)\} \,=\, B^4(2\gve) \setminus \overline{B}{}^4(\gve) .
$$
Then $\gl_0$ and~$\gl_\gd$ are Liouville forms on~$V \setminus Z_\gd$, 
tame along~$Z_\gd$, with equal residue $-\mu$.
By Proposition~\ref{p:singliouville}~(ii) there exist a smooth bounded closed $1$-form~$\vartheta$ on $V \priv Z_\delta$ 
such that 
$$
\lambda_\delta = \lambda_0 + \vartheta. 
$$
Since $H_1(V \priv Z_\delta)$ is generated by 
two small loops around~$Z_\delta$ on which $\vartheta$ integrates to~$0$,
$\vartheta$ is the derivative of 
a smooth function $f \colon V \priv Z_\delta \to \R$. 
Since $\vartheta$ is bounded, $f$ extends to a Lipschitz function on~$V$. 
Let now $\rho \colon [\eps,2\eps] \to [0,1]$ be a smooth function that equals~$1$ near~$\gve$
and vanishes near~$2\eps$, and define the smooth Liouville form~$\gl'$ on $B^4(2\gve) \setminus Z_\gd$ by
$$
\left\{
\begin{array}{lcl}
\lambda' |_{\{R \leq \eps\}}        &:=& \lambda_\delta \\ [0.3em]
\lambda' |_{\{\eps \leq R \leq 2\eps\}} &:=& \gl_0 + d(\rho(R)f). 
\end{array}
\right.
$$
Since $\rho(R) f$ is a Lip\-schitz function, smooth on $V \setminus Z_\gd$,
$\gl'$ is tame on $V \setminus Z_\gd$ with residue $-\mu$, 
and the same holds on $\{ R \leq \gve \}$.
Moreover, $\gl' = \gl_0$ near the boundary of $B^4(2\gve)$.
This finishes the construction of the interpolating model $\bigl(Z_\gd \cap B^4(2\gve), \gl' \bigr)$.

Now take two components $\Sigma_i,\Sigma_j$ of $(\bfsigma, \gl)$ with equal weights $\mu_i=\mu_j$ 
as in the proposition.
Since $\Sigma$ has normal crossings and $\lambda$ is tame along~$\Sigma$ by assumption, 
there exists for each $p \in \Sigma_i \cap \Sigma_j$ a Darboux chart 
$\phi_p \colon B^4(3\eps_p) \to (M,\omega)$ centered at~$p$ such that
$$
\begin{array}{l}
\phi_p^{-1}(\Sigma)=\phi_p^{-1}(\Sigma_i\cup \Sigma_j)=\{z_1z_2=0\}, \\ [0.4em]
\phi_p^*\lambda=\lambda_{\st}-(\mu_id\theta_1+\mu_jd\theta_2)=\lambda_{\st}-\mu_i(d\theta_1+d\theta_2). 
\end{array}
$$
Now choose the $\gve_p$ so small that the Darboux balls $\phi_p(B^4(3\gve_p))$
are disjoint, and apply to each ball $B^4(2\gve_p)$ the above model
interpolation. When transplanted to $M$, this yields a Liouville form~$\gl'$ on $M \setminus \Sigma'$,
tame along~$\Sigma'$, with correct residues. 
Also notice that in the open case the Liouville flow of~$\gl'$ is backward complete, 
because $\gl'=\gl$ outside $\bigcup_p \phi_p (B^4(2\gve))$.
\proofend

\section{Symplectic embeddings} \label{s:sympemb}

In this section we prove the symplectic embedding results
Theorems~\ref{t:Umut}, \ref{t:Remb}, \ref{t:anything}, and~\ref{t:main}.
The proofs are all based on the same principle, stated in Theorem~\ref{t:embliouville}, 
which reduces the proofs to finding convenient polarisations. 
We first prove Theorem~\ref{t:main}. 
While Theorem~\ref{t:Umut} is a direct corollary,
Theorem~\ref{t:Remb} deals with an unbounded domain, whence the proof needs some adjustments.
We finally prove Theorem~\ref{t:anything}, using explicit smooth polarisations of the ball.

\subsection{From symplectic morphisms to exact symplectic morphisms.}
In this subsection we prove two lemmas that show that the exactness condition is almost for free.
The first lemma is useful if the polarisation of the source is smooth.

\begin{lemma} \label{l:funny} 
Let $(\Sigma, d\ga)$ and $(\Sigma', d\ga')$ be compact connected symplectic surfaces 
of genus~$g$ and~$g'$, respectively. 
Let $b$ be the number of boundary components of~$\Sigma$. Assume that 
$$
\area \Sigma < \area \Sigma'  \quad \mbox{ and } \quad g+b-1 \leq g'.
$$
Then there exists an $(\alpha,\alpha')$-exact 
symplectic embedding~$\varphi \colon \Sigma \to \Sigmapcirc$.  
\end{lemma}

\begin{figure}[h]   
 \begin{center}
  \psfrag{Si}{$\Sigma$}
	\psfrag{S}{$\Sigma'$}
	\psfrag{f}{$f$}
  \leavevmode\includegraphics{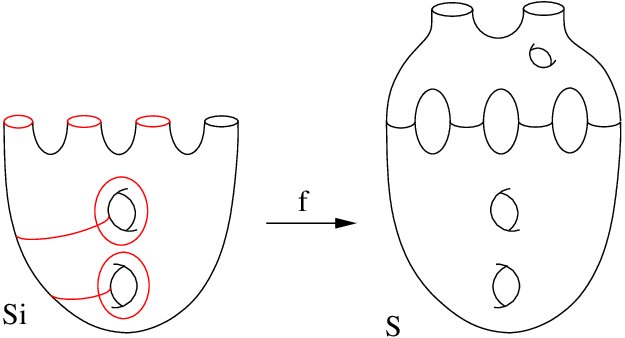}
 \end{center}
 \caption{The embedding $f \colon \Sigma \to \Sigmapcirc$ (in the figure it is the inclusion).}  \label{fig-surfacesf}
\end{figure}

\proof
A basis of $H_1(\Sigma,\ZZ)$ is given by oriented closed curves $\gamma_1, \dots, \gamma_{2g}$ 
around the $g$~holes
and $b-1$ oriented boundary components $\gamma_{2g+1}, \dots, \gamma_{2g+b-1}$,
see the red curves in Figure~\ref{fig-surfacesf}.
Since $g+b-1 \leq g'$, there exists an embedding $f \colon \Sigma \to \Sigmapcirc$
as shown in Figure~\ref{fig-surfacesf}.

Notice that $f_* \colon H_1(\Sigma) \to H_1(\Sigma', \pp \Sigma')$ is injective,
over $\ZZ$ and hence over~$\RR$.
In our context of oriented surfaces, the de Rham isomorphism gives 
$H_1(\Sigma;\RR)^* = H^1(\Sigma)$ and
$H_1(\Sigma',\partial \Sigma';\RR)^* = H^1_c(\Sigmapcirc)$, 
where both spaces on the right are the de Rham cohomology groups. 
Thus, $f^* \colon H^1_c(\Sigmapcirc) \to H^1(\Sigma)$ is surjective. 

Since the area of $\Sigma'$ is larger than the area of $\Sigma$, 
Moser's theorem allows to further assume that  
$f$ is area-preserving: $f^* d \ga' = d \ga$.
Write $\gamma_i' = f(\gamma_i)$.
Since the classes $[\gamma_i]$ are linearly independent in~$H_1(\Sigma)$ 
and since $f^*$ is surjective, 
there exists a compactly supported closed $1$-form~$\vartheta$ on~$\Sigma'$ 
with periods 
$$
\int_{\gamma_i'} \vartheta \,=\, \int_{\gamma_i} \ga - \int_{\gamma_i'} \ga', 
  \qquad  i=1, \dots, 2g+b-1. 
$$  
Let $\phi^t_\vartheta$ be the flow associated to the vector field~$X_\vartheta$ defined by 
$d\ga' (X_\vartheta,\cdot)=\vartheta$.
This flow is complete because $\vartheta$ is compactly supported,
and it is symplectic because $\vartheta$ is closed.
A classical computation for the flux shows that 
$$
\int_{\phi_\vartheta^t (\gamma_i')} \ga' \,=\, 
\int_{\gamma_i'} \ga' + t \int_{\gamma_i'} \vartheta , \qquad  i=1, \dots, 2g+b-1, 
$$
see e.g.\ \cite[$\S$\210.2]{McSa15}.
Thus $\phi_\vartheta^1 \circ f$ is an $(\alpha,\alpha')$-exact symplectic embedding of $\Sigma$ into~$\Sigma'$.
The proof of the proposition is complete.
\proofend

The next lemma deals with the exactness issue for singular polarisations.

\begin{lemma} \label{l:funny2}
Let $(\Sigma,d\alpha)$ and $(\Sigma',d\alpha')$ be compact symplectic multi-curves (with or without boundary) and let $\phi \colon \Sigma \to \Sigma'$ be a symplectic morphism 
whose restriction to each component of~$\Sigma$ is $(\alpha,\alpha')$-exact. 
Then there exists a Hamiltonian diffeomorphism~$f$ of~$\Sigma$ 
with support in an arbitrarily small neighbourhood of~$\sing\Sigma$, 
that is the identity near~$\sing \Sigma$, 
and is such that $\phi \circ f$ is an $(\alpha,\alpha')$-exact symplectic morphism. 
\end{lemma}

\proof
We can assume that $\Sigma$ has at least two components, so that each component~$\Sigma_i$
contains a singular point.
We fix in each $\Sigma_i$ a singular point~$x_i$.
If $\Sigma_i$ contains at least two singular points, we take for each singular point $x \neq x_i$
in $\Sigma_i$ two small open discs $D_x^i$ and $D_x^i{}'$ in~$\Sigma_i$ 
with $x \in D_x^i \subset \overline{D_x^i} \subset D_x^i{}'$, 
and a smooth path $\gamma_x^i$ from $x_i$ to $x$ in $\Sigma_i$. 
Define the real numbers
$$
\delta_i(x) := \int_{\gamma_x^i}\phi^*\alpha'-\alpha ,
$$  
which are responsible for the default of exactness of $\phi$. 
Let now $H_i \colon \Sigma_i \to \R$ be smooth functions with support in $\bigcup_x D_x^i{}'$ 
that equal $-\delta_i(x)$ on~$D_x^i$, and let $f \colon \Sigma \to \Sigma$ be the induced Hamiltonan diffeomorphism. It is well-defined, being the identity near $\sing \Sigma$ 
(and hence a symplectic morphism from $\Sigma$ to $\Sigma$). 
Moreover, 
$$
\int_{\gamma_x^i} (\phi\circ f)^*\alpha'=\int_{\gamma_x^i}\phi^*\alpha'-\delta_i(x)=\int_{\gamma_x^i}\alpha. 
$$ 
These equalities imply that for each $i$ the $(\phi\circ f)^*\alpha'$-action and the $\alpha$-action 
of every curve in~$\Sigma_i$ joining two singularities are equal. 
Indeed, for a curve $\gamma$ in~$\Sigma_i$ from $x$ to~$x_i$ (where $x=x_i$ is allowed), 
$\gamma \star \gamma_x^i$ provides a loop in~$\Sigma_i$, 
and since $\phi \circ f$ differs from $\phi$ by pre-composition with a Hamiltonian diffeomorphism 
and $\phi$ is $(\alpha,\alpha')$-exact on~$\Sigma_i$, so is $\phi\circ f$.
Hence 
\begin{eqnarray*}
0& = & \int_{\gamma\star \gamma_x^i} (\phi\circ f)^*\alpha'-\alpha \\
 & = & \int_\gamma (\phi\circ f)^*\alpha'-\alpha -\int_{\gamma_x^i} (\phi\circ f)^*\alpha'-\alpha \\
 & = & \int_\gamma (\phi\circ f)^*\alpha'-\alpha. 
\end{eqnarray*}
And if $\gamma$ joins two arbitrary singularities $x$ and $y$ in~$\Sigma_i$, then 
$\gamma \star \overline{\gamma_y^i} \star\gamma_x^i$ is a loop in~$\Sigma_i$, 
so the latter argument applies and shows the desired equality. 

We now claim that $\phi \circ f$ is an $(\alpha,\alpha')$-exact symplectic morphism between 
$\Sigma$ and~$\Sigma'$. Indeed, we have already observed that the special form of this diffeomorphism 
makes it $(\alpha,\alpha')$-exact on each component of~$\Sigma$. 
And if $\gamma$ is a loop in $\Sigma$ that is not contained in a single component of~$\Sigma$,
then (a time-translation of)~$\gamma$
is a concatenation $\gamma_1 \star \cdots \star \gamma_k$, where each $\gamma_j$
is a path in a component of $\Sigma$ connecting two singularities.
The integrals of $(\phi \circ f) \alpha'$ and~$\alpha$ agree on each $\gamma_j$
by the previous paragraph, and hence also on~$\gamma$.
\proofend

\begin{remark} \label{rk:funny2} 
{\rm 
The above proof shows that if $\phi$ is already known to preserve the action of all loops in~$\Sigma$ 
that pass only through singularities in a subset $\Gamma \subset \sing \Sigma$, 
then the Hamiltonian~$H$ can be taken with support in an arbitrarily small neighbourhood 
of $(\sing \Sigma) \priv \Gamma$.
}
\end{remark}

\subsection{The main embedding result (Theorem \ref{t:main})}
We first prove a special ``monotone" case of Theorem~\ref{t:main},
and then reduce the general case to this monotone case.

\subsubsection{The monotone case}

\begin{proposition} \label{p:monotone}
Let $\Gamma_1$ and $\Gamma_2$ be two regular grids in $D(A)$ and~$D(B)$, respectively, 
that cut $D(A)$ into $m$~topological discs of equal area~$a$ and 
$D(B)$ into $n$~topological discs of equal area~$b$. 
Then there exists an $(\ga_{\st},\ga_{\st})$-exact symplectic embedding
$$
\bigl( D(A) \times D(B) \bigr) \setminus (\Gamma_1 \times \Gamma_2)
\,\to\, Z^4(a+b) .
$$ 
\end{proposition}

\proof 
Let $\cp=\{p_1,\dots,p_m\}$, $\cq=\{q_1,\dots,q_n\}$,
$$
\bfsigma := 
\left\{ \bigl( \cp \times D (B),a \bigr), \bigl( D (A) \times \cq, b \bigr) \right\} ,
$$
and let $\lambda$ be the Liouville form on $\bigl( D(A) \times D(B) \bigr) \priv \Sigma$ 
provided by Lemma~\ref{l:poldim2} 
so that $(\bfsigma,\lambda)$
is an extendable Liouville polarisation of $D(A) \times D(B)$ 
with skeleton~$\Gamma_1 \times \Gamma_2$. 

\begin{figure}[h]   
 \begin{center}
    \psfrag{q1'}{$q_1'$}
		\psfrag{q2'}{$q_2'$}
	  \psfrag{q3'}{$q_3'$}
	  %\psfrag{q4'}{$q_4'$}
	  \psfrag{p1'}{$p_1'$}
		\psfrag{p2'}{$p_2'$}
	  \psfrag{S1}{\textcolor{blue}{$\Sigma_1'$}}
	  \psfrag{S2}{\textcolor{red}{$\Sigma_2'$}}
  \leavevmode\includegraphics{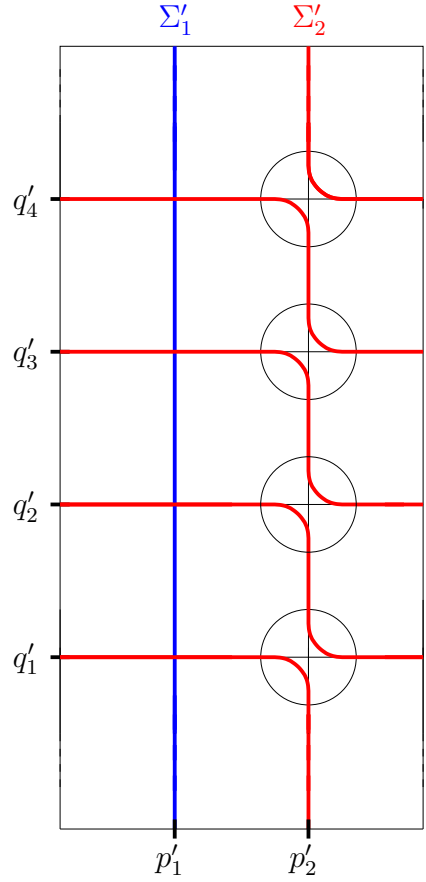}
 \end{center}
 \caption{$\Sigma_1'$ and $\Sigma_2'$ in $\Omega' = D (a+b) \times D (Kb)$ for $K=3$.}  \label{fig-S1S2}
\end{figure}

We now construct a Liouville polarisation of $\Om':=D(a+b) \times D(Kb)$ for $K$ large enough, 
to be fixed later in the proof. 
We refer to Figure~\ref{fig-S1S2} for an illustration of the construction.  
Take a regular grid (a smooth line) in $D(a+b)$ that divides $D(a+b)$ into topological discs 
of area $a$ and~$b$,
and take a point $p_1'$ in the disc of area~$a$ and~$p_2'$ in the disc of area~$b$. 
For a large integer~$K$ (to be chosen later), 
take a regular grid in~$D(Kb)$ that divides $D(Kb)$ into $K$ topological discs of area~$b$,
and choose a point~$q_j'$ in each disc.
Lemma~\ref{l:poldim2} yields a Liouville polarisation 
$$
\left(
\Bigl\{ 
\bigl( D (a+b) \times \{q_j'\}, b \bigr), \bigl( \{p_1'\} \times D (Kb),a \bigr), 
\bigl( \{p_2'\} \times D(Kb),b \bigr) 
\Bigr\} ,
\lambda' \right)
$$
of $\Om' = D(a+b) \times D(Kb)$, consisting of $K+2$ discs.

Define now $\Sigma_1' := \bigl\{ \{p_1'\} \times D (Kb) \bigr\}$ and let 
$\Sigma_2'$ be the symplectic curve obtained by resolving all the intersection points~$(p_2',q_j')$ 
of the $K+1$ remaining discs.
We can assume that the smoothening takes place in balls around these points that are disjoint
and do not intersect~$\partial \Om'$. 
Since $\lambda'$ has residue~$b$ along each component used to define~$\Sigma_2'$, 
Proposition~\ref{p:surgerform} shows that it can be modified in this union of balls 
to a Liouville form~$\lambda''$, tame along $\Sigma'_1 \cup \Sigma_2'$, 
with residue~$a$ along~$\Sigma'_1$ and $b$ along~$\Sigma'_2$. 
Since $\lambda''$ coincides with~$\lambda'$ near~$\partial \Om'$, its Liouville flow is also 
tangent to~$\partial \Om'$. In other terms, 
$$
\bigl( \bfsigma' := \{ (\Sigma_1',a),(\Sigma'_2,b) \},\lambda'' \bigr)
$$
is a Liouville polarisation of $\Om'$. 
Now $\Sigma_1'$ has area~$Kb$, $\Sigma_2'$ has area $K(a+2b)$, and 
$\# (\Sigma_1' \cap \Sigma_2') = K$. 
For $K$ large enough we therefore find 
$m$ disjoint topological discs~$D_1^i$ of area strictly larger than~$B$ in~$\Sigma_1'$ and 
$n$ disjoint topological discs~$D_2^j$ of area strictly larger than~$A$ in~$\Sigma_2'$ such that 
each $D_1^i$ intersects each $D_2^j$ exactly once. 
Hence there exists a symplectic morphism $\phi \colon \widehat \bfsigma \to \bfsigma'$.
Since the first cohomology of a disc vanishes, 
the restriction of~$\phi$ to each components of $\what \Sigma$
is $(\ga_{\st}, \ga_{\st})$-exact, so Lemma~\ref{l:funny2} shows that up to precomposing $\phi$ with 
a Hamiltonian diffeomorphism of~$\what \Sigma$, we can assume that $\phi$ is an
$(\alpha_\st,\alpha_\st)$-exact symplectic morphism.  
Theorem~\ref{t:embliouville} now guarantees an exact symplectic embedding 
of $\bigl( D(A) \times D(B) \bigr) \setminus \left( \Gamma_1 \times \Gamma_2 \right)$ 
into $D(a+b) \times D(Kb) \subset Z^4(a+b)$.
\proofend

\subsubsection{Reduction to the monotone setting}

\begin{proposition} \label{p:redmonotone}
Let $\Gamma_1 \subset D(A)$ and $\Gamma_2 \subset D(B)$ be regular grids whose complements are 
a union of $m$~topological discs of area~$a_i \leq a$ and 
of $n$~topological discs of area~$b_j \leq b$, respectively. 

Then there exists a regular grid $\Gamma'_1 \subset D(ma)$ 
that cuts $D(ma)$ into $m$~topological discs of area~$a$
and a regular grid $\Gamma_2' \subset D(nb)$ that cuts $D(nb)$ into 
$n$~topological discs of area~$b$, and an
$(\ga_{\st}, \ga_{\st})$-exact symplectic embedding  
$$
\bigl( D (A) \times D (B) \bigr) \setminus \bigl( \Gamma_1 \times \Gamma_2 \bigr)
\,\to\, 
\bigl( D (ma) \times D (nb) \bigr) \setminus \bigl( \Gamma_1' \times \Gamma_2' \bigr). 
$$
\end{proposition}

\proof
If $a_i=a$ and $b_j=b$ for all $i,j$, then there is nothing to prove. 
So after switching the factors and re-indexing
we can assume that $a_i<a$ for $1 \leq i \leq m' \leq m$ and $a_i=a$ for $i\geq m'+1$.
Define $n'\leq n$ for the $b_i$ in the same way, where however $n'=0$
if $b_i=b$ for all~$i$.

Write $\cp := \{p_1,\dots,p_m\} \subset D (A)$, $\cq := \{q_1,\dots,q_n\} \subset D (B)$, 
and let 
$$
\bigl( \bfsigma := \bigl\{ ( p_i \times D (B),a_i), ( D (A) \times q_j,b_j) \bigr\}, \lambda \bigr)
$$
be the extendable Liouville polarisation of $D (A) \times D (B)$ with skeleton $\Gamma_1 \times \Gamma_2$ provided by Lemma~\ref{l:poldim2}. 
%FF
We recall that in this case, an extension is obtained by taking somewhat larger discs, so
$\overline \Sigma$ is just the closure of $\Sigma$ in~$\C^2$, i.e.\ the closed discs.
Similarly, let
$\cp' := \{p_1',\dots,p_m'\} \subset D(ma)$, $\cq':=\{q_1',\dots,q_n'\} \subset D (nb)$, 
and let 
$$
\bigl( \bfsigma' := \bigl\{ (\cp' \times D(nb),a), ( D(ma) \times \cq', b) \bigr\}, \lambda' \bigr)
$$
be the Liouville polarisation of $D (ma) \times D (nb)$ 
with skeleton $\Gamma'_{\!1} \times \Gamma'_{\!2}$. 
Since we are free to choose $\Gamma'_{\!1}$ and $\Gamma'_{\!2}$, 
and since $D(A) \subset D(ma)$ and $D(B) \subset D(nb)$, 
we may as well assume that $p_i'=p_i$ and $q_j'=q_j$. 
This is not an important assumption but will simplify the proof at some point. 
We will nevertheless keep the primes in order to distinguish the source and the target.  
Lemma~\ref{l:double} provides disjoint open discs~$D_i^A$ in~$D(ma) \setminus \Gamma_1$, 
disjoint open discs~$D_j^B$ in~$D(nb) \setminus \Gamma_2$, 
points $p_i'^1, p_i'^2 \in D_i^A$ for $i\leq m'$, $q_j'^1, q_j'^2 \in D_j^B$ for $j\leq n'$, 
and a 
Liouville polarisation $(\bfsigma'',\lambda'')$ of $D(ma) \times D(nb)$ with
\begin{eqnarray*}
\bfsigma'' &\,:=\,  \Bigl\{  & 
 \bigl( p_i'^1 \times D (nb),a_i \bigr)_{i\leq m'}, \bigl( p_i'^2 \times D (nb), a-a_i \bigr)_{i\leq m'}, \\
&& 
  \bigl( D (ma) \times q_j'^1,b_j \bigr)_{j\leq n'}, \bigl( D (ma) \times q_j'^2, b-b_j \bigr)_{j\leq n'},\\
 & &   \bigl( p_i' \times D (nb),a \bigr)_{i\geq m'+1}, \bigl( D (ma) \times q_j',b \bigr)_{j\geq n'+1}\Bigr\}
\end{eqnarray*}
and such that 
\begin{eqnarray} \label{e:BBB}
&&\cb \Bigl( \bigcup_i D_i^A \times D (nb) \cup \bigcup_j D (ma) \times D_j^B,\lambda'' \Bigr) \\
&\subset&
 \bigl( D (ma) \times D (nb) \bigr) \setminus (\Gamma_1' \times \Gamma_2' ) .
\notag
\end{eqnarray}
As is clear from the proof of Lemma \ref{l:double} we can chose $p_i'^1=p_i'$ and $q_i'^1=q_i'$. 
We can find $\eps>0$ and two area preserving embeddings 
$$
\begin{array}{l}
\sigma \colon \overline{D(B)} \to \overline{D(nb)} \\ [0.2em]
\tau \colon D(A+\eps) \to D(ma)
\end{array}
$$  
such that $\sigma(q_j)=q_j'=q_j$ for all $j \leq n$ and 
the image $\im \sigma$ avoids all the~$q_j'^2$, 
and such that $\tau(p_i)=p_i'$ for all $i \leq m$ and 
$\im \tau$ avoids all the~$p_i'^2$. 
This is obvious for $\tau$ because $A<ma$ by assumption.
It is also true for~$\sigma$ because if $B<nb$ we can even extend $\sigma$ 
to $\hat \sigma \colon D(B+\eps)\to D(nb)$ with the required properties, 
and if $B=nb$ then $b_j=b$ for all $j$, so there is no~$q_j'^2$ to avoid, 
and $\sigma$ can be taken to be the identity. 
The two maps $\sigma$ and~$\tau$ induce an embedding 
$\phi \colon \Sigma\to \Sigma''$ defined by 
$$
\begin{array}{l}
\phi|_{p_i\times D(B)}=p_i'\times \sigma \\ [0.2em]
\phi|_{D(A)\times q_j}=\tau\times q_j'. 
\end{array}
$$  
There are now two cases.

\s \ni
\underline{Case $B<nb$:}
As already noticed, in this case $\sigma$ can be extended to 
an area preserving embedding
$\hat \sigma \colon D(B+\eps)\to D(nb)$, 
so $\phi$ can be extended to a smooth area preserving embedding 
$\what \Sigma \to \Sigma''$ 
which is a symplectic morphism between our polarisations. 
Its restriction to each component of $\Sigma$ is 
$(\alpha_\st,\alpha_\st)$-exact, since all components are discs. 
By Lemma~\ref{l:funny2}, $\phi$ can be modified to an $(\alpha_\st,\alpha_\st)$-exact symplectic morphism, 
and Theorem~\ref{t:embliouville} now implies that 
\begin{equation} \label{e:DD12}
\bigl( D (A) \times D (B) \bigr) \setminus (\Gamma_1 \times \Gamma_2) = \cb( \bfsigma, \lambda) \,\hraa\, \cb(\bfsigma'',\lambda'') .
\end{equation}
Moreover, since 
$$
\bfsigma'' \,\subset\,
\bigcup_i D_i^A \times D(nb) \cup \bigcup_j D(ma) \times D_j^B ,
$$
we obtain together with \eqref{e:BBB} the inclusion
\begin{eqnarray} \label{e:BDDGG}
\cb (\bfsigma'',\lambda'') &\subset& 
            \cb \Bigl( \bigcup_i D_i^A \times D(nb) \cup \bigcup_j D(ma) \times D_j^B \Bigr) \\
&\subset&  \bigl( D (ma) \times D (nb) \bigr) \setminus (\Gamma_1' \times \Gamma_2') .
\notag
\end{eqnarray}
Composing \eqref{e:DD12} with \eqref{e:BDDGG} we obtain the asserted exact symplectic embedding. 

\s \ni 
\underline{Case $B=nb$:} 
In this case, $\sigma=\id$ does not extend to a larger disc 
with image still in~$D(nb)$, so formally the above argument does not apply.  
We solve this difficulty by revisiting the proof of Theorem~\ref{t:embliouville} 
in this particular case: 
Although $\phi$ does not extend to an extended polarisation~$\what \Sigma$ 
as before, the product map $\psi := \tau \times \id$ extends $\phi$ 
to a small \nbd of~$\Sigma$ in~$D(A) \times D(B)$, 
still with image in~$D(ma) \times D(nb)$. 
The forms $\lambda$ on $D(A) \times D(B)$ and $\lambda''$ on~$D(ma) \times D(nb)$ are both split:  
$$
\lambda=\lambda_1 \oplus \lambda_2, \qquad \lambda''=\lambda_1''\oplus \lambda_2'',
$$
where $\lambda_1, \lambda_2$ are defined on 
$D(A) \priv \{p_i\}$, $D(B) \priv \{q_j\}$
and $\lambda_1'', \lambda_2''$ 
are defined on $D(ma) \priv \{p_i',p_i'^2\}$, $D(mb) \priv \{q_j'\}$. 
Since $B=mb$ and $q_j=q_j'$ we can take $\lambda_2=\lambda_2''$, 
and since $\lambda_1$ and $\lambda_1''$ have the same residues at $p_i$ and $p_i'=p_i$, 
we can make them coincide on a small \nbd of the~$p_i$. 
Since $\tau$ can clearly be taken to be the identity on a \nbd 
of the points~$p_i$, 
the map $\psi =\tau\times \id$ (defined on a \nbd of $\Sigma$) 
then pulls back $\lambda''$ to~$\lambda$. 
Thus the basic conjugacy procedure described in the proof of 
Theorem~\ref{t:embliouville} to embed 
$\cb(\bfsigma,\lambda)$ into $\cb(\bfsigma'',\lambda'')$ applies, 
and now the argument in the previous case applies. 
\proofend

\subsubsection{Proof of Theorem \ref{t:main}}
Let $\Gamma_1 \subset D(A)$ and $\Gamma_2 \subset D(B)$ be two regular grids that cut the discs into 
topological discs of area $\leq a$ and~$\leq b$, respectively. 
By Proposition~\ref{p:redmonotone}, 
$$
\bigr( D (A) \times D (B) \bigr) \setminus (\Gamma_1 \times \Gamma_2)
\,\hraa\, \bigl( D (ma) \times D (nb) \bigr) \setminus (\Gamma_a \times \Gamma_b)
$$ 
where $\Gamma_a$ and $\Gamma_b$ cut $D(ma)$ and~$D(nb)$ into topological discs of area~$a$ and~$b$, respectively. 
By Proposition~\ref{p:monotone}, we also have
$$
\bigl( D (ma) \times D (nb) \bigr) \setminus (\Gamma_a \times \Gamma_b) 
\,\hraa\, Z^4(a+b) .
$$
Composing these two maps we obtain the searched exact symplectic embedding  
$$
f \colon \bigl( D (A) \times D (B) \bigr) \setminus (\Gamma_1 \times \Gamma_2)
\hraa Z^4(a+b) .
$$ 
\proofend

%%%%%%%%%%%%%%%%%%%%%%%%%%%%%%%%%%%%%%%%%%%%%%%%%%%%%%%%%%%%%%%%%%%%%%%%%%%%%%%%%

\subsection{Proof of Theorem \ref{t:Remb}}
Let $\Gamma := \Gamma_1 \times \Gamma_1\subset \R^4$, where $\Gamma_1 \subset \R^2$ is the regular grid 
$$
\Gamma_1:=\bigcup_{(n,m) \in \Z^2} \{n\}\times \R\cup \R\times \{m\}.  
$$
Let also 
$$
\Sigma := \bigcup_{(n,m) \in \Z^2} \bigl\{ (n+\tfrac 12,m+\tfrac 12) \bigr\} \times \C \cup 
                                    \C \times \bigl\{ (n+\tfrac 12,m+\tfrac 12) \bigr\}.
$$

\begin{lemma} \label{l:lambdasquare} 
There exists a Liouville form $\lambda$ on $\R^4\priv \Sigma$ with the following properties:
\begin{itemize}
\item $\lambda$ has residue $1$ along $\Sigma$.

\s
\item 
For each $p \in \C^2 \priv \Sigma$, the Liouville trajectory $\Phi^t_{X_\lambda}(p)$ is defined for 
$t \in \; ]-\infty,t^+(p)[$, where $t^+(p) \in \; ]0,+\infty]$. 
If $t^+(p) < +\infty$, then $\lim_{t \to t^+(p)} \Phi^t_{X_\lambda}(p)$ exists and belong to~$\Sigma$.

\s
\item $\Gamma$ is the skeleton of $(\Sigma, \lambda)$: 
$$
\Gamma=\{p\;|\; t^+(p) = +\infty\} = \R^4 \priv \cb(\Sigma,\lambda). 
$$

\item $X_\lambda$ is tangent to the hyperplanes 
$\{x_1\in \Z\}$, $\{x_2\in \Z\}$, $\{y_1\in \Z\}$, $\{y_2\in \Z\}$. 
\end{itemize}
\end{lemma}

\proof
Again, the form $\lambda$ is a product: 
$$
\lambda_{(z,w)}(u,v) = \lambda_z(u)+\lambda_w(v),
$$
where $\lambda_z=\lambda_w$ is a Liouville form on $\R^2 \priv \bigl( \Z^2 + (\frac 12,\frac 12) \bigr)$, 
with $\Gamma_1$ as skeleton. 
We construct this Liouville form exactly as in Section~\ref{ss:dim2}, where compactness plays no role.
Alternatively, we can do the construction in the proof of Proposition~\ref{p:poldim1} 
on $[-2,2]^2$ in a $\ZZ^2$-periodic way.
We then obtain a smooth Liouville form on $[0,1]^2$ that extends to a smooth $\ZZ^2$-periodic
Liouville form on $\R^2 \priv \bigl( \Z^2 + (\frac 12,\frac 12) \bigr)$ with skeleton~$\Gamma_1$. 
\proofend

\begin{figure}[h!] 
\begin{center}
	\psfrag{12}{\tiny{$\frac 12$}}
  \psfrag{1}{\tiny{$1$}}
 \leavevmode\includegraphics{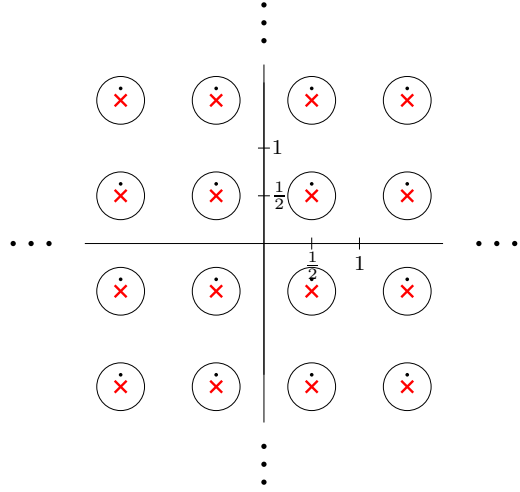}
\end{center}
\caption{The parametrisation of the curve $\Sigma_2'$. Each depicted disc $D(\gve)$ 
has area $\gve$ on the picture, but carries the symplectic form~$\omega_\gve$ of area~$2+\eps$. 
The black points are the punctures, and the red crosses are the intersection points with~$\Sigma_1'$.} \label{fig:Rembpict1}
\end{figure}

From Lemma~\ref{l:lambdasquare} we obtain the following two facts:
\begin{itemize}
\item 
For $N \in \NN$ define 
$$
\Sigma_N := \Sigma \2 \cap \, ]-N,N[^4 \quad \mbox{and} \quad \Gamma_N := \Gamma \times \2 ]-N,N[^4 .
$$ 
Then
$\bigl( (\Sigma_N,1), \lambda \bigr)$ is an extendable Liouville polarisation of~$]-N,N[^4$ 
whose skeleton is~$\Gamma_N$. 
Notice that $\Sigma_N \subset \Sigma_{N+1}$.  

\s
\item  
Define 
$$
Z(2) := \: ]-1,1[ \2 \times \2 ]0,1[ \2 \times \R^2 \quad \mbox{and}  \quad \Sigma_Z := \Sigma \cap Z(2) .
$$ 
Then $\bigl( (\Sigma_Z,1), \lambda \bigr)$ is a Liouville polarisation of~$Z(2)$. 
\end{itemize}
Now $\Sigma_Z$ is made of two planes $\{ (\pm \frac 12,\frac 12) \} \times \R^2$ 
and infinitely many discs 
$$
D_{n,m} := \: ]-1,1[ \2 \times \2 ]0,1[ \2 \times \bigl\{ (n+ \tfrac 12,m + \tfrac 12) \bigr\}, 
\quad n,m \in \Z. 
$$
As in the proof of Proposition~\ref{p:monotone}, 
we keep $\Sigma_1' := \bigl\{ ( -\frac 12,\frac 12) \bigr\} \times \R^2$ 
and glue all other components together to form a symplectic curve~$\Sigma_2'$, 
which is diffeomorphic to a plane with infinitely many punctures arranged on a lattice. 
Each puncture corresponds to a boundary component of~action~$2$ (see Figure~\ref{fig:Rembpict1}). 
Since $\lambda$ has residue~$1$ along each component of~$\Sigma_Z$, 
Proposition~\ref{p:surgerform} (applied with $U$ the union of small balls around the singular points
of~$\Sigma_2'$) yields
a Liouville form~$\lambda'$ on $Z(2) \priv (\Sigma_1'\cup \Sigma'_2)$ 
tame along $\Sigma' := \Sigma_1'\cup \Sigma_2'$ with residue~$1$ and
equal to $\lambda$ near~$\partial Z(2)$. 
In other terms, 
$$
\bigl( \{ (\Sigma'_1,1), (\Sigma'_2,1) \}, \lambda'\bigr)
$$ 
is a Liouville polarisation of $Z(2)$. 

\begin{lemma}\label{l:compatible}
For each $N$ there exists an $(\ga_{\st},\ga_{\st})$-exact 
symplectic morphism $\phi_N \colon \widehat \Sigma_N \to \Sigma'$ such that 
$$
\phi_{N+1} |_{\widehat \Sigma_N} = \phi_N . 
$$
\end{lemma}

\proof
For describing the construction of $\phi_N$ it is convenient to redescribe the 
components $\Sigma'_1$ and $\Sigma_2'$ of $\Sigma'$: 
\begin{itemize}
\item 
$\Sigma_1' = \{(-\frac 12,\frac 12)\} \times \R^2$ is symplectic with standard area form on the $\R^2$-factor. Its intersections with $\Sigma_2'$ are the points 
$$
I_{n,m} := (-\tfrac 12, \tfrac 12, n+\tfrac 12, m+\tfrac 12), \quad (n,m) \in \Z^2 .
$$ 
\item 
$\Sigma_2'$ is obtained by gluing $(\frac 12, \frac 12) \times \R^2$ with the topological discs 
$D_{n,m}$, $(n,m) \in \Z^2$. 
This can be done by replacing each round disc 
$(\frac 12, \frac 12) \times D \bigl( (n+\frac 12,m+\frac 12),\eps \bigr)$ 
by its gluing with the disc~$D_{n,m}$. 
These gluings can be parameterised by $\Z^2$-translations of 
a punctured disc $D(0,\eps) \priv \{(0, \frac \gve2 \}$ with an area form~$\om_\eps$ 
whose total area is~$2+\eps$, and that is standard near the boundary~$S^1(\eps)$. 
Each point~$I_{n,m}$ is again parameterised by $(n+\frac 12, m+\frac 12)$. 
\end{itemize}

We also define $\Sigma_1'(N) := \{ (-\frac 12,\frac 12) \} \times [-N,N]^2$ and 
$\Sigma_2'(N) \subset \Sigma_2'$ the part of $\Sigma_2'$ parameterised by 
$[-N,N]^2 \priv (\Z^2 + (\frac 12, \frac 12 + \frac \gve 2 ))$, 
and notice that they correspond exactly to the intersections of these curves with $Z(2)\cap [-N,N]^4$ 
(see Figure~\ref{fig:Rembpict1}).

Now divide $\Sigma_N$ into its $4N^2$ horizontal components 
$$
\Sigma_N^H(n,m) := \; ]-N,N[^2 \times \bigl\{ (n+\tfrac 12,m+\tfrac 12) \bigr\}, \quad (n,m) \in [-N,N-1]^2
$$
and its $4N^2$ vertical components 
$$
\Sigma_N^V(n,m) := \bigl\{ (n+\tfrac 12,m+\tfrac 12) \bigr\} \times \, ]-N,N[^2, \quad (n,m)\in [-N,N-1]^2. 
$$
Each such component is a disc of area $4N^2$, the horizontal and vertical components do not intersect pairwise, but each horizontal component intersects each vertical component exactly once, 
thus each component passes through $4N^2$~singular points of~$\Sigma_N$.

\begin{figure}[h!]
\begin{center}
  \psfrag{3}{\tiny{$3$}}
  \psfrag{9}{\tiny{$9$}}
	\psfrag{4}{\textcolor{green}{$\Sigma_1^V(0,0)$: area 4,}}
	\psfrag{4c}{\textcolor{green}{4 copies}}
  \psfrag{16}{\textcolor{blue}{$\Sigma_2^V(0,0)$: area 16}}
	\psfrag{16c}{\textcolor{blue}{16 copies}}
	\psfrag{18}{area $18^2 >16^2$}
 \leavevmode\includegraphics{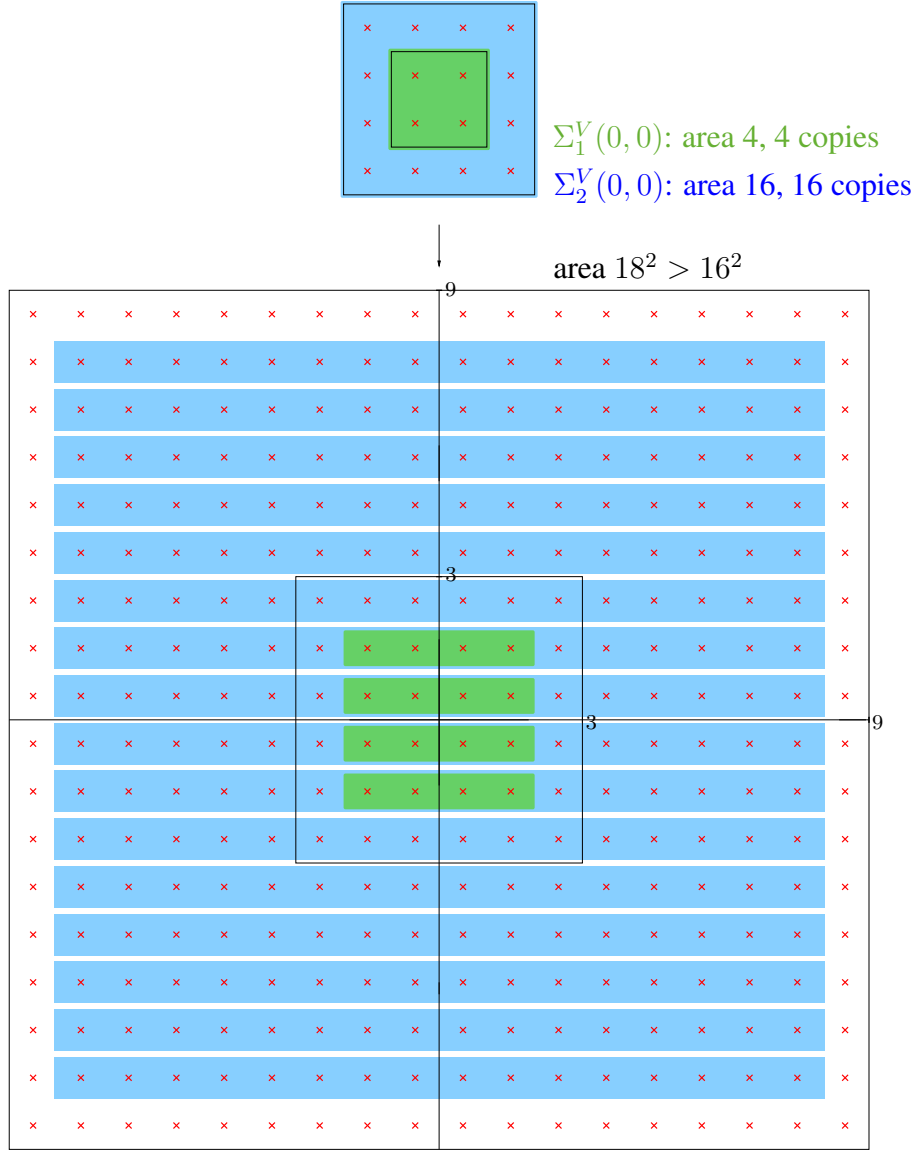}
\end{center}
\caption{The restriction of $\gf_1$ to $\Sigma_1^V$ and of $\gf_2$ to~$\Sigma_2^V$. 
Each band has height slightly less than~$1$, so that we can construct a larger
band of height less than~$1$ in the next step~$N+1$.}  \label{fig:Rembpict2}
\end{figure}

\begin{figure}[h!] 
\begin{center}
  \psfrag{3}{\tiny{$3$}}
  \psfrag{9}{\tiny{$9$}}
	\psfrag{4}{\textcolor{green}{$\Sigma_1^V(0,0)$: area 4,}}
	\psfrag{4c}{\textcolor{green}{4 copies}}
  \psfrag{16}{\textcolor{blue}{$\Sigma_2^V(0,0)$: area 16}}
	\psfrag{16c}{\textcolor{blue}{16 copies}}
	\psfrag{18}{area $18^2 >16^2$}
 \leavevmode\includegraphics{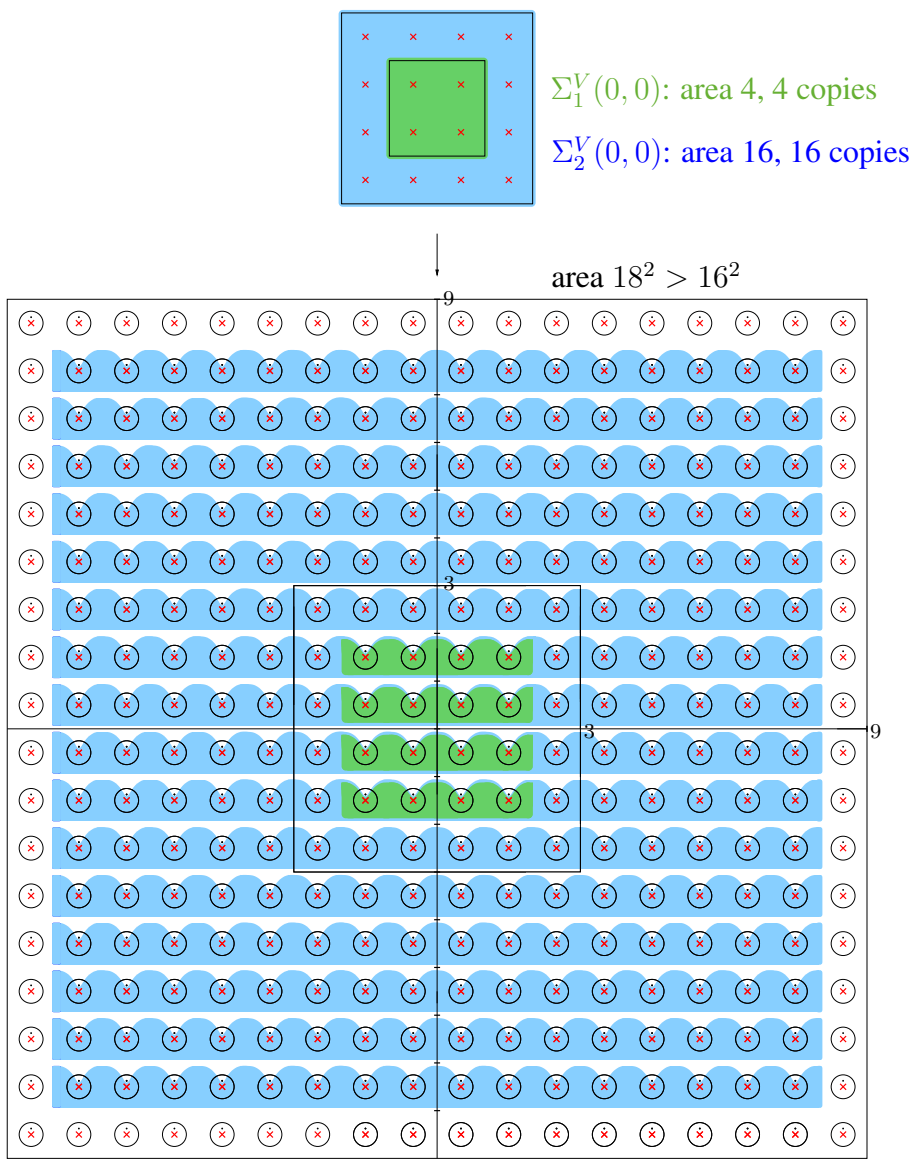}
\end{center}
\caption{The restriction of $\gf_1$ to $\Sigma_1^H$ and of $\gf_2$ to~$\Sigma_2^H$. 
These are almost the same embeddings as in the previous figure,  
except that they avoid the punctures. Note that the area $2+\gve$ of the discs 
is not even used since for $\gve >0$ small enough the area of~$\Sigma_2'(2N^2+1)$ without the discs
is already larger than the area of $\Sigma_N^H$.} \label{fig:Rembpict3}
\end{figure}

We first construct a symplectic morphism $\gf_N$ by sending a slightly larger disc around the closure of each $\Sigma_N^V(n,m)$ 
into $\Sigma_1'(2N^2+1)$, 
and by sending a slightly larger disc around the closure of each~$\Sigma_N^H(n,m)$ into 
$\Sigma'_2(2N^2+1)$. 
The areas of these surfaces are 
$$
\begin{array}{l}
\ca(\Sigma_1'(2N^2+1)) = (2(2N^2+1))^2 > 16N^4 = \ca(\Sigma^V_N) \\ [0.3em]
\ca(\Sigma_2'(2N^2+1)) = (2(2N^2+1))^2 + 2(2N^2+1)^2>16 N^4 = \ca(\Sigma_N^H),
\end{array}
$$
where $\Sigma_N^H := \coprod \Sigma_N^H(n,m)$ and $\Sigma_N^V := \coprod \Sigma_N^V(n,m)$, 
and
$$
\# \bigl( \Sigma_1'(2N^2+1) \cap \Sigma_2'(2N^2+1) \bigr) = (2(2N^2+1))^2 > 16N^4 = 
           \# (\Sigma_N^V \cap \Sigma_N^H). 
$$
We can therefore define $\gf_N$ by its restrictions to $\Sigma_1^V$ and~$\Sigma_1^H$ 
as illustrated in Figures~\ref{fig:Rembpict2} and~\ref{fig:Rembpict3},
where $\gf_N$ is explained in blue and $\gf_{N+1}$ in yellow for~$N=1$.
These are area preserving embeddings that take 
the points at the red crosses to points at the red crosses.
This defines symplectic morphisms  $\gf_N \colon \widehat \Sigma_N \to \Sigma'$,
and we can clearly choose $\gf_{N+1}$ such that $\gf_{N+1} |_{\widehat \Sigma_N} = \gf_N$.  

The restrictions of these symplectic morphisms to the components of $\what \Sigma_N$, which are discs, 
are ($\ga_{\st},\ga_{\st}$)-exact.
Dealing with the exactness statement is now a matter of applying Lemma~\ref{l:funny2} in a coherent way: 
This lemma guarantees that $\gf_N$ can be modified to an exact symplectic morphism~$\phi_N = \gf_N \circ f_N$ 
by precomposition with a Hamiltonian diffeomorphism~$f_N$ of~$\what \Sigma_N$ with compact support. 
Then $\gf_{N+1} \circ f_N|_{\what \Sigma_N} = \gf_N\circ f_N|_{\what \Sigma_N}$, 
so Remark~\ref{rk:funny2} shows that $f_{N+1}$ can be taken to coincide with~$f_N$ on~$\what \Sigma_N$. 
One now readily checks that the $\phi_N$ form a sequence of ($\ga_{\st},\ga_{\st}$)-exact symplectic morphisms such that 
$\phi_{N+1} |_{\widehat \Sigma_N}=\phi_N$.
\proofend

\noindent{\it Proof of Theorem \ref{t:Remb}.} 
We can now proceed as in the proof of Theorem~\ref{t:embliouville}.
First we extend each~$\phi_N$ from Lemma~\ref{l:compatible}
to a symplectic embedding 
$$
\psi_N \colon V_N \to \Op \bigl( \Sigma', Z(2) \bigr)
$$
of an open neighbourhood~$V_N$ of~$\widehat \Sigma_N$.
Set $U_N := V_N \cap [-N,N]^4$.
Since $\phi_{N+1} |_{\widehat \Sigma_N} = \phi_N$,
Moser's method shows that the maps $\psi_N$ can be chosen such that
$$
U_{N+1}\supset U_N \quad \text{ and } \quad \psi_{N+1}|_{U_N} = \psi_N |_{U_N} . 
$$
We can therefore define the $C^\infty$-smooth symplectic embedding
$$
\psi \colon \bigl( U, \Sigma \bigr) \to \bigl( Z(2),\Sigma' \bigr) 
$$
on $U := \bigcup U_N$
by $\psi |_{U_N} := \psi_N |_{U_N}$.

In the construction of the maps $\phi_N$ of Lemma~\ref{l:compatible}
we can assume that their images avoid a whole neighbourhood of the punctures.
After choosing $U$ smaller, if necessary, 
the closure of $\psi (U)$ is then disjoint from the boundary of~$Z(2)$.
Using Lemma~\ref{l:transferliouville2} as in the proof of Theorem~\ref{t:embliouville}
we correct~$\psi$ to a $C^1$-smooth symplectic embedding
$$
\Phi \colon \bigl( \widehat U, \Sigma \bigr) \to \bigl( Z(2), \Sigma' \bigr)
$$
defined on a neighourhood $\widehat U$ of~$\Sigma$ such that $\Phi (\widehat U) \subset \psi (U)$,
whence the closure of $\Phi (\widehat U)$ is also disjoint from $\partial Z(2)$.
Choosing $\widehat U$ smaller if necessary, we can also assume that $\widehat U$ retracts onto~$\Sigma$,
whose components are topological discs.
We therefore find a Liouville form~$\lambda''$ on~$Z(2)$ that coincides with~$\Phi_*\lambda$ 
on $\Phi (\widehat U)$ and with $\lambda'$ near~$\partial Z(2)$. 
Then the map 
\fonction{\Psi_1}{\cb(\Sigma,\lambda)}{Z(2)}{p}{\phi_{X_{\lambda''}}^{-t^+(p)+\eps(p)} \circ \Phi \circ \phi_{X_\lambda}^{t^+(p)-\eps(p)},}
defined for any small enough positive function~$\eps(p)$, 
is a $C^1$-smooth symplectic embedding. 

Since the $\phi_N$ are $(\ga_{\st} |_{\Sigma_N}, \ga_{\st} |_{\Sigma'})$-exact, 
$\Psi_1$ is $(\ga_{\st}, \ga_{\st})$-exact.
The asserted $(\ga_{\st}, \ga_{\st})$-exact 
symplectic embedding 
$\Psi \colon \cb(\Sigma,\lambda) = \RR^4 \setminus \Gamma \to Z(2)$ is now obtained from 
Lemma~\ref{le:Zehnder}.
\proofend

\subsection{Proof of Theorem \ref{t:anything}}
\label{ss:baby}

In this subsection we first give another proof of Theorem \ref{t:main} for the ball as domain, 
that does not use a singular polarisation but Biran's smooth polarisation.
We then use this result to prove Theorem~\ref{t:anything}.

\begin{theorem}  \label{t:baby1}
$B^4(1) \priv \Delta_k  \; \hraa \; Z^4\bigl( \tfrac 2k \bigr)$.
\end{theorem}

\proof 
The polarisation of $\CP^2$ of degree~$k$ in Example~\ref{ex:polp2}
restricts to a tame Liouville polarisation $(\Sigma_k,\lambda_k)$ 
of the ball 
whose skeleton is~$\Delta_k \cap B^4(1)$, so
$$
B^4(1) \priv \Delta_k   = \cb \bigl( \Sigma_k,\lambda_k \bigr). 
$$
One readily checks that this Liouville polarisation is extendable.
%PP \red{C'est un pour court, je sais. Pour ceci on aura l'appendix (tame et extendable)} 
%Non je ne crois pas qu''il y ait besoin d'appendice pour ça.
Moreover, the area, genus, and number of punctures of $\Sigma_k$ are well-known and given by 
$$
\ca_\om (\Sigma_k ) =k , \quad g(\Sigma_k) = \frac{(k-1)(k-2)}2, \quad b (\Sigma_k) =k. 
$$
On the other side, proceeding as in the proof of Proposition~\ref{p:monotone}, 
we take for any $A \in \NN$ two points~$p_i$ in $D(\frac 2k)$
and $kA$ points~$q_j$ in~$D(A)$, all with weight~$\frac 1k$,
and consider the union of the discs $D(\frac 2k) \times q_j$ and $p_i \times D(A)$.
Using Lemma~\ref{l:poldim2} and resolving all the $2kA$ intersections as in Proposition~\ref{p:surgerform}, 
we obtain a smooth Liouville polarisation $\Sigma_k'$ of $D(\frac 2k) \times D(A)$
with residue~$\frac 1k$. Then
$$
\ca_\om(\Sigma_k') = 4A \quad \mbox{ and } \quad g(\Sigma_k')=kA-1.
$$
When $A \geq \frac k2$, we have $\ca_\om(\Sigma'_k) > \ca_\om(\Sigma_k)$
and  $g(\Sigma_k') \geq g(\Sigma_k)+b(\Sigma_k)-1$ since $k \geq 2$,
so by Lemma~\ref{l:funny} there exists an 
$(\alpha_{\st}|_{\Sigma_k},\alpha_{\st}|_{\Sigma_k'})$-exact 
symplectic embedding 
of $\widehat \Sigma_k$ into~$\Sigma_k'$, and by 
Theorem~\ref{t:embliouville} there exists an $(\alpha_{\st},\alpha_{\st})$-exact  
symplectic embedding of
$\cb \bigl( \Sigma_k,\lambda_k \bigr) = B^4(1) \priv \Delta_k$
into~$Z^4(\frac 2k)$. 
\proofend

\begin{remark} \label{rem:volumefill} 
{\rm
When $k=2$ we can get a better estimate and find a (non-exact) symplectic embedding 
$$
\Phi \colon B^4(1) \priv \Delta_2 = B^4(1) \priv \R^2 \,\to\, E(\tfrac 12,2+\eps)
$$ 
into the ellipsoid $E(\frac 12,2+\eps)$, like~\cite{SSVZ24}. 
Indeed, $\Sigma_2$ is a sphere with two punctures.
Thus a cylindrical extension of~$\Sigma_2$ symplectically embeds into the disc $D(2+\eps)$ . 
Since by Example~\ref{ex:ext}, $(D(2+\eps),\frac 12)$ is a polarisation of~$\sdb(D(2+\eps),\frac  12)$, 
Theorem~\ref{t:embliouville} (in its version without exactness assumption and conclusion) 
shows that $B^4(1)\priv \Delta_2$ embeds into~$\sdb(D(2+\eps),\frac 12)$. 
And by \cite[Lemma 2.1]{Op07}, 
$\sdb(D(2+\eps),\frac 12)$ is symplectomorphic to $E(2+\eps,\frac 12)$. 
Note that the embedding $\Sigma_2 \to D(2+\eps)$ is not exact since the boundary circles of~$\Sigma_2$ 
in its extension both have area~$1$, and so $\Phi$ is not exact. 
See \cite[$\S$\26.4]{SSVZ24} for a related discussion.
\diam
}
\end{remark}

Recall that the open ellipsoid $E(a,b)$ is defined as 
$$
E(a,b) \,=\, \left\{ (z,w) \in \CC^2 \;\bigg|\; \frac{|z|^2}{a} + \frac{|w|^2}{b} < 1 \right\} .
$$

\begin{proposition}\label{p:embell+}
For $d,N \in \NN$ coprime there exists
for every $m \in \NN$ with $m \geq d$ and $(m, d, N) \neq (2,2,1)$ 
a symplectic embedding 
$$
B^4(1)\setminus \Delta_{m N d} \,\hraa\, E \bigl( \tfrac 1d, d+\tfrac 1N \bigr).
$$
\end{proposition}

\proof 
As above, consider the extendable Liouville polarisation $(\Sigma_{m Nd}, \lambda_{m Nd})$ 
of~$B^4(1)$ of degree~$m Nd$, whose skeleton is $\Delta_{m Nd}$.
We have that $\ca_\om(\Sigma_{m Nd}) = m Nd$, 
$$
g(\Sigma_{m Nd})=\frac 12  (m Nd-1)(m Nd-2) \quad \mbox{and} \quad 
b(\Sigma_{m Nd})=m Nd. 
$$
In order to produce the required Liouville polarisation of the target ellipsoid, consider the symplectic ramified covering 
\fonction{\Phi}{B^4(N(Nd^2+d))}{E(N,Nd^2+d)}{(R_1,\theta_1,R_2,\theta_2)}{\big(\frac {R_1}{Nd^2+d},(Nd^2+d)\theta_1,\frac{R_2}{N},N\theta_2\big).} 
One checks without difficulty that the Liouville polarisation $(\Sigma_m,\lambda_m)$  
of $B^4(N(Nd^2+d))$ of degree~$m$
projects under~$\Phi$ to a smooth Liouville polarisation of degree~$m$ 
of~$E(N,Nd^2+d)$ 
(smoothness requires that $N$ and $Nd^2+d$ are coprime, which is equivalent to our assumption that $N$ and $d$ are coprime). 
This polarisation can therefore be seen as a Liouville polarisation $(\Sigma',\lambda')$ of degree~$m Nd$ 
of $E(\frac 1d, d+\frac 1N)$. 
The polarizing curve~$\Sigma'$ has area 
$$
\ca_\om (\Sigma') \,=\, m(Nd+1) \,>\, \ca_\om (\Sigma_{m Nd}) ,
$$
and by the Riemann--Hurwitz formula its genus is 
$$
g(\Sigma') \,=\, \frac 12 \bigl( (m N-1)(m Nd^2+m d-1)-m+1 \bigr)
$$
which is $\geq g \bigl( \Sigma_{m Nd} \bigr) + b \bigl( \Sigma_{m Nd} \bigr) - 1$
by our assumption on $(m, d,N)$. 
By Lemma \ref{l:funny}, there exists an 
$(\alpha_{\st} |_{\Sigma_{m N d}},\alpha_{\st}|_{\Sigma'})$-exact symplectic embedding 
of $\widehat \Sigma_{m Nd}$
into~$\Sigma'$, 
and by Theorem~\ref{t:embliouville} an $(\alpha_{\st},\alpha_{\st})$-exact symplectic embedding of 
$\cb(\Sigma_{m Nd},\lambda_{m Nd})= B^4(1)\priv \Delta_{m Nd}$ 
into $E(\frac 1d, d+\frac 1N)$. 
\proofend

\m \ni
{\bf Proof of Theorem~\ref{t:anything}}.
It is shown in \cite[$\S$\26.2]{Sch05} 
that for every $\gve>0$ there exists $s_0 \in \RR$ such that
for every $s \geq s_0$ the ellipsoid $\lambda_s E(\frac 1s,s)$ 
of the same volume as the ball $B^4(a- \frac \gve 2)$
symplectically embeds into $(M,\go)$.
The claim now follows by precomposing this embedding with a scaling of a suitable embedding from 
Proposition~\ref{p:embell+}.
\proofend

The above proof gives no upper bound for $k$.
However, if $(M,\omega)$ is a rational symplectic manifold, or an affine part $M \setminus \Sigma$ therein, 
then by~\cite{Op07} there exists an explicit volume filling ellipsoid in $M$ that can be chosen to lie 
in the complement of~$\Sigma$.
Together with Proposition~\ref{p:embell+} one obtains an upper bound for~$k$.

%%%%%%%%%%%%%%%%%%%%%%%%%%%%%%%%%%%%%%%%%%%%%%%%%%%%%%%%%%%%%%%%%%%%%%%%%%%%%%%%%%%%%%%%%%%%%%%%%%

\section{Lagrangian rigidity} \label{s:lagrig}

\ni
{\bf Proof of Theorem~\ref{t:lagrig}.}
Let $\iota \colon L \hra D(A) \times D(B)$ 
be a Lagrangian embedding of a closed surface.   
Assume that there exists a Hamiltonian diffeomorphism~$\phi$ of~$\RR^4$
such that 
$$
\phi(L) \subset \bigl( D(A) \times D(B) \bigr) \priv \bigl( \Gamma_{\!\leq a} \! \times \Gamma_{\!\leq b} \bigr) 
$$
Since $\phi(L)$ is compact and disjoint from  
$\bigl( \Gamma_{\!\leq a} \cup \pp D(A) \bigr) \times \bigr( \Gamma_{\!\leq b} \cup \pp D(B) \bigr)$,
we can modify $\Gamma_{\!\leq a}$ and $\Gamma_{\!\leq b}$ 
%near their vertices 
%FF 
near the vertices of their closure 
to regular grids
$\Gamma'_{\!\leq a}$ and $\Gamma'_{\!\leq b}$
that still divide $D(A)$ and~$D(B)$ into topological discs of area~$\leq a$ and~$\leq b$,
and such that still
$$
\phi(L) \subset \bigl( D(A) \times D(B) \bigr) \priv \bigl( \Gamma'_{\!\leq a} \times \Gamma'_{\!\leq b} \bigr) .
$$
By Theorem~\ref{t:main} there exists an exact symplectic embedding
$$
\psi \colon \bigl( D(A) \times D(B) \bigr) \priv \bigl( \Gamma'_{\!\leq a} \times \Gamma'_{\!\leq b} \bigr)
\,\to\, Z^4(a+b).
$$
Then $\psi \circ \phi \circ \iota \colon L \to L' \subset Z^4(a+b)$ is a Lagragian embedding.
Since $\phi$ and~$\psi$ are exact, the action classes
$\iota^*[\alpha_{\st}]$ and $(\psi \circ \phi \circ \iota)^*[\alpha_{\st}]$ in~$H^1(L;\R)$ coincide.
Since by Stokes' theorem the minimal area~\eqref{e:areamin} is also the minimal action
$$
A_{\min}(L) \,=\, 
\inf \left\{ \int_\gamma \ga_{\st} \mid [\gamma] \in \pi_1(L), \, \int_\gamma \ga_{\st} >0 \right\} ,
$$
we obtain that 
$$
A_{\min}(L') = A_{\min}(L) .
$$
And since $L' \subset Z^4(a+b)$, its displacement energy $e(L)$ in $\C^2$ is  
$<a+b$, so by Chekanov's result from~\cite{Ch98}, 
$A_{\min} (L') \leq e(L) < a+b$.
Altogether,  $A_{\min}(L) < a+b$, as we wished to prove.
\proofend

%%%%%%%%%%%%%%%%%%%%%%%%%%%%%%%%%%%%%%%%%%%%%%%%%%%%%%%%%%%%%%%%%%%%%%%%%%%%%%%%%%%%%%%%%%%%%%

\section{Legendrian barriers} \label{s:legbar}

Theorem \ref{t:legbarriers1} on the existence of short Reeb chords between Legendrian curves 
is obtained from the Lagrangian rigidity result of the previous paragraph via 
a construction of Mohnke~\cite{Mo01}, 
that associates to the Reeb trajectory of a Legendrian knot a Lagrangian torus. 
We first review this construction. 

\begin{lemma} \label{l:mohnke} 
Let $(M^3,\xi,\alpha)$ be a contact manifold, 
$\bigl( SM = M \times \R_{>0}, d(R \alpha) \bigr)$ its symplectisation, 
$\Lambda \subset M$ a Legendrian knot, and $X \subset M$ any subset. 
Assume that there is no Reeb chord of length~$\leq T$ from
$\Lambda$ to $\Lambda \cup X$. 

Then there exists a Lagrangian torus $L$ in 
$(M \priv X) \times (0,1] \subset SM$ 
with $\ca_{\min}(L,SM)=T$.
\end{lemma}

\proof
Consider the map
$$
\iota \colon \Lambda \times (0,1] \times [0,\infty) \,\to\, SM,
\quad
(p,\tau, t) \mapsto \left( \Phi_\ga^t(p),\tau \right) ,
$$
where $\Phi^t_\ga$ is the Reeb-flow on $(M,\alpha)$.
The restriction of $\iota$ to every band $\{p\} \times (0,1] \times [0,\infty)$,
$p \in \Lambda$,
is symplectic for the forms $d \tau \wedge dt$ and~$d(R\alpha)$ 
on the domain and the target, respectively.
Since there is no Reeb chord of length $\leq T$ from 
$\Lambda$ to $\Lambda \cup X$, we find $\gve >0$ such that
$\iota$ is an embedding of $\Lambda \times (0,1] \times [0,T+\gve)$ into
$$
V \,:=\, (M \priv X) \times (0,1] \,\subset\, SM .
$$ 
Take a closed disc $D_\gamma$ in $(0,1] \times [0,T+\eps]$ of area~$T$ and smooth oriented boundary~$\gamma$. 
Then $\iota(\Lambda\times \gamma)$ is an embedded Lagrangian torus in~$V$.  
The actions of the generators $[\Lambda]$ and~$[\gamma]$ of~$\pi_1(L)$
are $\int_\Lambda \ga = 0$ and $\int_\gamma \ga = T$, so $\ca_{\min}(L, SM) = T$.
Moreover, the symplectic disc~$\iota(D_\gamma)$ lies in~$V$, has boundary on~$L$, 
and has area~$T$. 
\proofend

Returning to the setting of Theorem \ref{t:legbarriers1},
we notice that when $S$ is the smooth boundary of a starshaped domain $U \subset \R^4$, 
with contact form $\lambda_S = \ga_{\st}|_S$, 
then the exact symplectomorphism
$$
\begin{array}{lcl}
\bigl( S \times \R_{>0}, R \2 \lambda_{S} \bigr) & \lra & 
\bigl( \R^4 \priv \{0\}, \ga_{\st} \bigr) \\ [0.2em]
(s,R) &\longmapsto & R \, s
\end{array}
$$
identifies $S\times (0,1]$ with $U \priv \{0\}$ and 
$X \times(0,1]$ with the part in $U \priv \{0\}$ of the cone over~$X$ 
centered at the origin of~$\R^4$. 

\m \ni
{\it Proof of Theorem \ref{t:legbarriers1}:} 
Let $S$ be the smooth boundary of a starshaped domain $U \subset C^4(1) \subset \R^4$. 
Arguing by contradiction, assume that $\Lambda$ is a Legendrian knot in~$S$ with no Reeb chord of length~$\leq \delta_1+\delta_2$ 
from $\Lambda$ to $\Lambda \cup \Lambda_\delta$, where 
$\Lambda_\delta = (\Gamma_{\! \delta_1} \! \times \Gamma_{\! \delta_2}) \cap S$.
Since $\Gamma_{\! \delta_1}$ and $\Gamma_{\! \delta_2}$ are radial, 
$\Gamma_{\! \delta_1} \! \times \Gamma_{\! \delta_2}$ lies in the cone over $\Lambda_\delta$.
Hence Lemma~\ref{l:mohnke} provides a Lagrangian torus
$$
L \,\subset\, U \setminus \left( \Gamma_{\! \delta_1} \! \times \Gamma_{\! \delta_2} \right)
\,\subset\, C^4(1) \setminus \left( \Gamma_{\! \delta_1} \! \times \Gamma_{\! \delta_2} \right)
$$ 
with $\ca_{\min}(L) = \delta_1+\delta_2$.
%and a symplectic disc in $\RR^4 \priv \left( \Gamma_{\! \delta_1} \! \times \Gamma_{\! \delta_2} \right)$ with boundary on~$L$ and area~$\ca_{\min}(L, SM)$.  
This is a contradiction to Theorem~\ref{t:lagrig}.
\proofend

\paragraph{\bf An illustration.} 
To get a feeling for the phenomenon of Legendrian barriers, we explictely work out
the case where $S$ is the round sphere $S^3(1) = \pp B^4(1)$ 
and $\delta_1=\delta_2 = \frac 1k$.
Recall that 
$$
\Delta_k \,=\, \bigcup_{0 \leq i,j \leq k-1} \xi^i \R_{\geq 0} \times \xi^j \R_{\geq 0}
$$
where $\xi$ is the $k$\2th root of unity $e^{2\pi i/ k}$.
Hence the Legendrian graph $\Lambda_k := \Delta_k \cap S^3(1)$ 
is the union of the $k^2$ Legendrian quarter-circles
$$
Q_{i,j} \,:=\, 
\left(  \xi^i \R_{\geq 0} \times \xi^j \R_{\geq 0} \right) \cap S^3(1).
$$
Note that $Q_{i,j}$ and $Q_{i',j'}$ are disjoint if $i\neq i'$ and $j \neq j'$,
and intersect at one end-point if either $i \neq i'$ or~$j \neq j'$.
We group these $k^2$ quarter-circles in $k$ sets
$$
Q_j \,:=\, \coprod_{i=0}^{k-1} Q_{i,i+j}.
$$
Since the Reeb flow on $S^3(1)$ is the Hopf flow
$$
\Phi_R^t(z_1,z_2) \,=\, e^{2\pi it} (z_1,z_2) \,=\, 
\left( e^{2\pi it} z_1, e^{2\pi it}z_2 \right),
$$
we have $\Phi_R^{\frac ik}(Q_{0,j}) = \xi^i \, Q_{0,j} = Q_{i,i+j}$, 
that is, $\Phi_R^{\frac 1k}$ cyclically acts on the components of~$Q_j$.

The full sweep out of $Q_{0,j}$ under the (backward) Reeb flow is the Lagrangian surface 
$$
L_j \,:=\, \bigcup_{t \in \RR} e^{-2\pi i t} \, Q_{0,j} \,=\,
\bigcup_{t \in [0,1]} e^{-2\pi i t} \, Q_{0,j} \,=\,
\bigcup_{t \in [0,\frac 1k]} e^{-2\pi i t} \, Q_j .
$$
Therefore, if $\pi \colon S^3 \to S^2$ denotes the Hopf fibration 
(whose fibers are the Reeb trajectories), then
$$
L_j \,=\, \pi^{-1}(h_j), \quad \text{where } 
h_j=\pi(Q_{0,j}) 
\,=\, \pi (Q_j). 
$$
A computation shows that $h_j$ is a half-great circle joining the north pole to the south pole of~$S^2$. 
Since the area of the reduced space $S^2$ is~$1$, 
it follows that the pairs $(h_j,h_{j+1})$ each bound an open disc~$D_j$ 
of area~$\frac 1k$. 
%
%
%Another straightforward computation shows that ...
%$(\pi^{-1}(D_j),\alpha_\st)\approx (D(\frac 1k)\times S^1,\alpha_0=\lambda_\st+d\theta)$.
%
%PP \footnote{Quel calcul? Formellement, ca ne peut pas être un calcul, comme l'identification
%de $D_j$ avec $D(\frac 1k)$ n'est pas explicite.
%On a, par la Hopf projection, un diffeomorphism $\pi^{-1}(D_j) \to D_j \times S^1$. 
%Le calcul consiste-t'il en vérifiant que le push forward de $\ga_{\st}$ 
%est de la forme $\lambda (x,y) \oplus d t$ ?
%Il suivrait alors, après un symplectomorphism $D_j \to D(\frac 1k)$, 
%que la forme sur $D(\frac 1k) \times S^1$ a la forme $(\ga_{\st} + dh) \oplus d t$.
%Il faut alors encore tuer le $dh$. Ca se fait en appliquant Moser à $\gl_t := \gl_{\st} \oplus t dh$; il n'a pas de problèmes avec le bord, comme le Moser vector field est parallel au Reeb field $\pp_t$.} 
%
It is also not hard to see that 
there exists a diffeomorphism $\pi^{-1}(D_j) \to D(\frac 1k)\times S^1$
that takes $\lambda_{S^3(1)}$ to $\ga_{\st} \oplus dt$.
Altogether, we see that 
$$
L \,:=\, \bigcup_j L_j \,=\, 
\bigcup_{t \in [0,\frac 1k]} \Phi^{-t}_R (\Lambda_k) 
$$
cuts $S^3(1)$ into $k$ connected components, each of which can be identified with 
$\bigl( D(\frac 1k) \times S^1, \ga_\st + dt \bigr)$, where $S^1 = \RR/\ZZ$. 

Now take a Legendrian knot $\Lambda$ in~$S^3(1)$ such that there is no Reeb chord 
of length~$\leq \frac 1k$ from~$\Lambda$ to~$\Lambda_k$. 
Equivalently, $\Lambda$ is disjoint from 
$\bigcup_{t \in [0,\frac 1k]} \Phi^{-t}_R (\Lambda_k)$,
which is~$L$. 
So $\Lambda$ lies in one of the connected components of~$S^3(1) \priv L$
and can be seen as a Legendrian knot in 
$\bigl( D(\frac 1k) \times S^1, \ga_\st + dt \bigr)$.
But a very classical and elementary argument 
(see e.g.\ \cite[p.\ 192]{AS19}) 
shows that such a knot has a Reeb chord of length~$\leq \frac 1k$
(in fact, of length $< \frac 1{2k}$). 

Summing up this discussion, in the case of the round sphere $S = S^3(1)$, 
where the Reeb flow is explicit,   
the barrier property of~$\Lambda_k$ directly follows from the fact that the 
$\frac 1k$-negative Reeb chords starting at~$\Lambda_k$ disconnect~$S$ into 
pieces all of whose Legendrian knots have small Reeb chords. 
There is no reason, however, that this disconnectedness property remain true 
when the Reeb flow is modified by taking the contact form on an arbitrary starshaped domain in~$C^4(1)$. 
It is therefore remarkable that the Legendrian rigidity result holds true.

%%%%%%%%%%%%%%%%%%%%%%%%%%%%%%%%%%%%%%%%%%%%%%%%%%%%%%%%%%%%%%%%%%%%%%%%%%%%%%%%%%%%%%%

\end{document}